\pdfoutput=1
\documentclass[a4paper, 11pt, abstracton, titlepage, oneside]{scrartcl}

\makeatletter
\usepackage{pgfplots}
\usepgfplotslibrary{groupplots} % For creating subfigures

\usepackage[T1]{fontenc}
\usepackage[utf8]{inputenc}
\usepackage[onehalfspacing]{setspace}
\usepackage[headsepline]{scrlayer-scrpage}
\clearscrheadfoot % reset standard format
\usepackage{layout}
\usepackage{geometry}
\usepackage{adjustbox}
\usepackage{authblk}
\usepackage[titletoc,title]{appendix}

\usepackage[hyphens]{url}
\usepackage{color}

\usepackage{amsmath,amssymb,amsfonts,amsthm, mathtools, bm}
\usepackage{thmtools, thm-restate}
\usepackage{nicefrac}
\usepackage{array}

\usepackage[normal]{threeparttable}
\usepackage{threeparttablex}
\usepackage{tabularx}
\usepackage{longtable}
\usepackage{booktabs}
\usepackage{colortbl}
\usepackage{multirow}
\usepackage{makecell}
\usepackage{pdflscape}

\usepackage{standalone}

\usepackage[acronym]{glossaries}

\usepackage{subcaption}
\usepackage{caption}
\usepackage{floatrow}
\usepackage{graphicx}
\usepackage{placeins}

\usepackage[ruled,vlined,linesnumbered,resetcount]{algorithm2e}
\SetKw{Continue}{continue}
\SetKwComment{Comment}{/* }{ */}
\SetKwInput{KwInput}{Input}
\SetKwInOut{KwOutput}{return}
\usepackage{float}
\usepackage{algpseudocode}

\usepackage{tikz}
\usetikzlibrary{arrows.meta}
\usetikzlibrary{math}
\usetikzlibrary{shapes}
\usepackage{pgfplots}
\usepgfplotslibrary{groupplots}
\usepgfplotslibrary{dateplot}
\usepackage{tikzscale}

\usepackage{enumerate}
\usepackage{multicol}
\usepackage[]{pdfcomment}
\usepackage{xparse}
\usepackage{twoopt}
\usepackage{etoolbox}

\usepackage{eurosym}
\usepackage{ellipsis}
\usepackage{natbib}
\usepackage{paralist}
\usepackage{ulem}

\graphicspath{{Figures/}}
\usepackage{wrapfig}

\usepackage{optidef}

\usepackage{bm}
\usepackage{bbm}

\usepackage{pgfplots}
\usepgfplotslibrary{groupplots} % For creating subfigures
\usepackage{import}
\pgfplotsset{compat=1.17}
\usepgfplotslibrary{statistics}

\usepackage{subcaption}

\usepackage{svg}
\usepackage{enumitem}

\usepackage{afterpage} % Add to your preamble if not already included

\AtBeginDocument{
	\newgeometry{
		textheight=9in,
		textwidth=6.5in,
		top=1in,
		headheight=14pt,
		headsep=25pt,
		footskip=30pt
	}
}
\newsavebox{\abstractbox}
\renewenvironment{abstract}
{\begin{lrbox}{0}\begin{minipage}{\textwidth}
			\begin{center}\normalfont\sectfont\abstractname\end{center}\quotation}
		{\endquotation\end{minipage}\end{lrbox}%
	\global\setbox\abstractbox=\box0 }

\makeatletter
\expandafter\patchcmd\csname\string\maketitle\endcsname
{\vskip\z@\@plus3fill}
{\vskip\z@\@plus2fill\box\abstractbox\vskip\z@\@plus1fill}
{}{}
\makeatother

\addtokomafont{captionlabel}{\footnotesize\bfseries} % style for caption labels
\addtokomafont{caption}{\footnotesize\bfseries} % and the caption
\bibpunct[, ]{(}{)}{,}{a}{}{,}%
\counterwithin*{equation}{section}

\newenvironment{myfont}{\fontfamily{ccr}\selectfont}{\par}
\DeclareTextFontCommand{\textmyfont}{\myfont}

\AtBeginDocument{%
	\abovedisplayskip=7.5pt plus 0pt minus 0pt
	\abovedisplayshortskip=7.5pt plus 0pt minus 0pt
	\belowdisplayskip=7.5pt plus 0pt minus 0pt
	\belowdisplayshortskip=7.5pt plus 0pt minus 0pt
}
\newcolumntype{L}[1]{>{\raggedright\let\newline\\\arraybackslash\hspace{0pt}}p{#1}}
\newcolumntype{C}[1]{>{\centering\let\newline\\\arraybackslash\hspace{0pt}}p{#1}}
\newcolumntype{R}[1]{>{\raggedleft\let\newline\\\arraybackslash\hspace{0pt}}p{#1}}
\def\checkmark{\tikz\fill[scale=0.4](0,.35) -- (.25,0) -- (1,.7) -- (.25,.15) -- cycle;}

\renewcommand{\emph}[1]{\textit{#1}}

\begin{document}
    
\setlength{\belowdisplayskip}{0pt} \setlength{\belowdisplayshortskip}{0pt}
\setlength{\abovedisplayskip}{0pt} \setlength{\abovedisplayshortskip}{0pt}

\emergencystretch 3em
% acronyms and abbrevs.
\newacronym{fe}{FE}{first echelon}
\newacronym{se}{SE}{second echelon}
\newacronym{2e-lrpmd-ds}{2E-LRPMD-DS}{two-echelon location routing problem with mobile depots and direct shipment}
\newacronym{2e-lrp-ds}{2E-LRP-DS}{two-echelon location routing problem with direct shipment}
\newacronym{alns}{ALNS}{adaptive large neighborhood search} 
\newacronym{flp}{FLP}{facility location problem} 
\newacronym{vrp}{VRP}{vehicle routing problem} 
\newacronym{cflp}{CFLP}{capacitated facility location problem}
\newacronym{2e-lrp}{2E-LRP}{two-echelon location routing problem}
\newacronym{lrp}{LRP}{location routing problem}
\newacronym{milp}{MILP}{mixed integer linear programming}
\newacronym{2e-vrp}{2E-VRP}{two-echelon vehicle routing problem}
\newacronym{2e-clrp}{2E-CLRP}{two-echelon capacitated location routing problem}
\newacronym{vrp}{VRP}{vehicle routing problem}
\newacronym{clrp}{CLRP}{capacitated location routing problem}
\newacronym{vrpis}{VRPIS}{vehicle routing problem with intermediate stops}
\newacronym{tsp}{TSP}{traveling salesman problem}
\newacronym{ip}{IP}{integer program}
\newacronym{grasp}{GRASP}{greedy randomized adaptive search process}
\newacronym{ga}{GA}{genetic algorithm}
\newacronym{ils}{ILS}{iterated local search}
\newacronym{mip}{MIP}{mixed integer programming}
\newacronym{ts}{TS}{tabu search}
\newacronym{bnc}{B\&C}{branch-and-cut}
\newacronym{lrp-mppd-2e}{LRP-MPPD-2E}{two-echelon multi-products location-routing problem with pickup and delivery}
\newacronym{2e-lrptwtrs}{2E-LRPTWTRS}{two-echelon location-routing problem with time windows and transportation resource sharing}
\newacronym{mp-2elrp}{MP-2ELRP}{multi-period two-echelon location routing problem}
\newacronym{2ecd-ms}{2ECD-MS-TDD}{two-echelon dispatching
model with mobile satellites}
\newacronym{vns}{VNS}{variable neighborhood search}
\newacronym{lns}{LNS}{large neighborhood search}
\newacronym{vrptwmd}{VRPTWMD}{vehicle routing problem with time windows and mobile depots}
\newacronym{2e-mtvrpds}{2E-MTVRPDS}{two-echelon multi-trip vehicle routing problem with a dynamic satellite}
\newacronym{lip}{LIP}{location-inventory model}
\newacronym{lirp}{LIRP}{location-inventory-routing model}
\newacronym{fe}{FE}{first echelon}
\newacronym{lsp}{LSP}{Logistic service provider}
\newacronym{sc}{SC}{set cover}
\newacronym{scp}{SCP}{set cover problem}
\newacronym{sa}{SA}{simulated annealing}
\newacronym{2e-cls}{2E-CLS}{two-echelon city logistics systems}
\newacronym{2e-lrpmd}{2E-LRPMD}{two-echelon location routing problem with mobile depots}
\newacronym{pv}{PV}{primary vehicle}
\newacronym{sv}{SV}{secondary vehicle}
\newacronym{mdvrp}{MDVRP}{multi-depot vehicle routing problem}
\newacronym{sttrpsd}{STTRPSD}{single truck and trailer routing problem with satellite depots}
\newacronym{pf}{PF}{progressive filtering}
\newacronym{lrpif}{LRPIF}{location-routing problem with intra-route facilities}
\newacronym{osm}{OSM}{openstreetmap}
\newacronym{cw}{CW}{Clarke-Wright}
\newacronym{kpi}{KPI}{key performance indicator}

% alns parameters
\newacronym{maxIter}{$\eta^{max}$}{the maximum number of iterations}
\newacronym{iter}{$\iota$}{iteration}
\newacronym{lambda}{$\lambda$}{the decay parameter}
\newacronym{reheat_iter}{$\eta^{reheat}$}{the number of iterations without improvement}
\newacronym{weight}{$\omega^{1},\omega^{2},\omega^{3},\omega^{4}$}{the operator weight update parameters}
\newacronym{neigh_count}{$\Pi$}{the number of neighbors for each node}
\newacronym{t_init}{$\tau_{0}$}{the initial temperature}
\newacronym{dod}{$\Omega$}{the number of customers to be removed in destroy operators}
\newacronym{t_final}{$\tau_{final}$}{the final temperature}
\newacronym{ls_counter}{$\upsilon$}{the local search improvement method}
\newacronym{temp_alpha}{$\alpha$}{the temperature decay value}
\newacronym{sc_iter}{$\eta^{SC}$}{the number of iterations after which the set cover problem applied}
\newacronym{currentTemp}{$\tau$}{current temperature}
\newacronym{temp_update}{$g(\tau)$}{the temperature update function}

%%%%%%%%%%%%%%%%%%%%%%%%%%%%%%%%%%%%%%%%%%%%%%%%%%%%%%%%%%%%%%%%%%%%%%%%%%%%%%%%%%%%%%%%%%%%%%%%%%%%%%%%%%%%%%%%%%%%%%%%%%
% Commands

\newcommand{\mNodeFirst}[0]{\ensuremath{i}}
\newcommand{\mNodeSecond}[0]{\ensuremath{j}}
\newcommand{\mCustomer}[0]{\ensuremath{j}}
\newcommand{\mSatellite}[0]{\ensuremath{t}}
\newcommand{\mVehicleFirst}[0]{\ensuremath{k}} %it was t before, I changed it because I now use t for satellites 
\newcommand{\mVehicleSecond}[0]{\ensuremath{k}}
\newcommand{\mArc}[2]{\ensuremath{a^{#1}_{#2}}}
\newcommand{\mTypeHub}[0]{\ensuremath{m}}
\newcommand{\mTypeCargo}[0]{\ensuremath{h}}
\newcommand{\mTrip}[0]{\ensuremath{w}}
\newcommand{\mDepot}[0]{\ensuremath{m}}
\newcommand{\sDepot}[0]{\ensuremath{s}}

% sets
\newcommand{\mGraph}[0]{\ensuremath{\mathcal{G}}}
\newcommand{\mSetVertices}[0]{\ensuremath{\mathcal{V}}}
\newcommand{\mSetCustomers}[0]{\ensuremath{\mathcal{C}}}
\newcommand{\mSetHubLocations}[0]{\ensuremath{\mathcal{D}}}
\newcommand{\mSetDepots}[1]{\ensuremath{\mathcal{D}^{\mathrm{#1}}}}
\newcommand{\mSetMiniNetwork}[0]{\ensuremath{\mathfrak{N}}}
\newcommand{\mGraphArc}[0]{\ensuremath{\mathcal{A}}}
\newcommand{\mGraphVertex}[0]{\ensuremath{\mathcal{V}}}
\newcommand{\mSetArcs}[1]{\ensuremath{\mathcal{A}}^{#1}}

\newcommand{\mSetTypeVehicle}[1]{\ensuremath{\mathcal{K}^{#1}}}
\newcommand{\mTypeVehicle}[0]{\ensuremath{k}}
\newcommand{\mSetTypeHubs}[0]{\ensuremath{\mathcal{M}}}
\newcommand{\mSetTypeCargo}[0]{\ensuremath{\mathcal{H}}}
\newcommand{\mSetTrips}[0]{\ensuremath{\mathcal{W}}}
\newcommand{\mSetTypeVehicleECVOnly}[0]{\ensuremath{\mathcal{K}^{\mathtt{ECV}}}}

\newcommand{\mSetArcsVehicle}[2]{\ensuremath{\mathcal{A}^{#1}_{#2}}}
\newcommand{\mSetArcsTrucks}[0]{\ensuremath{\mSetArcsVehicle{\text{t}}}}
\newcommand{\mSetArcsFreighters}[0]{\ensuremath{\mSetArcsVehicle{\text{f}}}}

%parameters
\newcommand{\mCostLocation}[2]{\ensuremath{c^{#1}_{#2}}}
\newcommand{\mCostArc}[2]{\ensuremath{c^{#1}_{#2}}}
\newcommand{\mFixedCost}[1]{\ensuremath{f_{#1}}}

\newcommand{\mDemand}[2]{\ensuremath{d^{#1}_{#2}}}
\newcommand{\mEarliestBOS}[1]{\ensuremath{e_{#1}}}
\newcommand{\mLatestBOS}[1]{\ensuremath{l_{#1}}}
\newcommand{\mServiceTime}[2]{\ensuremath{s^{#1}_{#2}}}
\newcommand{\mHorizonStart}[0]{\ensuremath{E}}
\newcommand{\mHorizonEnd}[0]{\ensuremath{L}}

\newcommand{\mDecisionHubLocation}[2]{\ensuremath{y^{#1}_{#2}}}
\newcommand{\mDecisionArc}[2]{\ensuremath{x^{#1}_{#2}}}
\newcommand{\mDecisionVehicle}[2]{\ensuremath{z^{#1}_{#2}}}
\newcommand{\mWaitingTime}[2]{\ensuremath{r^{#1}_{#2}}}
\newcommand{\mArrival}[2]{\ensuremath{v^{#1}_{#2}}}
\newcommand{\mDecisionArrival}[2]{\ensuremath{z^{#1}_{#2}}}
\newcommand{\mSyncEchelons}[2]{\ensuremath{u^{#1}_{#2}}}

\newcommand{\mIncomingArcs}[1]{\ensuremath{\delta^{+}(k,#1)}}
\newcommand{\mOutgoingArcs}[1]{\ensuremath{\delta^{-}(k,#1)}}

\newcommand{\mOriginOfArc}[1]{\ensuremath{o(#1)}}
\newcommand{\mDestinationOfArc}[1]{\ensuremath{d(#1)}}

\newcommand{\mFreightCapacity}[2]{\ensuremath{Q^{#1#2}}}
\newcommand{\mHubCapacity}[2]{\ensuremath{Q^{#1#2}}}
\newcommand{\mResidualCargo}[3]{\ensuremath{q^{#1#2}_{#3}}}

\newcommand{\mCargoCapacity}[2]{\ensuremath{Q^{#1}_{#2}}}
\newcommand{\mCargoDemand}[2]{\ensuremath{d^{#1}_{#2}}}

\newcommand{\mBeginOfService}[2]{\ensuremath{\alpha^{#1}_{#2}}}
\newcommand{\mDeparture}[2]{\ensuremath{\beta^{#1}_{#2}}}
\newcommand{\mTravelTimeOfArc}[3]{\ensuremath{\tau^{#1}_{#2#3}}}

\newcommand{\mBatteryCapacity}[1]{\ensuremath{G^{#1}}}
\newcommand{\mResidualBattery}[2]{\ensuremath{g^{#1}_{#2}}}
\newcommand{\mLoadFlow}[2]{\ensuremath{f^{#1}_{#2}}}
\newcommand{\mBatteryConsumptionOfArc}[2]{\ensuremath{\gamma^{#1}_{#2}}}

\newcommand{\mHubLocationAvailable}[2]{\ensuremath{a^{#1}_{#2}}}
\newcommand{\mCargoVehicleCompatibility}[2]{\ensuremath{\eta ^{#1}_{#2}}}
\newcommand{\mCargoSatelliteCompatibility}[2]{\ensuremath{\theta ^{#1}_{#2}}}
\newcommand{\mHubTransferTime}[1]{\ensuremath{t^{m}}}
\newcommand{\mHubCargoCapacity}[2]{\ensuremath{Q^{#1#2}}} % cargo hold + supported cargo type

%alns values and parameters
\newcommand{\currentSol}[0]{\ensuremath{s}}
\newcommand{\newSol}[0]{\ensuremath{\bar{s}}}
\newcommand{\bestSol}[0]{\ensuremath{s^{*}}}
\newcommand{\scSol}[0]{\ensuremath{\hat{s}}}
\newcommand{\repairSol}[0]{\ensuremath{\hat{s}}}
\newcommand{\repairRoute}[0]{\ensuremath{\hat{r}}}
\newcommand{\openDepots}[0]{\ensuremath{\bar{\mathcal{D}}}}
\newcommand{\removedCustomers}[0]{\ensuremath{\mathcal{C}^{R}}}
\newcommand{\mSetNetwork}[0]{\mathcal{N'}}
\newcommand{\route}[2]{\ensuremath{r^{#1}_{#2}}}
\newcommand{\node}[1]{\ensuremath{v_{#1}}}
\newcommand{\md}[0]{\ensuremath{m}}
\newcommand{\customer}[0]{\ensuremath{i}}
\newcommand{\customerj}[0]{\ensuremath{j}}

\newcommand{\neighborCount}[0]{\ensuremath{\Pi}}

\newcommand{\abs}[1]{\left\lvert #1 \right\rvert}

%decompositions
\newcommand{\baseflp}[0]{\ensuremath{FLP^{B}}}
\newcommand{\basekmeans}[0]{\ensuremath{kMeans^{B}}}
\newcommand{\decompflp}[0]{\ensuremath{FLP^{D}}}
\newcommand{\decompkmeans}[0]{\ensuremath{kMeans^{D}}}

% \newcommand{\myalgorithm}{%
% \footnotesize
% \begin{algorithm}[H]
%     % \caption{\footnotesize Cluster-first-route-second based Decomposition Approach}
%     \label{alg:decomposition}
%     \KwIn{Instance $I$}
%     \SetKwFunction{generateClusters}{GenerateFirstLevelCluster}
%     \SetKwFunction{kmeans}{k-means}
%     \SetKwFunction{tsp}{SolveTSP}
%     \SetKwFunction{solve}{SolveCFLP}
%     \SetKwFunction{combine}{Aggregate}
%     \SetKwFunction{alns}{ALNS}
        
%     \BlankLine
%     $\mathbb{C} \gets \generateClusters(I)$ \label{generate} \\
%     \ForEach{$\mathbb{C}_{i} \in \mathbb{C}$}{
%         \BlankLine
%         \kmeans($\mathbb{C}_{i}$) \label{second_level_cluster}\\
%         \ForEach{$C_{ij} \in \{C_{i1}, C_{i2}, \ldots, C_{iN}\}$}{
%             \BlankLine
%             \tsp{$\mathbb{C}_{ij}$} \label{solve_tsp}
%         }
%     }
%     \BlankLine
%     % \textbf{Step 4: Aggregate Solutions}\;
%     $s^{agg} \gets \combine(C_{ij})$ \label{aggregate}
    
%     \BlankLine
%     % \textbf{Step 5: Apply \gls{alns}}\;
%     $\bestSol\gets \alns({s^{agg}})$ \label{decomp_alns}
    
%     \BlankLine
%     \Return{\bestSol}
%     \end{algorithm}}

% then the title
%----------------------------------------------------------------------------------
% titlepage template
% author: Maximilian Schiffer
% version: 1.0 / 11.05.2017
%----------------------------------------------------------------------------------

\title{\large Solving Large-Scale Two-Echelon Location Routing Problems in City Logistics}

% and the authors
\author[1]{\normalsize Banu Ulusoy Dereli}
\author[1,2]{\normalsize Gerhard Hiermann}
\author[1,3]{\normalsize Maximilian Schiffer}
\affil{\small 
	TUM School of Management, Technical University of Munich, Germany
	
	\scriptsize banu.dereli@tum.de\\	

        \small 
	\textsuperscript{2}Institute of Production and Logistics Management, Johannes Kepler University Linz, Austria
	
	\scriptsize gerhard.hiermann@jku.at
	
	\small
	\textsuperscript{3}Munich Data Science Institute, Technical University of Munich, Germany
	
	\scriptsize schiffer@tum.de}

% if you like - a date
\date{}

% in case you have a headline - otherwise outcomment
\lehead{\pagemark}
\rohead{\pagemark}

% finally your abstract
\begin{abstract}
\begin{singlespace}
{\small\noindent Logistic service providers increasingly focus on two-echelon distribution systems to efficiently manage thousands of deliveries in urban environments. Effectively operating such systems requires designing cost-efficient delivery networks while addressing the challenges of increasing e-commerce demands. In this context, we focus on a \glsentrylong{2e-lrpmd-ds}, where decisions involve locating micro-depots, and designing first and second-level routes. Our model also incorporates the flexibility of direct shipments from the main depot to customers.

To solve such large-scale problems efficiently, we propose a metaheuristic approach that integrates a set cover problem with an adaptive large neighborhood search (ALNS). Our ALNS approach generates a set of promising routes and micro-depot locations using destroy and repair operators while using a local search for intensification. We then utilize the set cover problem to find better network configurations. Additionally, we present a decomposition-based cluster-first, route-second approach to solve large-scale instances efficiently. We show the efficacy of our algorithm on well-known benchmark datasets and provide managerial insights based on a case study for the city of Munich. Our decomposition approach provides comparable results while reducing computational times by a factor of 15. Our case study results show that allowing direct shipment can reduce total costs by 4.7\% and emissions by 11\%, while increasing truck utilizations by 42\%. We find that integrating both stationary and mobile micro-depots, along with allowing direct shipments, can reduce total costs by 5.9\% compared to traditional two-echelon delivery structures.\\
\smallskip}


\end{singlespace}
\end{abstract}

% don't forget to make the tile
\maketitle

% and your chapters
\section{Introduction}
Cities are experiencing increases in urbanization at a historically exceptional rate. In 2030, 60\% of the world's total population is expected to live in cities. As urbanization increases, congestion and pollution grow exponentially~\citep{wef2020}\unskip.~By 2023, global retail e-commerce has reached 5.8 trillion US dollars, and it is expected to increase by 39\% until 2027~\citep{Statista2024}. This leads to cities encountering challenges in handling large volumes of parcel deliveries in city logistics.~In this context, it is essential to establish good distribution networks to efficiently and effectively transport goods in urban areas.

\glspl{lsp} traditionally manage urban deliveries through a central depot located on the outskirts of a city, distributing goods directly to customers. However, such a network structure bears several inefficiencies, particularly with increasing delivery volumes, urban traffic congestion, and last-mile delivery challenges. To tackle these issues, \glspl{lsp} increasingly locate smaller micro-depots within city centers, which allows for significantly reduced costs compared to direct deliveries from a single main depot. In response to these challenges, \glspl{lsp} are shifting their focus to \gls{2e-cls}~\citep{Schiffer2022}.

In two-echelon networks, the main depot operates outside the city, while multiple micro-depots bring distribution points closer to customers. Large vehicles transport parcels from the main depot to micro-depots, and smaller vehicles handle last-mile deliveries. By operating \gls{2e-cls}, \glspl{lsp} can also reduce environmental impact by limiting the use of large trucks in city centers and favoring sustainable delivery vehicles, e.g., cargo bikes. Studies highlight the advantages of cargo bikes in \gls{2e-cls}, particularly in lowering transportation costs and $CO_{2}$ emissions~\citep{Sheth2019, Fontaine2021}. 

In a two-echelon distribution setting, most studies assume that the customers can only be served through the satellites, such as micro-depots, rather than directly by the main depot. However, depending on the customer and micro-depot locations, delivering via cargo bikes from the micro-depots may result in higher routing costs. Moreover, the cargo bikes might not be able to transport the parcels due to capacity or safety reasons.~To mitigate these problems, direct transportation can increase flexibility in the distribution network, yielding total cost reductions and increased customer satisfaction with on-time delivery.  

Due to daily changes in delivery points~\citep{Febransyah2022}, designing an efficient two-echelon system can be challenging. In this context, \glspl{lsp} face a trade-off: increasing the number of stationary micro-depots can lower routing costs but raise overall expenses due to high setup costs in city centers~\citep{Crainic2010}. To this end, mobile micro-depots, such as swap body containers or mobile parcel lockers, provide a cost-effective and flexible alternative. These micro-depots serve as temporary inner-city distribution hubs, reducing setup costs and allowing easy adaptation to changing delivery requirements. Accordingly, there is a need to explore the benefit of integrating mobile micro-depots into two-echelon networks.

Against this background, we focus on solving the \gls{2e-lrp} with mobile depots and direct shipment in which we consider locating both stationary and mobile micro-depots while additionally allowing for direct shipment from the main depot to customers. We propose a metaheuristic that integrates a \glsentrylong{sc} problem with an \gls{alns} algorithm and develop a decomposition-based cluster-first-route-second approach to solve large-scale instances fast and efficiently.

\subsection{State of the Art} \label{lit-rev}
 Our work extends the \gls{2e-lrp} by conjointly considering locating stationary and mobile micro-depots, further allowing for direct shipments from the main depot. In the following, we concisely review related literature, mostly focusing on a city logistics context. 
 
 To the best of our knowledge, the introduction of the \gls{2e-lrp} dates back to \cite{Madsen1980} and has gained significant interest in the research community since then. \cite{Boccia2010} proposed a \glsentrylong{ts} approach that decomposes the problem into a \gls{cflp} and a \glsentrylong{mdvrp}. \cite{Nguyen2012} developed an \glsentrylong{ip} to formulate the \gls{2e-lrp} with a single depot and proposed a \glsentrylong{grasp}, as well as a multi-start \glsentrylong{ils} algorithm \citep{Nguyen2010}. \cite{Contardo2012} proposed a \glsentrylong{bnc} algorithm for solving the \gls{2e-lrp} based on a two-index vehicle flow formulation to solve small and medium-sized instances optimally and complement it by developing an ALNS algorithm. \cite{Rahmani2016} developed clustering-based approaches to solve the \gls{2e-lrp} with pickup and delivery with multi-products. \cite{Mirhedayatian2021} focused on the \gls{2e-lrp} with synchronization of both echelons and proposed a decomposition-based heuristic. \cite{Sörensen2021} introduced a heuristic framework to solve the \gls{2e-lrp}'s routing subproblem, and combined it with \glsentrylong{pf} to remove unpromising depot configurations.

Our problem setting combines a \gls{2e-lrp} with a \glsentrylong{vrpis} (\cite{Schiffer2019}), in which we determine the location decisions of both stationary and mobile micro-depots. In the context of city logistics, most studies only consider locating stationary micro-depots on a strategic level. Only recently have some works considered the location of mobile depots in a two-echelon distribution setting to properly demonstrate daily fluctuations in customer demand and locations. \cite{Lan2022} proposed a \glsentrylong{2ecd-ms} in which the locations of the mobile depots change according to customers' demands and trucks directly dispatch the customers. \cite{SchifferWalther2018} proposed a \glsentrylong{lrpif}, which can be charging stations as well as pick-up or unloading stations for freight or waste. \cite{Hof2021} proposed an \gls{alns} and a path relinking approach to solve an intra-route resource replenishment problem with mobile depots. \cite{Sutrisno2023} studied a \gls{2e-lrp} with mobile satellites and proposed a clustering-based simultaneous neighborhood search. \cite{Tian2023} studied a \gls{2e-lrp} in the context of city logistics in which they recommended satellite locations from the set of customers. 

The second main component of our problem is a direct shipment from the main depot to customers. In their review paper, \cite{Sluijk2023} highlighted the potential impact of direct deliveries by \glsentrylong{fe} vehicles on solution costs. \cite{Guastaroba2016} introduced hybrid networks where freight can either pass through an intermediate facility or be delivered directly from the main depot to customers. \cite{Anderluh2017} classified customers into bike and van-customers and addressed the \gls{2e-vrp} by synchronizing vans and cargo bikes at satellite depots. Similarly, \cite{Lan2022} allowed goods to be dispatched directly from a central depot to customers, further exploring the role of direct deliveries. \cite{Song2023} examined a \glsentrylong{lip} and a \glsentrylong{lirp}, comparing supply chain structures with and without direct shipments from suppliers to retailers. Additionally, \cite{Mokhtarinejad2015} modeled an integrated vehicle routing and scheduling problem for cross-docking systems, incorporating direct shipments from manufacturers to customers. They proposed a machine learning-based heuristic method and a \gls{ga} approach to address large-scale instances.

Finally, we aim to solve large-scale problems. \cite{Madsen1980} worked on an instance size with up to 4500 customers, where they developed three problem-specific heuristics. However, this work did not consider capacity constraints and establishment costs for satellites. In \cite{Drexl2015} and \cite{Cuda2015}, the number of customers and satellites considered in the \gls{2e-lrp} setting were up to 200 and 20, respectively. For the \gls{lrp}, \cite{SchneiderLöffler2017} generated instances with up to 600 customers and 30 depot locations. Accordingly, studies in the literature that solve large-scale instances with thousands of customers in a \gls{2e-lrp} setting are missing so far. 

Table \ref{table:Table1} shows the most related works that capture the specific components of the \gls{2e-lrp}. All of these existing works consider either locating mobile or stationary micro-depots, but not both. Only three of them allow for direct shipments in their problem setting. These existing works so far solve up to 200 customers but have not been applied to large instances. As can be seen, our work is the first to study the \gls{2e-lrp} with both mobile and stationary micro-depots, as well as allowing for direct shipments from the main depot to customers, proposing an algorithm that solves large-scale problems. 

\begin{table}[htbp]
  \centering
    \caption{Related works on \gls{2e-lrp} and \gls{lrp}.}
    \small
    \resizebox{0.9\textwidth}{!}{
    \begin{tabular}{lp{0.2cm}p{0.2cm}p{0.2cm}p{0.2cm}p{0.2cm}p{0.2cm}p{0.2cm}p{0.2cm}p{0.2cm}p{0.2cm}p{0.2cm}}
    \toprule    
    & \rotatebox{60}{\cite{Mirhedayatian2021}}
    & \rotatebox{60}{\cite{Cheng2022}}
    & \rotatebox{60}{\cite{Lan2022}}
    & \rotatebox{60}{\cite{Hof2021}}
    & \rotatebox{60}{\cite{Nguyen2010}}
    & \rotatebox{60}{\cite{He2019}}
    & \rotatebox{60}{\cite{Wang2021}}
    & \rotatebox{60}{\cite{Rahmani2015}}
    & \rotatebox{60}{\cite{Voigt2022}}
    & \rotatebox{60}{\cite{Sutrisno2023}}
    & \rotatebox{60}{\textbf{Our Work}} \\
    \midrule
    Mobile micro-depot & \multicolumn{1}{c}{-} & \multicolumn{1}{c}{-} & \multicolumn{1}{c}{\checkmark} & \multicolumn{1}{c}{\checkmark} & \multicolumn{1}{c}{-} & \multicolumn{1}{c}{\checkmark} & \multicolumn{1}{c}{-} & \multicolumn{1}{c}{-} & \multicolumn{1}{c}{-} & \multicolumn{1}{c}{\checkmark} & \multicolumn{1}{c}{\checkmark} \\
    Stationary micro-depot & \multicolumn{1}{c}{\checkmark} & \multicolumn{1}{c}{\checkmark} & \multicolumn{1}{c}{-} & \multicolumn{1}{c}{-} & \multicolumn{1}{c}{\checkmark} & \multicolumn{1}{c}{-} & \multicolumn{1}{c}{\checkmark} & \multicolumn{1}{c}{\checkmark} & \multicolumn{1}{c}{\checkmark} & \multicolumn{1}{c}{-} & \multicolumn{1}{c}{\checkmark} \\
    Location decision & \multicolumn{1}{c}{\checkmark} & \multicolumn{1}{c}{\checkmark} & \multicolumn{1}{c}{\checkmark} & \multicolumn{1}{c}{\checkmark} & \multicolumn{1}{c}{\checkmark} & \multicolumn{1}{c}{\checkmark} & \multicolumn{1}{c}{\checkmark} & \multicolumn{1}{c}{\checkmark} & \multicolumn{1}{c}{\checkmark} & \multicolumn{1}{c}{\checkmark} & \multicolumn{1}{c}{\checkmark} \\
    Customer time window & \multicolumn{1}{c}{\checkmark} & \multicolumn{1}{c}{\checkmark} & \multicolumn{1}{c}{-} & \multicolumn{1}{c}{\checkmark} & \multicolumn{1}{c}{-} & \multicolumn{1}{c}{-} & \multicolumn{1}{c}{\checkmark} & \multicolumn{1}{c}{-} & \multicolumn{1}{c}{-} & \multicolumn{1}{c}{-} & \multicolumn{1}{c}{\checkmark} \\
    Vehicle time window & \multicolumn{1}{c}{\checkmark} & \multicolumn{1}{c}{-} & \multicolumn{1}{c}{-} & \multicolumn{1}{c}{-} & \multicolumn{1}{c}{-} & \multicolumn{1}{c}{-} & \multicolumn{1}{c}{-} & \multicolumn{1}{c}{-} & \multicolumn{1}{c}{-} & \multicolumn{1}{c}{-} & \multicolumn{1}{c}{\checkmark} \\
    Satellite capacity & \multicolumn{1}{c}{-} & \multicolumn{1}{c}{\checkmark} & \multicolumn{1}{c}{-} & \multicolumn{1}{c}{\checkmark} & \multicolumn{1}{c}{\checkmark} & \multicolumn{1}{c}{-} & \multicolumn{1}{c}{\checkmark} & \multicolumn{1}{c}{-} & \multicolumn{1}{c}{\checkmark} & \multicolumn{1}{c}{-} & \multicolumn{1}{c}{\checkmark} \\
    Vehicle capacity & \multicolumn{1}{c}{-} & \multicolumn{1}{c}{\checkmark} & \multicolumn{1}{c}{\checkmark} & \multicolumn{1}{c}{\checkmark} & \multicolumn{1}{c}{\checkmark} & \multicolumn{1}{c}{\checkmark} & \multicolumn{1}{c}{\checkmark} & \multicolumn{1}{c}{\checkmark} & \multicolumn{1}{c}{\checkmark} & \multicolumn{1}{c}{\checkmark} & \multicolumn{1}{c}{\checkmark} \\
    Direct shipment & \multicolumn{1}{c}{-} & \multicolumn{1}{c}{\checkmark} & \multicolumn{1}{c}{\checkmark} & \multicolumn{1}{c}{-} & \multicolumn{1}{c}{-} & \multicolumn{1}{c}{-} & \multicolumn{1}{c}{-} & \multicolumn{1}{c}{\checkmark} & \multicolumn{1}{c}{-} & \multicolumn{1}{c}{-} & \multicolumn{1}{c}{\checkmark} \\
    \bottomrule
    \end{tabular}
    }
  \label{table:Table1}
\end{table}

\subsection{Contribution}

To close the research gap outlined above, we aim to study the benefit of a two-echelon distribution structure that is as versatile as possible. To this end, we, for the first time, study a \gls{2e-lrp} that allows the placing of both stationary and mobile micro-depots while at the same time allowing for direct shipments from the first echelon depot to customers. Considering all of these distribution options allows us to precisely study the benefit of mobile micro-depots as well as the benefit of hybrid operations where customers can be served from both echelons. 

To solve scenarios of realistic size, we propose a problem-specific \gls{alns} algorithm that, beyond tailored operators, maintains a pool of promising mini-network configurations throughout the search and uses a set covering component to guide the search. To improve computational times for large-scale instances, we propose an additional decomposition scheme that utilizes a cluster-first, route-second approach, yielding a tangible trade-off between computational complexity and solution quality.

Beyond verifying our algorithm's performance on established benchmark data sets, we present a case study for the city of Munich, Germany, that allows for managerial analyses, which remains the focus of our results discussion. Our algorithm solves up to 2000 customers in Munich in minutes, and the results show that allowing direct shipment can reduce overall costs by 4.7\% and total emissions by 11\%. It also increases truck utilizations by 42\%. The decomposition approach also works well with large-scale instances and gives comparable solutions with 15 times faster run times than our baseline algorithm. Furthermore, we analyze trade-offs in different micro-depot configurations. The results show that locating stationary and mobile micro-depots while allowing direct shipments saves 5.9\% in total costs compared to traditional \gls{2e-lrp} with stationary micro-depots.

\vspace{-1em}\subsection{Organization}
The organization of this paper is as follows. Section \ref{problem-setting} details our problem setting. We then develop a metaheuristic solution approach in Section~\ref{alns}. In Section~\ref{exp-setting}, we detail benchmark instances and introduce our case study for the city of Munich. We discuss our computational findings in Section~\ref{results}. Finally, Section~\ref{conclusion} provides a conclusion and future research directions.
\glsreset{fe}
\glsreset{se}
\newcommand{\SumSet}[1]{\ensuremath{\sum_{\scriptstyle\mathclap{#1}}}}

\section{Problem Definition}\label{problem-setting}

\begin{figure}[t]
    \centering
    \includegraphics[clip, trim=3cm 4.1cm 1cm 2cm,width=0.9\textwidth]{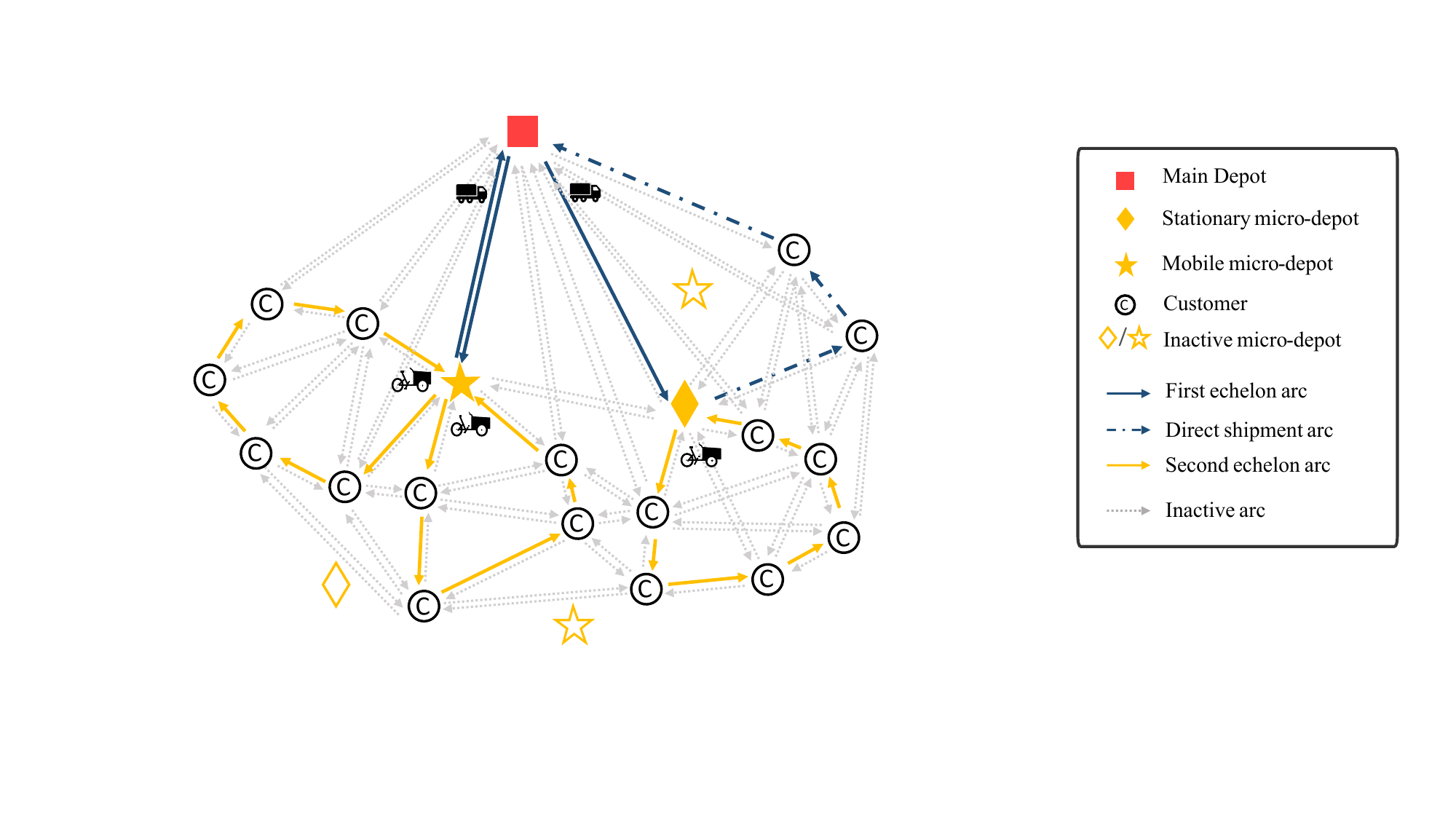}
    \caption{Illustration of the \gls{2e-lrp} with mobile depots and direct shipment network}
    \label{fig:2e-lrp-network}
\end{figure}

We focus on an \gls{lsp} that operates a \gls{2e-lrp} with mobile depots and direct shipment in the context of city logistics as shown in Figure~\ref{fig:2e-lrp-network}. The \gls{lsp} aims to transport parcels from a main depot to customers at minimum cost. To do so, the \gls{lsp} can use a two-echelon structure that further allows direct deliveries. For direct deliveries, the \gls{lsp} can utilize conventional trucks that start and end their routes at the main depot. Alternatively, the \gls{lsp} can transport parcels with conventional trucks to micro-depots, which act as hubs for last-mile delivery operations. From these micro-depots, secondary vehicles, e.g., environmentally friendly cargo bikes, perform last-mile customer delivery. 

To this end, the \gls{lsp} can utilize \textit{stationary} or \textit{mobile} micro-depots, which can both be located in the city center. \textit{Stationary micro-depots} are fixed logistics facilities that can once be strategically located in the city center. On the contrary, mobile micro-depots, e.g., swap bodies that can be preloaded at the main depot before being transported to the city center by truck, can be placed in the city center on demand. Both types of micro-depots have advantages and disadvantages. On the one hand, stationary micro-depots come at a higher cost due to their permanent construction but also accommodate larger parcel volumes, which can benefit operating a large fleet of secondary vehicles in areas with high and stable demand. On the other hand, mobile micro-depots have a significantly lower cost as they do not require permanent construction and provide flexibility to rearrange the micro-depot locations when delivery demand changes. Specifically, these modular units can be preloaded at the main depot each morning and transported to temporary parking locations near demand hotspots before being returned to the main depot in the evening.

In this setting, we put ourselves into the perspective of the \gls{lsp}, who aims to compute an optimal distribution plan for a representative, deterministic scenario. Computing such a distribution plan involves simultaneously determining the number and locations of stationary and mobile micro-depots, as well as the routing of \gls{fe} deliveries to micro-depots or customers and \gls{se} routes for last-mile deliveries.

In the following, we provide a pseudo-formal problem definition for conciseness and refer to a complete \gls{milp} definition to Appendix~~\ref{math-model} for brevity. 

\textit{Notation and solution representation:} To formally define our problem, let $\mGraph = (\mGraphVertex, \mGraphArc)$ be a directed graph consisting of a set of vertices $i \in \mGraphVertex$ and a set of arcs $(i,j) \in \mGraphArc$. The vertex set contains the main depot vertex 0, a set of potential micro-depot locations $\mSetDepots{}$ and a set of customers locations $\mSetCustomers$, such that $\mGraphVertex = \{0\}\cup\mSetDepots{}\cup\mSetCustomers$. We divide the arc set $\mGraphArc = \mSetArcs{f}\cup\mSetArcs{s}$ into two subsets $\mSetArcs{f}$ and $\mSetArcs{s}$, associated with routes in the \gls{fe} and \gls{se}, respectively. The set of potential micro-depots $\mSetDepots{}$ consists of stationary and mobile micro-depots. Each micro-depot $\mSatellite \in \mSetDepots{}$ has a limited capacity $\mCargoCapacity{}{\mSatellite}$, operating times $[\mEarliestBOS{\mSatellite}, \mLatestBOS{\mSatellite}]$ and a fixed opening cost $\mCostLocation{}{\mSatellite}$. Each customer $i \in \mSetCustomers$ has a demand $\mDemand{}{i}$, a service time $\mServiceTime{}{i}$ and a time-window $[\mEarliestBOS{i}, \mLatestBOS{i}]$. Each vehicle in the \gls{fe} $\mVehicleFirst \in \mSetTypeVehicle{f}$ and in the \gls{se} $\mVehicleSecond \in \mSetTypeVehicle{s}$ has limited capacity, denoted by $\mCargoCapacity{f}{\mVehicleFirst}$ and $\mCargoCapacity{s}{\mVehicleSecond}$. We use $\mFixedCost{\mVehicleFirst}$ to denote vehicle fixed costs. Traversing an arc $(i,j) \in \mGraphArc$ incurs a cost $\mCostArc{}{ij}$ depending on the distance traveled and vehicle used. We denote the binary variables for the location decision of micro-depots by $\mDecisionHubLocation{}{\mSatellite}$ and the vehicle decision on both echelons by $\mDecisionVehicle{}{\mTypeVehicle\mSatellite}$. Additionally, the binary variable $\mDecisionArc{}{\mNodeFirst\mNodeSecond\mTypeVehicle }$  indicates whether a vehicle $\mVehicleFirst$ travels on arc $(i,j)$.

A solution $\currentSol$ represents the sets of \gls{fe} and \gls{se} routes $\currentSol=\{\route{f}{1},...,\route{f}{\abs{\mSetTypeVehicle{f}}}, \route{s}{1},...,\route{s}{\abs{\mSetTypeVehicle{s}}}\}$. Each route in the \gls{fe} $\route{f}{k}=\{0,...,n, 0\}$ is a sequence of nodes $i \in \mGraphVertex$. Note that the starting and ending node is always the main depot. The \gls{se} route $\route{s}{k}=\{m,...,n, m\}$ consists of nodes $i \in \mGraphVertex\setminus\{0\}$. Here, the starting and ending node is always a micro-depot.

\textit{Objective function:} The \gls{lsp}'s objective is to minimize the total costs, which consist of fixed costs for opening stationary and mobile micro-depots as well as travel costs and fixed costs of vehicles used on both echelons. 
\begin{multline}
Z(s) = \SumSet{\mSatellite\in\mSetHubLocations}\mCostLocation{}{\mSatellite} \mDecisionHubLocation{}{\mSatellite} + \SumSet{\mVehicleFirst\in
	\mSetTypeVehicle{f}}~~~\SumSet{(\mNodeFirst,\mNodeSecond)\in\mSetArcs{f}} \mCostArc{}{\mNodeFirst\mNodeSecond} \mDecisionArc{}{\mNodeFirst\mNodeSecond\mVehicleFirst} + \SumSet{\mTypeVehicle\in
	\mSetTypeVehicle{s}}~~~\SumSet{(\mNodeFirst,\mNodeSecond)\in\mSetArcs{s}} \mCostArc{}{\mNodeFirst\mNodeSecond} \mDecisionArc{}{\mNodeFirst\mNodeSecond\mTypeVehicle }{}
	 + \SumSet{\mVehicleFirst\in\mSetTypeVehicle{f}} ~\SumSet{\mSatellite\in\mSetHubLocations} \mFixedCost{\mVehicleFirst} \mDecisionVehicle{}{\mVehicleFirst\mSatellite}
	 	 + \SumSet{\mTypeVehicle\in\mSetTypeVehicle{s}} ~ \SumSet{\mSatellite\in\mSetHubLocations} \mFixedCost{\mTypeVehicle} \mDecisionVehicle{}{\mTypeVehicle\mSatellite} \\
\label{obj: objective}
\end{multline}

\textit{Constraints:} A valid solution $\currentSol$ must comply with the following constraints. 

\begin{enumerate}[label=\roman*), itemsep=-0.5mm]
    \item The \gls{fe} vehicles start and end their routes in the main depot.
    \item The \gls{se} vehicles start and end their routes in the corresponding micro-depot.
    \item Vehicles can never exceed their capacity $\mCargoCapacity{}{\mVehicleFirst}$.
    \item Vehicles can depart and arrive at micro-depots within their time windows $[\mEarliestBOS{t}, \mLatestBOS{t}]$.
    \item Micro-depots can never exceed their capacity $\mCargoCapacity{}{\mSatellite}$.
    \item The customers must be served within their time windows $[\mEarliestBOS{i}, \mLatestBOS{i}]$. If a vehicle arrives early, it must wait until it serves the customer.
\end{enumerate}
Among all feasible solutions fulfilling these constraints, we seek a solution $\bestSol$ that minimizes objective function~\ref{obj: objective}.

One comment on our problem setting is in order: contrary to practice, where demand may change over time, our problem setting is based on a deterministic demand scenario. This assumption simplifies the problem formulation and allows for tractable optimization on large-scale instances. Clearly, expanding the problem setting studied in this paper to its stochastic counterpart remains a natural next step for future research. Despite this simplification, our results provide valuable insights into the potential benefits of integrating mobile micro-depots into urban logistics networks. The deterministic setting allows us to isolate and quantify the advantages of mobile micro-depots without the additional complexity of accounting for uncertainty, thereby laying a foundation for more advanced modeling approaches. In fact, if mobile micro-depots prove to be beneficial solely from a cost perspective in a deterministic scenario, one can see this as a strong indicator for potential cost savings under varying demand, where the savings potential is expected to be even higher when the mobile micro-depots flexibility comes into play.
\vspace{-1em}\section{Adaptive Large Neighborhood Search} \label{alns}

In this section, we introduce our \gls{alns}-based metaheuristic to solve \gls{2e-lrp} with mobile depots and direct shipment. \gls{alns} builds on the \gls{lns} framework proposed by \cite{Shaw1998}, where large neighborhoods are explored through a destroy and repair mechanism; in each iteration, a destroy operator removes a set of vertices from the current solution; subsequently, a repair operator constructs a new solution. This scheme allows to explore larger neighborhoods and helps the search process to escape local optima. \gls{alns}  enhances LNS by incorporating an adaptive selection mechanism that accounts for each operator's effectiveness whenever choosing destroy and repair operators in each iteration \citep{RopkePisinger2006}. 

Algorithm~\ref{alg:alns} shows the pseudocode of our metaheuristic. The algorithm begins with generating an initial solution (l.\ref{generate_initial_solution}). We then set the best solution $\bestSol$ to the initial solution $\currentSol$ (l.\ref{set_best}), and initialize a network pool $\mSetNetwork$ which we use to store route configurations found during the search~(l.\ref{init_network_pool}). 

In each search iteration, we select a destroy operator from the set of destroy operators (l.\ref{choose_destroy}) and apply it to obtain a partial solution $\newSol$ and a set of removed customers $\removedCustomers$ (l.\ref{destroy}). Subsequently, we choose a repair operator to repair the new partial solution $\newSol$ (l.\ref{choose_repair}-\ref{repair}) by reinserting the removed customers. After modifying the solution, we update the network pool $\mSetNetwork$ by adding the route configurations of $\newSol$ (l.\ref{update_network}). Finally, we apply local search to $\newSol$ to improve the \gls{se} routes (l.\ref{local_search_step}). 

Every \acrshort{sc_iter} iterations, we solve a \glsentrylong{sc} problem to find better network configurations and potentially update $\bestSol$ (l.\ref{sc_start}-\ref{sc_end}). We then apply a correction heuristic to find a feasible solution $\scSol$. If $\scSol$ is better than $\bestSol$, it is accepted as the new best solution. Otherwise, we continue with the next steps of our algorithm.

If we find an improved solution in terms of total cost, we update $\currentSol$ and $\mSetNetwork$. If the cost of $\newSol$, $f(\newSol)$, is less than $f(\bestSol)$, we also update $\bestSol$ (l.\ref{improvement_check}-\ref{update_best}). During the search, we use \gls{sa} to escape from local optima. If $\newSol$ is worse than $\currentSol$ but satisfies the \gls{sa} acceptance criterion, we accept it and continue the search from there (l.\ref{accept_criteria}). Finally, we update the current temperature \acrshort{currentTemp} and weights of the operators used in the current iteration (l.\ref{update_temp}-\ref{update_weight}). We repeat these steps until we reach the maximum number of iterations \acrshort{maxIter}.

In the following, we detail each algorithmic component (see Sections~\ref{init_section}-\ref{adaptive_mechanisms}) and a decomposition-based cluster-first, route-second approach (see Section~\ref{decomposition}). We refer the interested reader to Appendix \ref{ablation} for an ablation study in which we provide statistical evidence for the algorithmic components used.

\SetKwFunction{init}{GenerateInitialSolution}
\SetKwFunction{ChooseDestroyOperator}{ChooseDestroyOperator()}
\SetKwFunction{ChooseRepairOperator}{ChooseRepairOperator()}
\SetKwFunction{destroy}{Destroy}
\SetKwFunction{repair}{Repair}
\SetKwFunction{updateNetwork}{UpdateNetworkPool}
\SetKwFunction{setCover}{SetCover}
\SetKwFunction{correction}{CorrectionHeuristic}
\SetKwFunction{localSearch}{LocalSearch}
\SetKwFunction{UpdateTemperature}{UpdateTemperature}
\SetKwFunction{AdaptSearchParameters}{AdaptSearchParameters}
\SetKwFunction{accept}{Accept}

\begin{algorithm}[t!]
\small
\caption{ALNS-based Solution Algorithm}\label{alg:alns}
     $\currentSol \gets \init()$ \Comment{Section \ref{init_section}} \label{generate_initial_solution}\\
     $\bestSol \gets \currentSol$ ; \label{set_best}\\
     $\mSetNetwork \gets \emptyset$ ; \label{init_network_pool} \\
    \While{\acrshort{iter} $\leq$ \acrshort{maxIter}}{
        $\ChooseDestroyOperator$ \Comment{Section \ref{destroy_operators}}\label{choose_destroy}\\
        $(\newSol, \removedCustomers) \gets \destroy(\currentSol)$ ; \label{destroy}\\
        $\ChooseRepairOperator$ \Comment{Section \ref{repair_operators}}\label{choose_repair}\\ 
        $\newSol \gets \repair(\newSol, \removedCustomers)$ ; \label{repair} \\
        $\updateNetwork(\mSetNetwork, \newSol)$ ; \label{update_network}\\
        $\newSol \gets \localSearch(\newSol)$ \Comment{Section \ref{local_search_operators}}\label{local_search_step}\\
        \uIf{\acrshort{iter} $= 0 \pmod{$\acrshort{sc_iter}$}$}{\label{sc_start}
            $\bar{\mSetNetwork} \gets \setCover(\mSetNetwork)$ \Comment{Section \ref{set_cover_section}}\label{set_cover}\\
            $\scSol \gets \correction(\bar{\mSetNetwork})$ ; \label{correction}\\
            \uIf{f(\scSol) $\leq$ f(\bestSol)}{
                $\bestSol \gets \scSol$ ; \label{update_sc}\\
            }
            \uElse{ 
                $\updateNetwork(\mSetNetwork, \newSol)$ ;\label{sc_end}\\
            }
        }
        \uIf{f(\newSol) $\leq$ f(\currentSol)}{ \label{improvement_check}
            $\updateNetwork(\mSetNetwork, \newSol)$ ;\\
            $\currentSol \gets \newSol$ ; \label{update_current}\\
            \uIf{f(\newSol) $\leq$ f(\bestSol)}{ \label{best_check}
                $\bestSol \gets \newSol$ ; \label{update_best} \\
            }
        }
        \uElseIf{\accept(\newSol, \currentSol, \acrshort{currentTemp})}{
            $\currentSol \gets \newSol$ ; \label{accept_criteria}
        }
        \acrshort{currentTemp} $\gets \UpdateTemperature($\acrshort{currentTemp}, \acrshort{t_init}, \acrshort{t_final}, \acrshort{iter}$)$ \Comment{Section \ref{adaptive_mechanisms}} \label{update_temp} \\
        $\AdaptSearchParameters($\acrshort{weight}$)$ ; \label{update_weight}
    }
    \textbf{return} $\bestSol$
\end{algorithm}

\vspace{-1em}\subsection{Initial Solution Generation} \label{init_section}
To generate an initial solution, we design a two-step construction heuristic. First, we determine the initial open micro-depots and the customer assignments. We propose two approaches to determine the opened micro-depots: solving a \gls{cflp} and k-means clustering. When solving a \gls{cflp}, we decide on the micro-depot locations and the customer assignments. We define the objective function as a sum of the traveling cost from micro-depot $\mSatellite$ to customer $i$. A mathematical formulation of the \gls{cflp} can be found in Appendix \ref{cflp}. As a second approach, we apply k-means clustering~\citep{Hartigan1979} to cluster customer locations and select the micro-depot locations accordingly. 

While solving a \gls{cflp}, we inherently satisfy the capacity constraints at the micro-depots. However, if we use k-means clustering, we need to check for potential capacity violations. In case of any infeasibility, we repair the solution by reallocating the furthest customers from the current micro-depot to another available open micro-depot. After identifying initial open micro-depots, we generate the \gls{fe} routes by assigning each micro-depot to a \gls{fe} vehicle and then construct the \gls{se} routes using the \glsentrylong{cw} savings algorithm. 

\subsection{Destroy Operators}\label{destroy_operators}
We use eight destroy operators to remove customers and micro-depots from the solution. To this end, we distinguish between large operators that change the solution configuration by removing a micro-depot and small operators that only remove customers. In each iteration, our \gls{alns} randomly chooses the number of customers $q$ to be removed. Unless otherwise specified, we apply an upper bound $\Omega$ (see Section~\ref{parameter-tuning}) for large operators to limit $q$. For small operators, we limit $q$ to the maximum number of customers in the route. Our removal operators are as follows.

\begin{description}
    \item[Random Removal.] The random removal operator randomly selects a number of customers and removes them from the solution.

    \item[Random String Removal.]
    The random string removal operator randomly selects a \gls{se} route. Then, it selects a random string from the route and removes the customers in the selected string.
    
    \item[Furthest Customer Removal.]
    The furthest customer removal operator selects a random route and arranges its customers in descending order based on their distance from the assigned micro-depot. It then removes the furthest $q$ customers.
    
    \item[Micro-Depot Removal.]
    The micro-depot removal operator randomly chooses an open micro-depot and removes it from the solution \currentSol, along with the customers assigned and \gls{se} routes originating from that micro-depot. 
    
    \item[Partial Micro-Depot Removal.]
    The partial micro-depot removal operator is based on the micro-depot opening defined by \cite{Hemmelmayr2012}. It randomly chooses an open and a closed micro-depot. It then orders the assigned customers according to their vicinity to the closed micro-depot. Finally, it removes the closest $q$ customers from their current route, where~$q$ is a random number that does not exceed the total number of customers assigned to the open micro-depot. 
    
    \item[Partial Micro-Depot Swap.]
    The partial micro-depot swap operator randomly selects two open micro-depots. It then randomly selects $q_{1}$ and $q_{2}$  customers that are assigned to the respective micro-depots, where $q_{1}$ and $q_{2}$ are random numbers that do not exceed the total number of customers assigned to the respective micro-depots. Finally, it exchanges the micro-depot-customer assignments. 
    
    \item[\gls{fe} Route Removal.]
    The \gls{fe} route removal operator reduces the number of routes in the \gls{fe}. It selects the minimum utilized route based on Equation~(\ref{eq:fe_removal}), where $L(r)$ is the total load of route $r$ in the \gls{fe}. It removes the route and its micro-depot(s) from the solution. 
    
    \begin{equation}
    \route{1}{\text{min}} = \arg\min_{r \in R^{1}} L(r) \label{eq:fe_removal}
    \end{equation} 
    
    \item[\gls{se} Route Removal.]
    The \gls{se} route removal operator is similar to the \gls{fe} route removal. We select the route using Equation (\ref{eq:se_removal}) in the \gls{se} with minimum vehicle utilization in terms of total load and remove it from the solution. 
    \begin{equation}
    \route{2}{\text{min}} = \arg\min_{r \in R^{2}} L(r) \label{eq:se_removal}
    \end{equation} 
\end{description}

\subsection{Repair Operators}\label{repair_operators}

To repair the destroyed solution, by adding the removed customers, our \gls{alns} uses three repair operators: greedy and regret insertion as proposed in \cite{RopkePisinger2006}, as well as an additional repair operator that merges routes in the first echelon.

\begin{description}

\item[Greedy Insertion.]
The greedy insertion operator checks all possible insertion positions for each removed customer. Then, it inserts the customer into a route with minimum distance change. Since we have large removal operators that remove secondary vehicles or micro-depots from the solution, inserting the removed customers into the existing routes can sometimes be infeasible. In that case, we create a new route or even open a new micro-depot during the greedy insertion heuristic. 

Algorithm~\ref{alg:greedy} details the steps of the greedy insertion. We start with a destroyed solution $s$ that contains the subset of routes $\bar{R}$, the set of removed customers $\removedCustomers$, and the set of open micro-depots $\openDepots$. First, we sort the customers in ascending order according to their demand (l.\ref{greedy-sort}). For each customer in $\removedCustomers$, we check for the feasibility of each position at each route in $\bar{R}$. Then, we select the best insertion position that yields the minimum distance change (l.\ref{greedy-check-start}-\ref{greedy-check-end}). If there is no feasible position to insert the customer, we create a new secondary route and assign the customer to the route as described in the initial solution generation in Section~\ref{init_section} (l.\ref{greedy-create-route}). If any open micro-depot cannot serve the newly created route due to capacity limitations, we open a new micro-depot closest to the customer and assign the route to it (l.\ref{greedy-create-satellite}). 

\item[Regret Insertion.]
The regret insertion operator inserts the customers according to their regret value. Algorithm~\ref{alg:regret} shows the steps of regret insertion. We start with a destroyed solution $\newSol$ that contains the subset of routes $\bar{R}$, the set of removed customers $\removedCustomers$, and the set of open micro-depots $\openDepots$. For each customer in $\removedCustomers$, we check for feasibility for each position at each route in $\bar{R}$ (l.\ref{regret-check-start}-\ref{regret-check-end}). If there is no feasible insertion position, we create a new secondary route and assign the customer to it (l.\ref{regret-create-route}). If any open micro-depot cannot serve the newly created route due to capacity limitations, we open a new micro-depot closest to the customer and assign the route to it (l.\ref{regret-create-satellite}). Next, we calculate the regret value for each customer, defined as the difference in distance change between its best and second best insertion position (l.\ref{regret-value}). We then insert the customer with the highest regret value first using Equation \eqref{eqn:regret_function} (l.\ref{regret-select-highest}-\ref{regret-insert-highest}) and update the route and the customer's regret values again with the updated solution (l.\ref{regret-update}). Afterward, we continue the procedure until we have added all the customers who had been removed. 
\begin{equation}\label{eqn:regret_function}
    i := arg max_{i \in \removedCustomers} (\Delta f_{i}^{2} - \Delta f_{i}^{1})
\end{equation}

\SetKwFunction{feasibility}{CapacityCheck}
\SetKwFunction{twfeasibility}{TimeWindowCheck}
\SetKwFunction{createRoute}{CreateNewRoute}
\SetKwFunction{openSatellite}{OpenNewMicroDepot}
\SetKwFunction{sort}{Sort}
\SetKwFunction{insert}{Insert}
\SetKwFunction{update}{Update}
\SetKwFunction{position}{PositionAt}
\SetKwFunction{insertPosition}{InsertionPositions}

\begin{algorithm}[htbp]
\small
\caption{Greedy Insertion}\label{alg:greedy}
\KwInput{Destroyed solution $\newSol = \{\bar{R}, \removedCustomers, \openDepots\}$ where $\bar{R}$ is the subset of routes,  $\removedCustomers$ the set of removed customers and $\openDepots$ the set of opened micro-depots}
    $\sort(\removedCustomers)$; \label{greedy-sort}\\
    \For{$c \in \removedCustomers$}{
        $p^{*} \gets 0$; \\
        \For{$\bar{r} \in \bar{R}$ \label{greedy-check-start}} { 
            \uIf{\feasibility{$(c, \bar{r} ,\openDepots)$}}{
                \For{$p_{\bar{r}} \in \insertPosition(\bar{r})$} {
                    \uIf{\twfeasibility{$(c, p_{\bar{r}})$} and $f(p_{\bar{r}}) \leq f(p^{*})$}{
                        $p^{*} \gets  p_{\bar{r}}$; \label{greedy-check-end}\\
                    }
                }
            }
        }
        \uIf{$p^{*} = 0$}{
            $\bar{r} \gets \createRoute()$; \label{greedy-create-route}\\
            $\bar{R} \gets \bar{R} \cup \{\bar{r}\}$ ; \\
            $p^{*} \gets \position(\bar{r})$ ; \\
            \uIf{$\bar{r}$ cannot be assigned to existing micro-depots in $\openDepots$}{
                $\hat{t} \gets \openSatellite(c)$; \label{greedy-create-satellite}\\
                $\openDepots \gets \openDepots \cup \{\hat{t}\}$ ;\\
            }
        }
        $\insert(c,(p^{*},\bar{r}))$
    }
\end{algorithm}

\begin{algorithm}[htbp]
\small
\caption{Regret Insertion}\label{alg:regret}
\KwInput{Destroyed solution $\newSol = \{\bar{R}, \removedCustomers, \openDepots\}$ where $\bar{R}$ is the subset of routes, $\removedCustomers$ the set of removed customers, and $\openDepots$ the set of opened micro-depots} 

\For{$c \in \removedCustomers$}{
    $CP \gets \emptyset$; \\
    \For{$c \in \removedCustomers$}{
        $P \gets \emptyset$\;
        \For{$\bar{r} \in \bar{R}$} {  \label{regret-check-start}
            \uIf{\feasibility($c, \bar{r} ,\openDepots$)}{  
                \For{$p_{\bar{r}} \in \insertPosition(\bar{r})$} {
                    \uIf{$\twfeasibility(c, p_{\bar{r}})$}{
                        $P \gets  P \cup \{p_{\bar{r}}\}$; \label{regret-check-end}
                    }
                }
            }
        }
        \uIf{$P = \emptyset$}{
            $\bar{r} \gets \createRoute()$; \label{regret-create-route}\\
            $\bar{R} \gets \bar{R} \cup \{\bar{r}\}$ ; \\
            $P \gets P \cup \{\position(\bar{r})\}$ ; \\ 
            \uIf{$\bar{r}$ cannot be assigned to existing micro-depots in $\openDepots$}{
                $\hat{t} \gets \openSatellite(c)$; \label{regret-create-satellite}\\
                $\openDepots \gets \openDepots \cup \{\hat{t}\}$;
            }
        }
    }
    $\sort(P)$; \\
    $CP \gets CP \cup \{(\Delta f_{c}^{2} - \Delta f_{c}^{1}, p^{*}, c)\}$; \label{regret-value}
}
$(c^{*}, p^{*}, \bar{r}) \gets \max(CP)$; \label{regret-select-highest}\\
$\insert(c^{*}, (p^{*}, \bar{r}))$; \label{regret-insert-highest}\\
$\update(\newSol)$; \label{regret-update}
\end{algorithm}

\item[Merge Routes in First Echelon.]
The merge routes operator merges the primary routes in the first echelon after applying the destroy operator \textit{First Echelon Route Removal}.~The primary routes removed in the destroy operator contain the micro-depots and the customers directly served by the primary vehicles.~This operator aims to merge those removed primary routes with the remaining ones.~We start by sorting the total load of the removed micro-depots in descending order.~For each remaining primary route, we check the capacity constraints.~If there is enough capacity, we assign the micro-depot to the route.~We also try to insert the directly served customers into the remaining primary routes.~If there are still unserved micro-depots or customers, we create a new primary route. 

\end{description}

\subsection{Local Search}\label{local_search_operators}
We sequentially apply five local search operators at each iteration to intensify the search and improve the second echelon routes locally. We apply preprocessing to speed up the computation time when applying local search (see Appendix~\ref{speed-up}). To this end, we determine possible neighboring nodes by restricting the number of neighbors for each node and allowing $\neighborCount$ neighbors for each node. The nodes within the neighborhood are then sorted according to their distance to the selected node. We proceed by accepting the first improvement and continue with the next operator. The local search stops if no further improvement can be achieved by any operators. We use two intra-route operators to consider moves within the route and three inter-route operators to consider moves between different routes. For the intra-route operators, we only check for time-window violations. For the inter-route operators, we additionally check for capacity constraints. Since the micro-depot assignment of a customer may also change, we need to check for customer-micro-depot feasibility.

\begin{description}
\item[Relocate intra-route.]
The relocate intra-route operator works by selecting a customer within a route and relocating them to a different position within the same route, aiming to reduce the route's total cost. The relocation takes place only if the time-window constraints are met. For each customer $i$, the operator checks every neighboring customer $j$. If $j$ is in the same route as $i$, we remove the arcs $(i^{-}, i), (i, i^{+}), (j, j^{+})$ and add the arcs $(i^{-}, i^{+}), (j, i), (i, j^{+})$. If this reduces the total cost, the operator relocates $i$ next to $j$.  

\item[2-opt.]
The 2-opt operator checks every neighboring customer $j$ and reverses the chain in the route for each customer $i$. If $j$ is in the same route as $i$, we remove the arcs $(i, i^{+}), (j, j^{+})$ and add the arcs $(i, j), (i^{+}, j^{+})$. If this reduces the total cost, the operator places $j$ next to $i$ and reverses the direction of the customers between $i^{+}$ and $j$.

\item[Relocate inter-route.]
The relocate inter-route operator relocates a customer from its current route to another route. For each customer $i$, it tries to relocate it from route $R1$ next to its neighbor $j$, which is in another route $R2$. It removes the arcs $(i^{-}, i), (i, i^{+}), (j, j^{+})$ and adds the arcs $(i^{-}, i^{+}), (j, i), (i, j^{+})$. If this reduces the total cost, the operator relocates the customer into $R2$.

\item[2-opt*.]
The 2-opt* operator moves the customers served after customer $i$ on route $R1$ to be served after customer $j$ on route $R2$ and relocates the customers after $j$ on route $R2$ to route $R1$. To calculate the cost change, it removes the arcs  $(i, i^{+}), (j, j^{+})$, then adds the arcs $(i, j^{+}), (j, i^{+})$. The operator accepts the change if this adjustment results in a lower total cost and meets capacity and time-window constraints.

\item[Exchange.]
The exchange operator selects two customers in two different routes and exchanges their positions. It removes the arcs $(i^{-}, i), (i, i^{+}), (j^{-}, j), (j, j^{+})$ and adds the arcs $(i^{-}, j), (j, i^{+}),$ $(j^{-}, i), (i, j^{+})$. If there is an improvement in the total cost, the operator exchanges the customer positions.

\end{description}

\subsection{Set Cover Problem}\label{set_cover_section}
During the search, we solve a \glsentrylong{sc} problem in every \acrshort{sc_iter} iteration (\ref{alg:alns}.\ref{set_cover}) to find better network configurations. In this context, we define a mini-network $\mathfrak{n} \in \mSetMiniNetwork$, which consists of an open micro-depot $\mSatellite \in \mSetDepots{}$ and the customers $i \in \mSetCustomers$ assigned to it. We define binary variables $a_{i\mathfrak{n}}$ to denote whether a customer $i$ is included in the mini-network $\mathfrak{n}$ ($a_{i\mathfrak{n}} = 1$), or not ($a_{i\mathfrak{n}} = 0$). Similarly, $b_{\mathfrak{n}\mSatellite}$ denotes whether mini-network $\mathfrak{n}$ belongs to micro-depot $\mSatellite$. Finally, we use $f_{\mathfrak{n}}$ to denote the sum of fixed costs and travel costs of the mini-network. The \glsentrylong{sc} problem decides which mini-networks to use. Accordingly, objective (\ref{eq:sc_obj}) minimizes the total cost of the mini-networks chosen in the solution configuration. The model ensures that each customer is covered by at least one mini-network (\ref{eq:sc_1}), and for each micro-depot, at most one mini-network is selected (\ref{eq:sc_2}).

\begin{equation}
    \min  \SumSet{\mathfrak{n} \in \mSetMiniNetwork} f_{\mathfrak{n}} y_{\mathfrak{n}} \label{eq:sc_obj}
\end{equation} 
s.t.
\begin{align}
    \SumSet {\mathfrak{n} \in \mSetMiniNetwork} a_{i\mathfrak{n}} y_{\mathfrak{n}} \geq 1  \quad\quad \forall i\in C \label{eq:sc_1} \\
     \SumSet{\mathfrak{n} \in \mSetMiniNetwork} b_{\mSatellite\mathfrak{n}} y_{\mathfrak{n}} \leq 1 \quad\quad \forall \mSatellite\in \mSetDepots{}  \label{eq:sc_2}
\end{align}

\paragraph{Correction Heuristic.} 
After solving the set cover problem, we generate a new solution by selecting the \gls{se} routes from the network pool, $\mSetNetwork$, which is generated from the selected mini-networks (\ref{alg:alns}.\ref{correction}). In some cases, a customer may be visited multiple times within the solution because several selected \gls{se} routes serve the same customer. To mitigate this shortcoming, we use a correction heuristic to remove the extra visits and generate a feasible solution. 

The correction heuristic works as follows: for each customer assigned to more than one route, we compute the distance change when the customer is removed from the route. We keep the customer on the route with the largest savings and remove it from all other routes. We repeat the procedure for all customers visited more than once. We then generate the first echelon routes as described in Section~\ref{init_section}. We calculate the costs of the new solution and output the resulting solution. We continue the improvement step using the resulting solution. If the newly found solution is the new best solution, we update \bestSol.

\paragraph{Network Pool Size Management.}
The network pool $\mSetNetwork$ consists of the mini-networks and the respective \gls{se} routes. We update the network pool during the search whenever we modify the solution after applying the destroy and repair or local search operators (\ref{alg:alns}.\ref{update_network}). We add newly discovered routes to the network pool during the search. However, the size of the pool becomes excessively large after several iterations. In some cases, the route is not changed completely, but the customers' order changes. We implement a filtering mechanism to manage the pool size and prevent redundancy. When inserting a new set of routes associated with a mini-network into the pool, we compare it with the existing routes. If the new route includes the same customers as an existing route, we evaluate the total cost of the mini-network, calculated as the sum of travel and vehicle costs, and keep the route with the lower cost. This approach ensures that the network pool remains efficient and maintains a manageable size.

\subsection{Adaptive Mechanisms}\label{adaptive_mechanisms}

Our algorithm uses the \gls{sa} acceptance criterion to accept deteriorating solutions. We accept a worse solution if $e^{-\frac{f(\newSol) - f(s^*)}{\tau_{i}}} > U(0, 1)$ where $f(\newSol)$ and $f(s^*)$ are the objective values of $\newSol$ and $s^*$, and $\tau_{i}$ is the current temperature. We decrease the temperature in each iteration using the linear temperature update function as follows:

\begin{equation}
    \text{\acrshort{currentTemp}}_{i} = \text{\acrshort{t_init}} - i \frac{\text{\acrshort{t_init}} - \text{\acrshort{t_final}}}{\text{\acrshort{maxIter}}}
\end{equation}\\

\noindent where \acrshort{t_init} and \acrshort{t_final} are predefined initial and final temperatures, respectively.

At each iteration, we select one destroy and one repair operator. The probability $\phi_{j}$ for choosing the operator $\rho_{j}$ is calculated using roulette wheel selection (Equation~\eqref{eqn:roulette}). In the beginning, each destroy and repair operator has the same weight. The weight update procedure for each operator is based on the work of \cite{RopkePisinger2010}. ALNS updates the weights of repair and destroy operators dynamically based on the success of the previous iterations in Equation~\eqref{eqn:weight_update}, where $\lambda \in [0,1]$ is the smoothing factor and $\omega$ is the score value. We define four different score values as proposed in \cite{RopkePisinger2010}. Specifically, our ALNS assigns $\omega_{1}$ if the solution is a new global best, $\omega_{2}$ if it improves the current solution, $\omega_{3}$ if the new solution is accepted, and $\omega_{4}$ if the new solution is rejected. 

\begin{equation}\label{eqn:roulette}
    \phi_{j} = \frac{\rho_{j}}{\sum_{k \in \Omega} \rho_{k}}
\end{equation}

\begin{equation}\label{eqn:weight_update}
\rho_{j} = \lambda\rho_{j} + (1 - \lambda)\omega
\end{equation}

\subsection{Decomposition} \label{decomposition}
%\vspace{-0.3cm}
Although our \gls{alns} algorithm incorporates several speed-up techniques to enhance runtime performance, solving large-scale instances still requires considerable computational time to achieve high-quality solutions. Existing decomposition strategies such as cluster-first-route-second reduce computational runtimes by dividing the problem into smaller subproblems. However, these methods often sacrifice solution quality as they neglect the complexity of the overall problem. Therefore, it is essential to improve runtimes while maintaining the quality of the final solution through a hybrid approach.

To overcome this limitation, we introduce a hybrid decomposition-based cluster-first, route-second approach that first decomposes the problem into smaller subproblems, focusing on micro-depot assignments and route optimization. After improving the routes within these subproblems, we aggregate them to construct a good-quality initial solution. This approach replaces the construction phase of our initial algorithm and provides a good starting solution for the \gls{alns}, allowing it to run for fewer iterations while maintaining a high solution quality. 

Figure~\ref{fig:decomposition} illustrates our hybrid decomposition-based cluster-first-route-second approach to solve large-scale instances. We decompose the problem into two levels: micro-depots and routes. After improving each route, we create an aggregated solution and apply our \gls{alns}. 
Algorithm~\ref{alg:decomposition} shows the detailed steps of our decomposition approach. We determine the initial open micro-depots as described in Section~\ref{init_section}(l.\ref{generate}). At the first level of decomposition, we have the micro-depot and customers assigned to it at each cluster. Next, we apply k-means clustering to find smaller clusters that will then determine the routes (l.\ref{second_level_cluster}). For each small cluster, we solve a \gls{tsp} (l.\ref{solve_tsp}). Here, we initialize the routes by solving the \glsentrylong{cw} savings algorithm and apply local search intra-route operators defined in Section~\ref{local_search_operators} to improve the routes. We then aggregate the routes and create one large solution (l.\ref{aggregate}). Finally, we apply our proposed ALNS with fewer number of iterations to the aggregated solution (l.\ref{decomp_alns}).

\vspace{2em}\hspace{-0.5cm}\begin{minipage}{0.5\textwidth}
    \centering
    \includegraphics[clip, trim=0.5cm 5cm 5cm 2cm, width=1.1\textwidth]{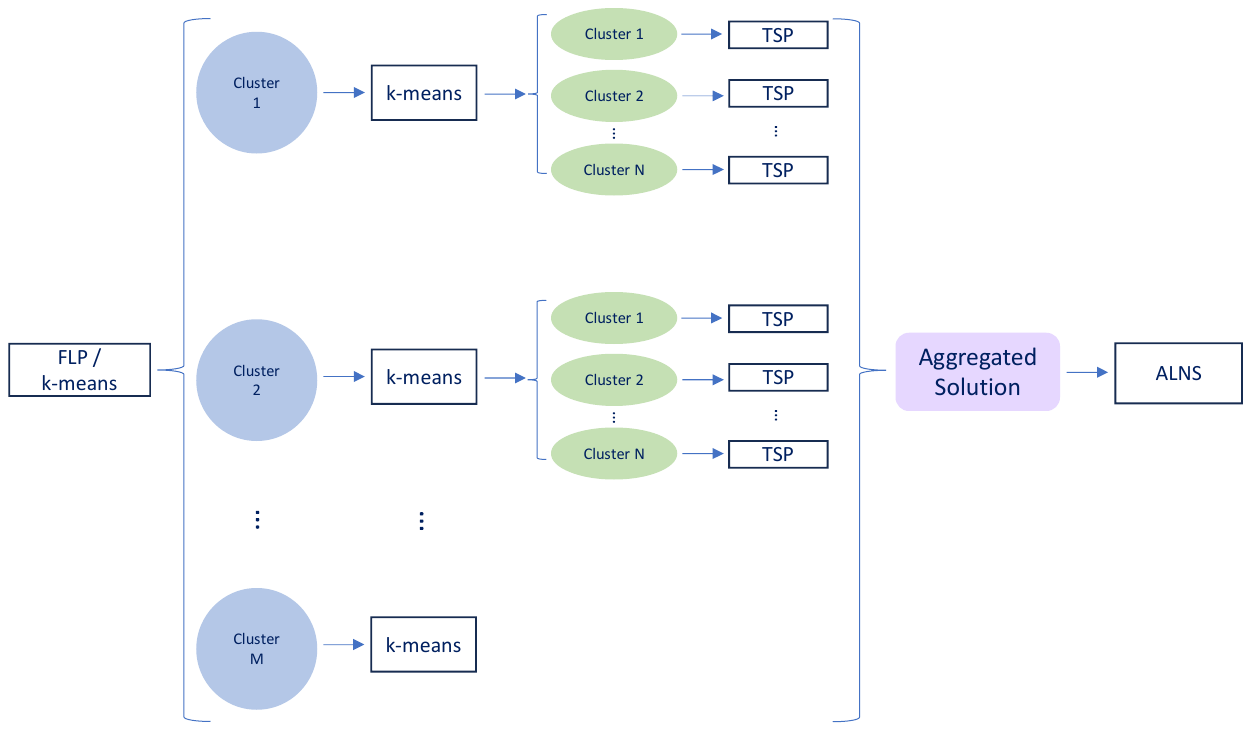}
    \captionof{figure}{Visualization of decomposition approach.} % Use captionof
    \label{fig:decomposition}
\end{minipage}
\hfill
\begin{minipage}{0.45\textwidth}
    \centering
    \footnotesize
    \begin{algorithm}[H]
        \caption{\small Decomposition approach.}
        \label{alg:decomposition}
        \KwIn{Instance $I$}
        \SetKwFunction{generateClusters}{GenerateFirstLevelCluster}
        \SetKwFunction{kmeans}{k-means}
        \SetKwFunction{tsp}{SolveTSP}
        \SetKwFunction{solve}{SolveCFLP}
        \SetKwFunction{combine}{Aggregate}
        \SetKwFunction{alns}{ALNS}
        
        \BlankLine
        $\mathbb{C} \gets \generateClusters(I)$ \label{generate} \\
        \ForEach{$\mathbb{C}_{i} \in \mathbb{C}$}{
%            \BlankLine
            \kmeans($\mathbb{C}_{i}$) \label{second_level_cluster}\\
            \ForEach{$C_{ij} \in \{C_{i1}, C_{i2}, \ldots, C_{iN}\}$}{
                \BlankLine
                \tsp{$\mathbb{C}_{ij}$} \label{solve_tsp}
            }
        }
        \BlankLine
        $s^{agg} \gets \combine(C_{ij})$ \label{aggregate}
        
%        \BlankLine
        $\bestSol\gets \alns({s^{agg}})$ \label{decomp_alns}
        
%        \BlankLine
        \Return{\bestSol}
    \end{algorithm}
\end{minipage}

\section{Experimental Design} \label{exp-setting}
To validate the performance of our \gls{alns} algorithm, we use benchmark datasets to provide a comparison against existing algorithms on basic \gls{2e-lrp} instances. Additionally, we generate case study instances for the city of Munich to derive managerial insights. Finally, we report our parameter tuning to find the best-performing parameters for our algorithm.

\subsection{Benchmark Instances} \label{benchmark_instances}
We use two 2E-LRP benchmark datasets proposed by \cite{Nguyen2012} to validate our algorithm's performance.~The first dataset is a modified version of Prodhon's 2E-LRP instances that consists of 30 instances with the following features: number of customers~$n\in\{20,50,100,200\}$ with uniform integer demands between 11 and 20, number of satellites $m\in\{5,10\}$, number of clusters $\beta\in\{1,2,3\}$, and vehicle capacity $Q \in \{70,150\}$. In addition to these features, the modified instances contain the main depot location, located at the origin $(0,0)$, and first echelon vehicles. The first echelon vehicle capacities are the maximum capacity of the satellites multiplied by $1.5$. The second dataset, Nguyen's 2E-LRP instances, contains 24 instances with the following features: number of customers $n \in \{ 25,50,100,200 \}$, number of satellites $m \in \{ 5,10 \}$, first-echelon vehicle capacity $Q_{1} \in \{ 750,850 \}$ and second-echelon vehicle capacity $Q_{2} \in \{ 100,150 \}$. Notably, customer demands follow a normal distribution with mean $\mu = 15$ and variance $\sigma^{2} = 25$. Customer locations are normally distributed ($\mathcal{N}$) or follow a multi-normal distribution ($\mathcal{MN}$).

For both datasets, the costs $c_{ij}$ correspond to the Euclidean distances rounded to the next integer, multiplied by 100. As also discussed in \cite{Vidal2018}, we use the following equations to calculate the distance matrices for the first (\ref{dm_first}) and second (\ref{dm_second}) echelon.

\begin{align}
&\mathit{d}_{AB} = \sqrt{(x_A - x_B)^2 + (y_A - y_B)^2} \times 100 \times 2 \label{dm_first} \\
&\mathit{d}_{AB} = \sqrt{(x_A - x_B)^2 + (y_A - y_B)^2} \times 100 \label{dm_second}
\end{align}

\subsection{Case Study Instances}
\begin{figure}[ht]
\begin{subfigure}[t]{.45\textwidth}
    \hspace{-1.2cm} % Move everything left
    \begin{minipage}[t]{\textwidth} % Ensures caption aligns with the image
        \centering
        \includegraphics[clip, trim=0cm 0cm 0cm 3.4cm, width=0.7\textwidth]{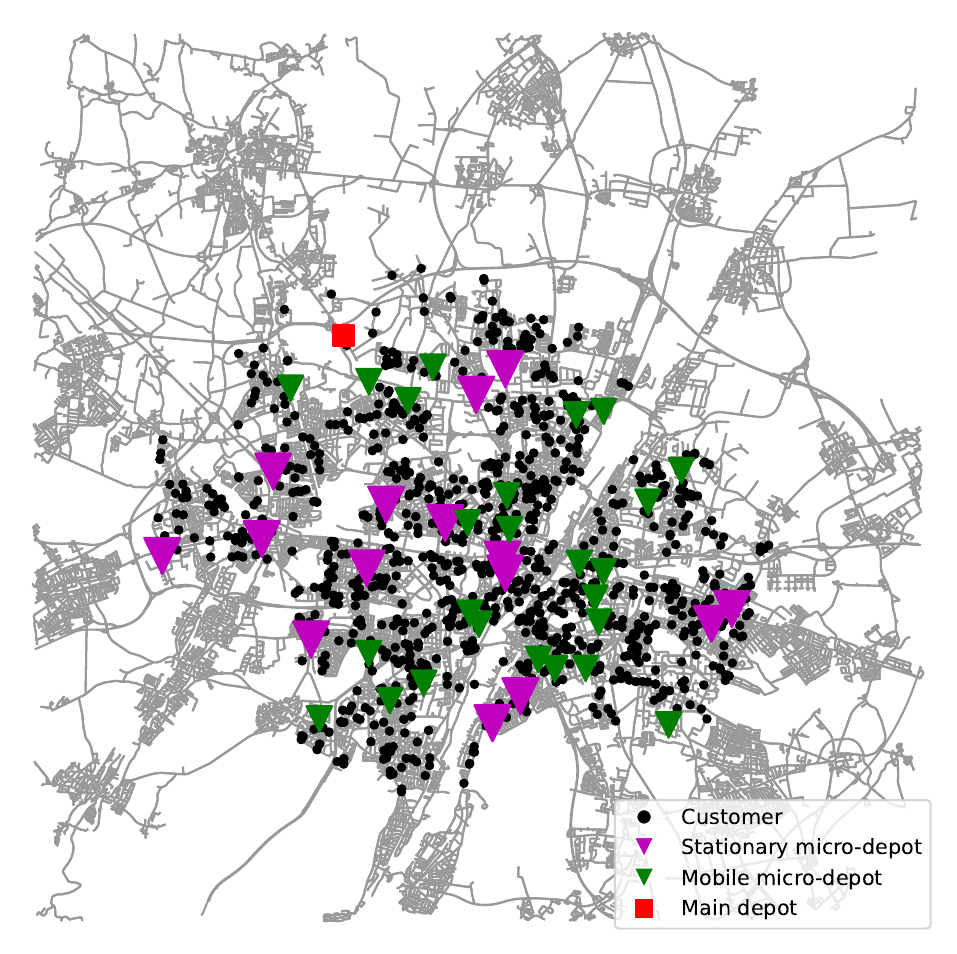}
        \caption{(a) Map of Munich showing the main depot, micro-depots and 2000 customer locations.}
        \label{fig:munich_cs}
    \end{minipage}
\end{subfigure}
    \begin{subfigure}[t]{.45\textwidth}
        \centering
        \includegraphics[width=0.7\textwidth]{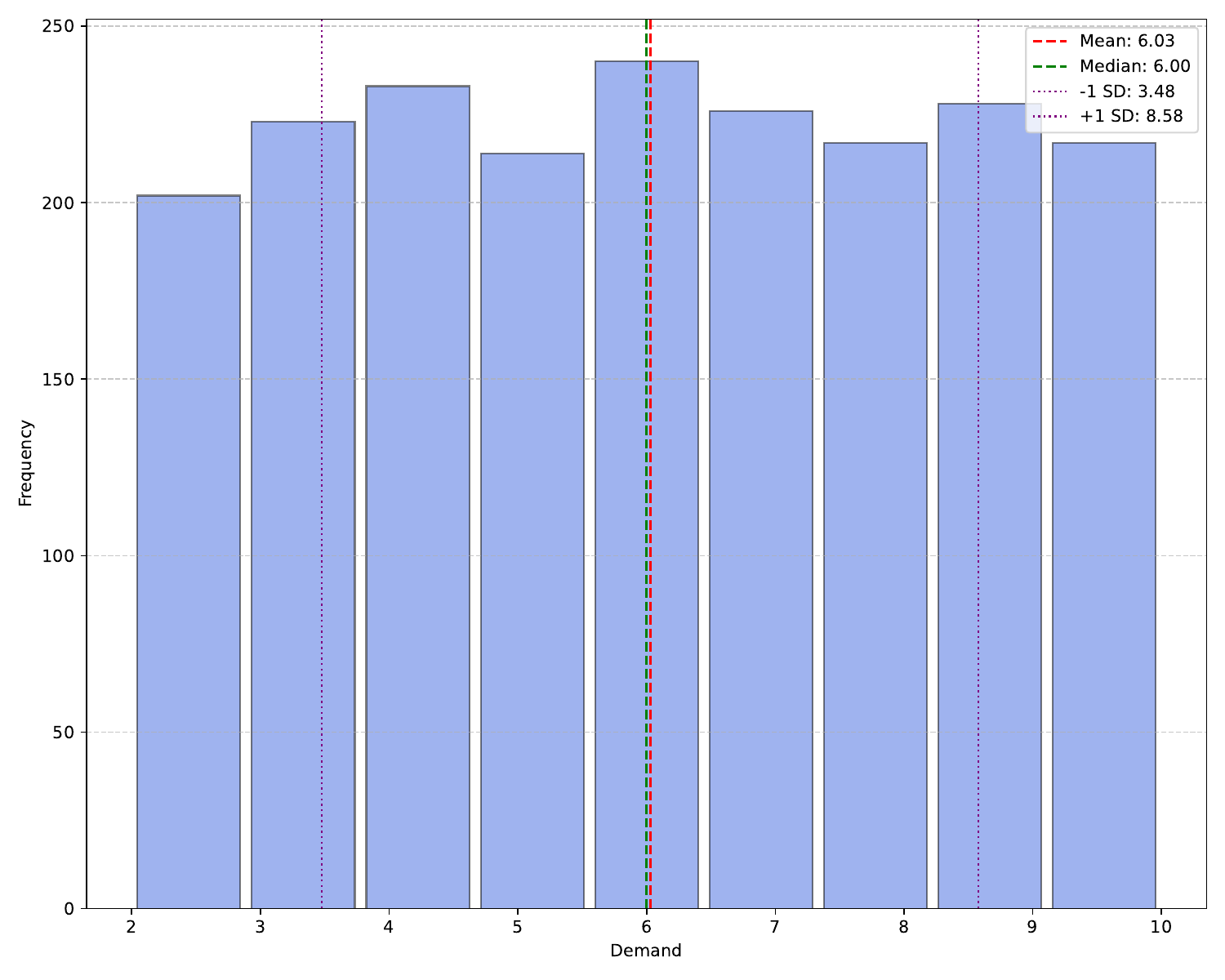}
        \caption{(b) Demand distribution of the customers.}
        \label{fig:demand_distribution}
    \end{subfigure}
\caption{Case study map for the city of Munich and demand distribution.}
\end{figure}

We generate instances for the city of Munich with 1000 and 2000 customers, 40 potential capacitated micro-depot locations, and one main depot outside the city center. Figure \ref{fig:munich_cs} represents the map of Munich with customer and depot locations. The number of customers and micro-depots assigned to each district depends on its population share. Subsequently, we randomly select the locations within the district. Customer demands follow a uniform distribution with mean $\mu = 6.03$ and standard deviation $\sigma = 2.55$ as shown in Figure~\ref{fig:demand_distribution}. We assign large trucks to be used in the first echelon and cargo bikes as city freighters in the second echelon. Instead of Euclidean distances, we generate the distance matrix for trucks and cargo bikes based on the actual street network in Munich using \gls{osm}~\citep{osm}.

Table~\ref{tab:Table2} represents the input parameters used for the case study. We take the costs and capacities for micro-depots and vehicles from \cite{Sheth2019}, \cite{Qi2018}, and \cite{Moeckl2020}. Based on these, we assume the costs of using a truck and cargo bike are \euro{480}/day and \euro{280}/day, respectively. The truck capacity is 865 cubic feet (cu. ft.), and the cargo bike capacity is 77 cu. ft. The capacities of stationary and mobile depots differ from each other. The former is the same as the truck capacity and is 865 cu. ft. and the latter is assumed to be 500 cu. ft. We assume the cost of locating a micro-depot is \euro{156}/day. Rather than assigning a fixed cost to stationary micro-depots, we model the cost as a uniform distribution within the range of \euro{600} to \euro{650}. Using this distribution, the assigned costs have a mean of \euro{628} and a standard deviation of \euro{15.08}.

\begin{table}[htbp]
\small
  \centering
  \begin{threeparttable}
  \caption{Input parameters of the case study}
    \begin{tabular}{lllr}
    \toprule
    \textbf{Description} & \textbf{Unit} & \textbf{Value} & \textbf{Source} \\
    \midrule
    Cost of truck & \euro{}/day & 480   & \cite{Sheth2019} \\
    Cost of cargo bike & \euro{}/day & 280   & \cite{Sheth2019} \\
    Capacity of truck & cu. ft. & 865   & \cite{Sheth2019} \\
    Capacity of cargo bike & cu. ft. & 77    & \cite{Sheth2019} \\
    Capacity of stationary micro-depot & cu. ft. & 865   & \cite{Sheth2019} \\
    Capacity of mobile micro-depot & cu. ft. & 500   &   \cite{Rossi2021}$^{*}$ \\
    Cost of stationary micro-depot & \euro{}/day & U(600,650) &  \cite{Sheth2019}\\
    Cost of mobile micro-depot & \euro{}/day & 156   &  \cite{Sheth2019}\\
    Fuel cost of truck & \euro{}/lt &  1.8  & \cite{Qi2018} \\
    Fuel consumption of truck per km & km/lt &  3  & \cite{Qi2018} \\
    Cost of cargo bike & \euro{}/km & 0.47   & \cite{Moeckl2020} \\
    $CO_{2}$ emission of truck & kg CO$_{2}$/km & 0.597   &  \cite{Qi2018} \\
    $CO_{2}$ emission of cargo bike & g CO$_{2}$/km &  0.079   &  \cite{Temporelli2022} \\
    \bottomrule
    \end{tabular}%
    \begin{tablenotes}
        \footnotesize
            \item[*]The unit is converted to cu.ft. for consistency.
    \end{tablenotes}
    \label{tab:Table2}%
  \end{threeparttable}
\end{table}

\vspace{-0.3cm}\subsection{Parameter Tuning}\label{parameter-tuning}
We conducted hyperparameter tuning to fine-tune the parameters of our \gls{alns} outlined by \cite{RopkePisinger2006}. To do so, we selected Prodhon's and Nguyen's \gls{2e-lrp} instances with more than 100 customers. To evaluate each parameter, we followed a sequential tuning approach: we varied one parameter at a time while fixing the remaining parameters at their initial values, as indicated by the superscript $i$ in Table~\ref{tab:Table3}. After identifying the best value for a parameter, we fixed it at its selected value and proceeded to tune the following parameter in the sequence. Table~\ref{tab:Table3} presents the tested values for each parameter and the final values used in the \gls{alns}. The bold values represent the final tuned parameters that yielded the best performance.

\glsreset{maxIter}
\glsreset{lambda}
\glsreset{weight}
\glsreset{t_init}
\glsreset{t_final}
\glsreset{temp_alpha}
\glsreset{reheat_iter}
\glsreset{neigh_count}
\glsreset{dod}
\glsreset{ls_counter}
% \glsreset{sc_iter}

Our \gls{alns} includes the following parameters: \gls{maxIter}, \gls{lambda}, \gls{weight}, \gls{t_init}, \gls{t_final}, \gls{temp_alpha}, \gls{temp_update}, \gls{reheat_iter}, \gls{neigh_count}, \gls{dod}, \gls{ls_counter} and the number of iterations after which the set cover problem is applied~(\acrshort{sc_iter}).~It is worth noting that changing \gls{dod} and \gls{reheat_iter} did not affect the solution quality. Therefore, we chose the maximum \gls{dod} value to allow for variation during the search. We further removed the temperature reheating step from our algorithm as it did not contribute to the solution quality. In the experiments, we observed the same solution quality when solving the \glsentrylong{sc} problem in every 500 and 2000 iteration. Since we aim for faster computation times, we selected \acrshort{sc_iter} as 2000.

\begin{table}[htbp]
\small
  \centering
  \caption{Parameter setting for the ALNS}
  \resizebox{1\textwidth}{!}{
    \begin{tabular}{ccccccccc}
    \toprule
    \gls{maxIter} & $50,000^{i}$ & \textbf{100,000} & $150,000$ &       & \textcolor[rgb]{ 0,  0,  0}{\gls{temp_update}} & \textbf{linear}$^{i}$ & sigmoid & hyperbolic \\
    \gls{lambda} &   \textbf{0.8}    &   $0.9^{i}$    &   $0.99$    &       & \gls{reheat_iter} & $10$ & $100$   & \textbf{0}$^{i}$ \\
    \gls{weight} &   $(50,10,1,0.5)$    &   $(100,50,10,0.5)^{i}$    &    \textbf{$(10,5,1,0)$}   &       & \gls{neigh_count} & \textbf{$0.2|\mathcal{C}|$} & $0.5|\mathcal{C}|$ & $1|\mathcal{C}|^{i}$ \\
    \gls{t_init} & \textbf{$10$} & $50$    & $100^{i}$   &       & \gls{dod} & $20$ & $0.5|\mathcal{C}|$ & \textbf{$1|\mathcal{C}|$}$^{i}$ \\
    \gls{t_final} & $0.01$  & $0.1^{i}$   & \textbf{$1$} &       & \gls{ls_counter} & \textbf{$1$} & $100$   & best$^{i}$ \\
    \gls{temp_alpha} & \textbf{$0.80$} & $0.90^{i}$  & $0.99$  &       & \acrshort{sc_iter} & $500^{i}$   & \textbf{$2000$} & $10000$ \\
    \bottomrule
    \end{tabular}%
    }
  \label{tab:Table3}%
\end{table}
\vspace{-2em}\section{Results} \label{results}
All experiments were conducted on a desktop computer equipped with an Intel(R) Core(TM) i9-9900, 3.1 GHz CPU, and 16 GB of RAM, running Ubuntu 20.04. We implemented our \gls{alns} as a single thread code in C\# .NET 6.0 and used Gurobi 10.0.0 to solve \gls{cflp} and \glsentrylong{sc} problems. To evaluate the performance of our algorithm, we present results from well-known benchmark datasets for the \gls{2e-lrp} and our proposed decomposition scheme in Section~\ref{computational_results}. Additionally, we include case study instances to show the scalability of our solution in Section~\ref{managerial_results}. For this case study, we further analyze the effects of direct shipment and conduct sensitivity analysis to derive managerial insights.

\subsection{Computational Analysis}\label{computational_results}
\paragraph{Benchmark Comparison.}
In this section, we evaluate the performance of our \gls{alns}; we performed a comprehensive comparative analysis against the best-known solutions derived from \cite{Sörensen2021} on Prodhon's and Nguyen's \gls{2e-lrp} benchmark datasets described in Section~\ref{benchmark_instances}. We perform ten runs for each instance using our \gls{alns}. From these runs, we record the best objective value achieved $(z^{*})$, the average of the best objective values across all runs ($\bar{z}$), the best gap to the BKS in percentage $\Delta(\%)$, and the average gap in percentage $\bar{\Delta}(\%)$. We calculate the gap between a solution with a value of z and the BKS using the formula $(z - BKS)/BKS *100$. Additionally, we record the average runtime per run ($\bar{t}$) in seconds. 

Table~\ref{tab:Table4} presents the detailed results on Prodhon's \gls{2e-lrp} instances. Out of 30 instances, our algorithm finds better solutions than the best-known solutions for two instances and matches the best-known solutions for seven instances. On average, our algorithm finds solutions that have 1.28\% gap to the best-known solutions over 30 instances. The average gap in percentage $\bar{\Delta}$ is below 5\% for most instances, where the overall average for $\bar{\Delta}$ is 2.45\%. Overall, we find comparable results against the best-known solutions in less than three minutes. For smaller instances containing 20 and 50 customers, the average computation time of our algorithm is less than a minute.

Table~\ref{tab:Table5} shows the detailed results on Nguyen's \gls{2e-lrp} instances. Out of 24 instances, our algorithm finds the best-known solutions for three instances. On average, our algorithm finds solutions that have 2.87\% gap to the best-known solutions for 24 instances. The average gap in percentage $\bar{\Delta}$ is below 5\% for most instances. Overall, we can also find comparable results against the best-known solutions in less than three minutes for Nguyen's dataset. For smaller instances, the average computation time of our algorithm is less than a minute.

\begin{table}[htbp]
  \centering
  \caption{Comparison of Prodhon's \gls{2e-lrp} best-known solutions derived from \citet{Sörensen2021} and the ALNS best results out of 10 runs and the average runtimes.}
  \resizebox{1\textwidth}{!}{
    \begin{tabular}{rrrrrrrrlrrrrrr}
    \toprule
          &       & \multicolumn{5}{c}{ALNS} &       &       &       & \multicolumn{5}{c}{ALNS} \\
\cmidrule{3-7}\cmidrule{11-15}    \multicolumn{1}{c}{Instance} & \multicolumn{1}{c}{BKS} & \multicolumn{1}{c}{$z^{*}$} & \multicolumn{1}{c}{$\bar{z}$}  & \multicolumn{1}{c}{$\Delta (\%)$} & \multicolumn{1}{c}{$\bar{\Delta}\%$}  & \multicolumn{1}{c}{$\bar{t}$(s)} &       & \multicolumn{1}{c}{Instance} & \multicolumn{1}{c}{BKS} & \multicolumn{1}{c}{$z^{*}$} & \multicolumn{1}{c}{$\bar{z}$}  & \multicolumn{1}{c}{$\Delta (\%)$} & \multicolumn{1}{c}{$\bar{\Delta}\%$}  & \multicolumn{1}{c}{$\bar{t}$(s)} \\
    \midrule
    \multicolumn{1}{l}{20-5-1} &      89,075  &      89,075  &      90,978  & \textbf{0.00} & 2.14  & 15.10 &       & 100-5-2b &    194,728  &    194,729  &    195,113  & \textbf{0.00} & 0.20  & 91.28 \\
    \multicolumn{1}{l}{20-5-1b} &      61,863  &      62,537  &      64,017  & 1.09  & 3.48 & 12.69 &       & 100-5-3 &    244,071  &    244,329  &    245,434  & 0.11  & 0.56  & 95.78 \\
    \multicolumn{1}{l}{20-5-2} &      84,478  &      84,478  &      85,696  & \textbf{0.00}  & 1.44 & 19.68 &       & 100-5-3b &    194,110  &    196,433  &    197,052  & 1.20  & 1.52  & 80.12 \\
    \multicolumn{1}{l}{20-5-2b} &      60,838  &      60,838  &      60,838  & \textbf{0.00} & \textbf{0.00} & 13.40 &       & 100-10-1 &    351,243  &    362,628  &    365,787  & 3.24  & 4.14  & 384.82 \\
    \multicolumn{1}{l}{50-5-1} &    130,843  &    137,225  &    137,245  & 4.88  & 4.89  & 43.90 &       & 100-10-1b &    297,167  &    306,231  &    313,750  & 3.05  & 5.58  & 546.68 \\
    \multicolumn{1}{l}{50-5-1b} &    101,530  &    103,016  &    107,903  & 1.46  & 6.28  & 31.30 &       & 100-10-2 &    304,438  &    305,129  &    308,031  & 0.23  & 1.18  & 117.74 \\
    \multicolumn{1}{l}{50-5-2} &    131,825  &    138,364  &    138,364  & 4.96  & 4.96  & 48.70 &       & 100-10-2b &    263,873  &    263,618  &    264,799  & \textbf{-0.10} & 0.35  & 134.32 \\
    \multicolumn{1}{l}{50-5-2b} &    110,332  &    113,126  &    116,877  & 2.53  & 5.93  & 29.40 &       & 100-10-3 &    310,148  &    317,486  &    324,930  & 2.37  & 4.77  & 271.40 \\
    \multicolumn{1}{l}{50-5-2BIS} &    122,599  &    124,902  &    128,216  & 1.88  & 4.58  & 56.50 &       & 100-10-3b &    260,328  &    265,812  &    274,825  & 2.11  & 5.57  & 451.96 \\
    \multicolumn{1}{l}{50-5-2BBIS} &    105,696  &    108,535  &    108,567  & 2.69  & 2.72  & 36.30 &       & 200-10-1 &    548,703  &    557,322  &    561,718  & 1.57  & 2.37  & 311.74 \\
    \multicolumn{1}{l}{50-5-3} &    128,379  &    128,379  &    129,218  & \textbf{0.00} & 0.65  & 38.30 &       & 200-10-1b &    445,301  &    446,300  &    452,753  & 0.22  & 1.67  & 283.56 \\
    \multicolumn{1}{l}{50-5-3b} &    104,006  &    104,006  &    104,006  & \textbf{0.00} & \textbf{0.00} & 37.00 &       & 200-10-2 &    497,451  &    498,827  &    500,045  & 0.28  & 0.52  & 291.82 \\
    \multicolumn{1}{l}{100-5-1} &    318,134  &    318,779  &    320,878  & 0.20  & 0.86  & 112.68 &       & 200-10-2b &    422,668  &    422,616  &    423,820  & \textbf{-0.01} & 0.27  & 313.22 \\
    \multicolumn{1}{l}{100-5-1b} &    256,878  &    256,888  &    259,007  & \textbf{0.00} & 0.83  & 161.60 &       & 200-10-3 &    527,162  &    530,991  &    533,617  & 0.73  & 1.22  & 290.68 \\
    \multicolumn{1}{l}{100-5-2} &    231,305  &    231,530  &    232,888  & 0.10  & 0.68  & 81.32 &       & 200-10-3b &    401,672  &    416,392  &    418,204  & 3.66  & 4.12  & 327.08 \\
    \midrule
          &       &       &       &       &       &       &       &  &       &       &    Average   & \textbf{1.28} & \textbf{2.45} & \textbf{157.67} \\
    \bottomrule
    \end{tabular}%  
    }  
  \label{tab:Table4}%
\end{table}%

% Table generated by Excel2LaTeX from sheet 'Nguyen_2016_BKS_final'
\vspace{-0.5em}\begin{table}[htbp]
  \centering
  \caption{Comparison of Nguyen's \gls{2e-lrp} best-known solutions derived from \citet{Sörensen2021} and the ALNS best results out of 10 runs and the average runtimes.}
  \resizebox{1\textwidth}{!}{
    \begin{tabular}{rrrrrrrrrrrlrrr}
    \toprule
          &       & \multicolumn{5}{c}{ALNS} &       &       &       & \multicolumn{5}{c}{ALNS} \\
\cmidrule{3-7}\cmidrule{11-15}    \multicolumn{1}{c}{Instance} & \multicolumn{1}{c}{BKS} & \multicolumn{1}{c}{$z^{*}$} & \multicolumn{1}{c}{$\bar{z}$}  & \multicolumn{1}{c}{$\Delta (\%)$} & \multicolumn{1}{c}{$\bar{\Delta}\%$}  & \multicolumn{1}{c}{$\bar{t}$(s)} &       & \multicolumn{1}{c}{Instance} & \multicolumn{1}{c}{BKS} & \multicolumn{1}{c}{$z^{*}$} & \multicolumn{1}{c}{$\bar{z}$}  & \multicolumn{1}{c}{$\Delta (\%)$} & \multicolumn{1}{c}{$\bar{\Delta}\%$}  & \multicolumn{1}{c}{$\bar{t}$(s)} \\
    \midrule
    \multicolumn{1}{l}{25-5N} &        80,370  &        84,892  &        87,488  & 5.63  & 8.86  & 16.4  &       & \multicolumn{1}{l}{100-5N} &     193,228  &     205,798  & \multicolumn{1}{r}{    206,801 } & 6.51  & 7.02  & 111.3 \\
    \multicolumn{1}{l}{25-5NB} &        64,562  &        64,562  &        64,562  & \textbf{0.00} & \textbf{0.00} & 15.2  &       & \multicolumn{1}{l}{100-5Nb} &     158,927  &     164,287  & \multicolumn{1}{r}{    164,412 } & 3.37  & 3.45  & 101.9 \\
    \multicolumn{1}{l}{25-5MN} &        78,947  &        78,947  &        78,947  & \textbf{0.00} & \textbf{0.00} & 18.2  &       & \multicolumn{1}{l}{100-5MN} &     204,682  &     205,394  & \multicolumn{1}{r}{    215,990 } & 0.35  & 5.52  & 101.8 \\
    \multicolumn{1}{l}{25-5MNb} &        64,438  &        64,438  &        64,438  & \textbf{0.00} & \textbf{0.00} & 18.1  &       & \multicolumn{1}{l}{100-5MNb} &     165,744  &     166,115  & \multicolumn{1}{r}{    168,240 } & 0.22  & 1.51  & 109.4 \\
    \multicolumn{1}{l}{50-5N} &     137,815  &     138,646  &     138,794  & 0.60  & 0.71  & 37.6  &       & \multicolumn{1}{l}{100-10N} &     209,952  &     233,230  & \multicolumn{1}{r}{    239,568 } & 11.09 & 14.11 & 145.6 \\
    \multicolumn{1}{l}{50-5Nb} &     110,094  &     112,737  &     112,737  & 2.40  & 2.40  & 35.7  &       & \multicolumn{1}{l}{100-10Nb} &     155,489  &     166,212  & \multicolumn{1}{r}{    169,424 } & 6.90  & 8.96  & 119.5 \\
    \multicolumn{1}{l}{50-5MN} &     123,484  &     128,793  &     128,793  & 4.30  & 4.30  & 36.4  &       & \multicolumn{1}{l}{100-10MN} &     201,275  &     209,255  & \multicolumn{1}{r}{    214,382 } & 3.96  & 6.51  & 124.3 \\
    \multicolumn{1}{l}{50-5MNb} &     105,401  &     106,313  &     106,313  & 0.87  & 0.87  & 37.1  &       & \multicolumn{1}{l}{100-10MNb} &     170,625  &     174,085  & \multicolumn{1}{r}{    176,567 } & 2.03  & 3.48  & 116.1 \\
    \multicolumn{1}{l}{50-10N} &     115,725  &     117,431  &     117,431  & 1.47  & 1.47  & 34    &       & \multicolumn{1}{l}{200-10N} &     343,232  &     361,536  & \multicolumn{1}{r}{    370,589 } & 3.57  & 7.97  & 453 \\
    \multicolumn{1}{l}{50-10Nb} &        87,315  &        87,686  &        87,686  & 0.42  & 0.42  & 40.2  &       & \multicolumn{1}{l}{200-10Nb} &     256,171  &     268,646  & \multicolumn{1}{r}{    278,550 } & 4.87  & 8.74  & 390 \\
    \multicolumn{1}{l}{50-10MN} &     135,519  &     139,664  &     139,664  & 3.06  & 3.06  & 34.5  &       & \multicolumn{1}{l}{200-10MN} &     323,801  &     338,782  & \multicolumn{1}{r}{    342,492 } & 4.63  & 5.77  & 426 \\
    \multicolumn{1}{l}{50-10MNb} &     110,613  &     111,290  &     114,389  & 0.61  & 3.41  & 55.5  &       & \multicolumn{1}{l}{200-10MNb} &     287,076  &     293,021  & \multicolumn{1}{r}{    306,601 } & 2.07  & 6.80  & 523 \\
    \midrule
          &       &       &       &       &       &       &       &       &       &       & Average & \textbf{2.87} & \textbf{4.39} & \textbf{129.20} \\
    \bottomrule
    \end{tabular}%
    }
  \label{tab:Table5}%
\end{table}%

\vspace{-1em}\paragraph{Decomposition Performance.}
This section evaluates the performance of the baseline \gls{2e-lrp} with mobile depot algorithms with direct shipment and our proposed decomposition-based cluster-first route-second models, emphasizing the trade-off between solution quality and computational efficiency. The baseline approaches include $\baseflp$, which initializes the solution using \acrshort{flp}, and $\basekmeans$, which employs a k-means clustering for initialization. Similarly, the decomposition models are defined as $\decompflp$ and $\decompkmeans$, based on the same respective initialization methods. 

Figure~\ref{fig:Figure4} compares these four approaches on total cost, computational time in seconds, number of micro-depots, and number of \gls{se} routes. Although the decomposition models show slightly higher total costs, on average 2.9\% more than the baseline algorithms, they offer substantial computational benefits, as illustrated in Figures~\ref{fig:Figure4_a} and \ref{fig:Figure4_b}. The decomposed models operate 15 times faster than the baseline approaches. To better understand the impact on the solution configuration, we assess the total number of open micro-depots in Figure~\ref{fig:Figure4_c}. The decomposition approaches produce configurations comparable to the baselines, with $\decompkmeans$ showing greater consistency by requiring fewer open micro-depots. Additionally, Figure~\ref{fig:Figure4_d} highlights the number of \gls{se} routes, demonstrating that although the decomposition models result in more routes, they provide more consistent and less variable solutions than the baseline methods.
\input{figures/decomposition_results}
In conclusion, the decomposition strategies offer a clear trade-off: reduced computational runtimes and improved consistency in solution configurations at the expense of modestly increased total costs. Among the decomposition models, $\decompflp$ delivers the most balanced performance across the evaluated metrics, including total cost, number of \gls{se} routes, and computational time.

\subsection{Managerial Insights}\label{managerial_results}
We consider a case study for the city of Munich, Germany, to derive managerial insights using the parameters presented in Table \ref{tab:Table2}. In this section, we compare two distribution network settings: the \acrshort{2e-lrp} with mobile depots and \acrshort{2e-lrp} with mobile depots and direct shipment. We provide the system analysis in \ref{system-analysis} followed by a sensitivity analysis in \ref{sensitivity-analysis} and the benefit of mobile micro-depots in~\ref{depot_configuration}.

\subsubsection{Effect of Direct Shipment} \label{system-analysis}
Our analysis compares the \gls{2e-lrp} with mobile depots and \gls{2e-lrp} with mobile depots and direct shipment settings. Figure~\ref{fig:Figure5} illustrates the comparison between these settings regarding the number of trucks, cargo bikes, and both stationary and mobile micro-depots.~The results show that the \gls{2e-lrp} with mobile depots and direct shipment reduces the number of truck routes from 20 to 19, cargo bike routes from 168 to 166, and mobile micro-depots from 12 to 11. However, it does not change the number of stationary micro-depots. This adjustment in routes and micro-depots underscores the efficiency gains and cost benefits of incorporating direct shipment into the logistics network.

\definecolor{RoyalBlue}{RGB}{65, 105, 225} % RGB values for Royal Blue
\definecolor{BlueGreen}{RGB}{13, 152, 186} % RGB values for Blue Green

% \vspace{-3cm}
\begin{figure}[htbp]
        \centering
        \begin{subfigure}[t]{0.20\textwidth}
            \centering
            \begin{tikzpicture}
            \centering
            %\node[anchor=north west] at (0.9, 0) {\small (a)};
                \begin{axis}[
                    width=4cm, height=4cm,
                    ybar,
                    ymin=14,
                    ymax=22,
                    enlarge x limits=0.8,
                    bar width=8pt,
                    symbolic x coords={\gls{2e-lrpmd}, \gls{2e-lrpmd-ds}},
                    xtick=\empty,
                    % nodes near coords,
                    % every node near coord/.append style={font=\fontsize{4}{5}\selectfont} % Rotate the labels vertically
                    tick label style={font=\footnotesize},
                    yticklabel style={font=\footnotesize},
                ]
                    % \addplot+[ybar, fill=teal] coordinates {(\gls{2e-lrpmd}, 20) (\gls{2e-lrpmd-ds}, 19)};
                    \addplot+[ybar, fill=RoyalBlue!50, draw=RoyalBlue!50] coordinates {(\gls{2e-lrpmd}, 20)};
                    \addplot+[ybar, fill=BlueGreen!50, draw=BlueGreen!50] coordinates {(\gls{2e-lrpmd-ds}, 19)};
                \end{axis}
                \node[anchor=north west] at (0.9, -0.15) {\small \centering \textbf{(a)}};
            \end{tikzpicture}

        \end{subfigure}
        \hspace{0.5cm}
        % Subfigure 2
        \begin{subfigure}[t]{0.20\textwidth}
            \centering
            \begin{tikzpicture}
                \begin{axis}[
                    width=4cm, height=4cm,
                    ybar,
                    ymin=162,
                    ymax=170,
                    enlarge x limits=0.5,
                    bar width=8pt,
                    symbolic x coords={\gls{2e-lrpmd}, \gls{2e-lrpmd-ds}},
                    xtick=\empty,
                    % nodes near coords,
                    tick label style={font=\footnotesize},
                    yticklabel style={font=\footnotesize},
                ]
                    % \addplot+[ybar, fill=teal] coordinates {(\gls{2e-lrpmd}, 168) (\gls{2e-lrpmd-ds}, 166)};
                    \addplot+[ybar, fill=RoyalBlue!50, draw=RoyalBlue!50] coordinates {(\gls{2e-lrpmd}, 168)};
                    \addplot+[ybar, fill=BlueGreen!50, draw=BlueGreen!50] coordinates {(\gls{2e-lrpmd-ds}, 166)};
                \end{axis}
                \node[anchor=north west] at (0.9, -0.15) {\small \centering \textbf{(b)}};
            \end{tikzpicture}
        \end{subfigure}
        \hspace{0.5cm}
        % Subfigure 3
        \begin{subfigure}[t]{0.20\textwidth}
            \centering
            \begin{tikzpicture}
                \begin{axis}[
                    width=4cm, height=4cm,
                    ybar,
                    ymin=2,
                    ymax = 10,
                    enlarge x limits=0.5,
                    bar width=8pt,
                    symbolic x coords={\gls{2e-lrpmd}, \gls{2e-lrpmd-ds}},
                    xtick=\empty,
                    % nodes near coords,
                    tick label style={font=\footnotesize},
                    yticklabel style={font=\footnotesize},
                ]
                    % \addplot+[ybar, fill=teal] coordinates {(\gls{2e-lrpmd}, 8) (\gls{2e-lrpmd-ds}, 8)};
                    \addplot+[ybar, fill=RoyalBlue!50, draw=RoyalBlue!50] coordinates {(\gls{2e-lrpmd}, 8)};
                    \addplot+[ybar, fill=BlueGreen!50, draw=BlueGreen!50] coordinates {(\gls{2e-lrpmd-ds}, 8)};
                \end{axis}
                \node[anchor=north west] at (0.9, -0.15) {\small \centering \textbf{(c)}};
            \end{tikzpicture}
        \end{subfigure}
        \hspace{0.5cm}
        % Subfigure 4
        \begin{subfigure}[t]{0.20\textwidth}
            \centering
            \begin{tikzpicture}
                \begin{axis}[
                    width=4cm, height=4cm,
                    ybar,
                    ymin=0,
                    ymax = 15,
                    enlarge x limits=0.5,
                    bar width=8pt,
                    symbolic x coords={\gls{2e-lrpmd}, \gls{2e-lrpmd-ds}},
                    xtick=\empty,
                    % nodes near coords,
                    tick label style={font=\footnotesize},
                    yticklabel style={font=\footnotesize},
                ]
                    % \addplot+[ybar, fill=teal] coordinates {(\gls{2e-lrpmd}, 12) (\gls{2e-lrpmd-ds}, 11)};
                    \addplot+[ybar, fill=RoyalBlue!50, draw=RoyalBlue!50] coordinates {(\gls{2e-lrpmd}, 12)};
                    \addplot+[ybar, fill=BlueGreen!50, draw=BlueGreen!50] coordinates {(\gls{2e-lrpmd-ds}, 11)};
                \end{axis}
                \node[anchor=north west] at (0.9, -0.15) {\small \centering \textbf{(d)}};
            \end{tikzpicture}
        \end{subfigure}
        
        % Horizontal legend spanning across the figure
       \vskip 0.5em % Add vertical spacing
        \begin{tikzpicture} 
            \begin{axis}[%
                hide axis,
                xmin=10,
                xmax=50,
                ymin=0,
                ymax=0.4,
                    legend style={
                    draw=white!15!black,
                    legend cell align=left,
                    legend columns=2, % Set 2 columns
                    column sep=0.5cm, % Adjust column spacing
                    row sep=0.5em % Adjust row spacing
                }
                ]
                \addlegendimage{fill=RoyalBlue!50, draw=RoyalBlue!50, line width=1.5pt} 
                \addlegendentry{\footnotesize\acrshort{2e-lrpmd}};
                \addlegendimage{fill=BlueGreen!50, draw=BlueGreen!50, line width=1.5pt}
                \addlegendentry{\footnotesize\acrshort{2e-lrpmd-ds}};
            \end{axis}
        \end{tikzpicture}

        \vskip 0.4em % Add vertical spacing
        \caption{ Comparison of \acrshort{2e-lrp} with mobile micro-depots (MD) and \acrshort{2e-lrpmd} with direct shipment (DS). Number of (a) trucks, (b) cargo bikes, (c) stationary micro-depots, (d) MDs.}

    \label{fig:Figure5}
\end{figure}
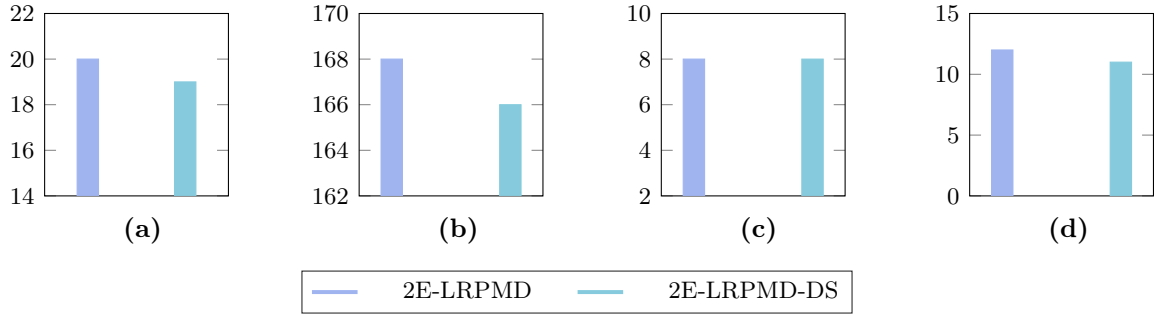
\FloatBarrier
\textbf{Insight 1.} The \gls{2e-lrp} with mobile depots and direct shipment setting achieved a 5\% reduction in truck routes, a 1.2\% reduction in cargo bike routes, and an 8\% decrease in the number of mobile micro-depots, highlighting the potential for cost efficiency through direct shipment integration.

Figure~\ref{fig:Figure6} illustrates the case study comparing the total costs, total emissions, truck and cargo bike utilization, and total distances traveled in the first and second echelons for the best solution obtained out of 10 runs. Figure~\ref{fig:Figure6_a} shows the total costs split into micro-depot, travel, truck, and cargo bike costs. Although there is an increase of 2.4\% in the micro-depot opening costs, the costs of trucks and cargo bikes decrease by 15\% and 2\%, respectively. The traveling cost increases slightly due to increased distance. Overall, there is a reduction in total costs of 4.7\% for the \gls{2e-lrp} with mobile depots and direct shipment. We also show in Figure~\ref{fig:Figure6_b} that allowing for direct shipments reduces total emissions by 11\%. In Figures~\ref{fig:Figure6_c} and \ref{fig:Figure6_d}, we compare the total truck and cargo bike utilizations. These plots illustrate that with the direct shipment, we can achieve increased truck and cargo bike utilizations of 42\% and 4\%, respectively. Additionally, we observe that the loads are more evenly distributed among the vehicles in both echelons. Finally, the total distances in both echelons are shown in Figures~\ref{fig:Figure6_e} and \ref{fig:Figure6_f}. Since customers are also served by trucks in the \gls{2e-lrp} with mobile depots and direct shipment, the distance traveled is increased by 8\% in the first echelon. For the second echelon, we also see an increase in the distances, resulting from using fewer cargo bikes and increased utilization for each cargo bike.

\definecolor{RoyalBlue}{RGB}{65, 105, 225} % RGB values for Royal Blue
\definecolor{BlueGreen}{RGB}{13, 152, 186} % RGB values for Blue Green
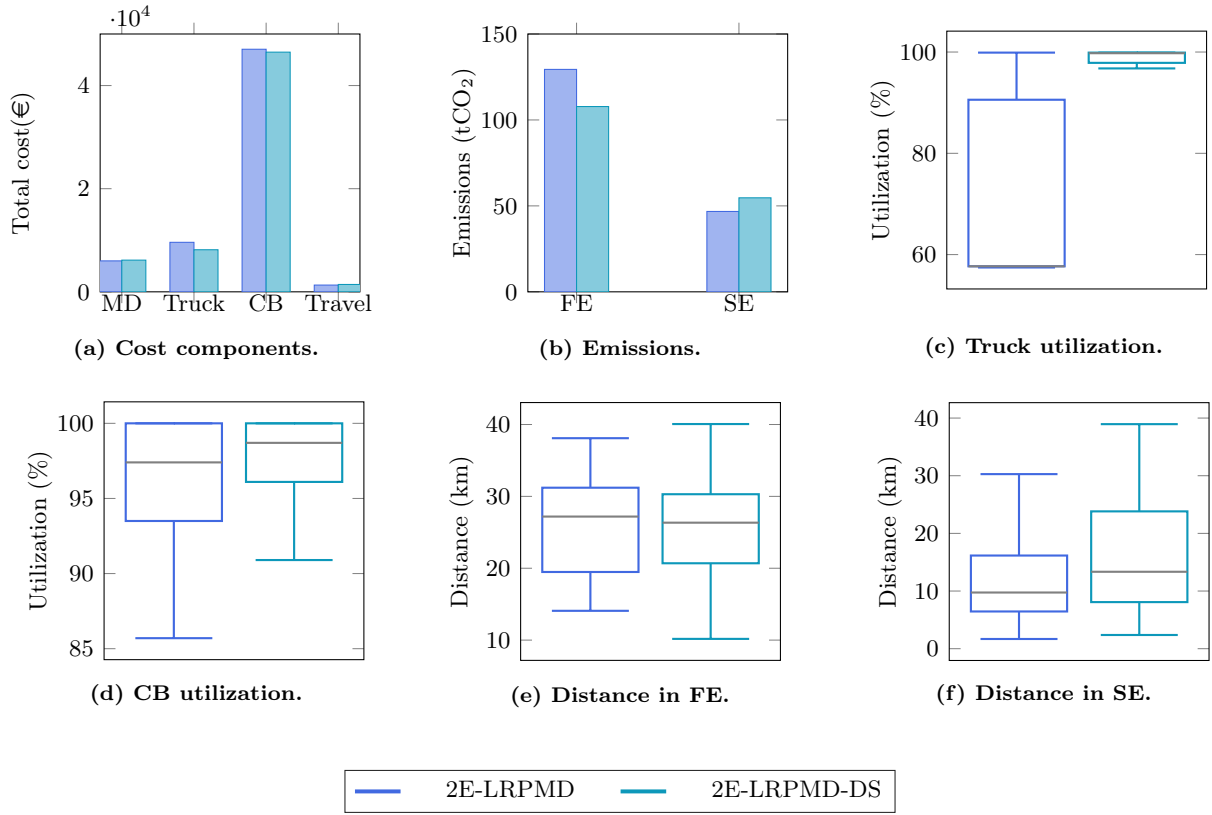
\begin{figure}[htbp]
    \centering
    % \scalebox{0.9}{%
    % First row
    \makebox[\textwidth][c]{
        \begin{subfigure}[t]{0.33\textwidth}
            \centering
                    \begin{tikzpicture}
            \begin{axis}[
                ybar=0pt,
                bar width=9pt,
                width=5cm,
                height=5cm,
                enlarge x limits=0.1,
                ylabel style={font=\footnotesize, yshift=0.3cm},
                ylabel={ Total cost(\euro{})},
                yticklabel style={font=\footnotesize},
                % symbolic x coords={Micro-depot, Truck, Cargo Bike , Travel},
                symbolic x coords={MD, Truck, CB, Travel},
                xtick=data,
                ymin=0,
                ymax = 50000,
                x tick label style={rotate=0, anchor=center, font=\footnotesize}, % Rotate x-axis labels
                legend style={at={(0.5,-0.2)},
                    anchor=north,legend columns=-1},
            ]
            % Data for 2E-LRPMD
            \addplot+[RoyalBlue, fill=RoyalBlue!50] 
                coordinates {(MD,6001.0) (Truck,9600.0) (CB,47040.0) (Travel,1309.54812)};
            % Data for 2E-LRPMD-DS
            \addplot+[BlueGreen, fill=BlueGreen!50] 
                coordinates {(MD,6147.0) (Truck,8160.0) (CB,46480.0) (Travel,1427.21807)};           
            \end{axis}
        \end{tikzpicture}
            \subcaption{(a) Cost components.}
            \label{fig:Figure6_a}
        \end{subfigure}        
        % \hspace{0.5cm} % Space between subfigures
        \begin{subfigure}[t]{0.33\textwidth}
            \centering
                        \begin{tikzpicture}
        \begin{axis}[
            ybar=0pt,
            bar width=12pt,
            width=5cm,
            height=5cm,
            enlarge x limits=0.3,
            ylabel style={font=\footnotesize, yshift=-0.2cm},
            ylabel={Emissions (tCO\textsubscript{2})},
            yticklabel style={font=\footnotesize},
            symbolic x coords={FE, SE},
            xtick=data,
            ymin=0,
            ymax=150, % Adjust this based on the range of your data
            x tick label style={rotate=0, anchor=center, font=\footnotesize}, % Rotate x-axis labels
            legend style={at={(0.5,-0.2)},
                anchor=north,legend columns=-1},
        ]
        
        % Data for 2E-LRPMD
        \addplot+[RoyalBlue, fill=RoyalBlue!50] 
            coordinates {(FE,129.406) (SE,46.760472)};

        % Data for 2E-LRPMD-DS
        \addplot+[BlueGreen, fill=BlueGreen!50] 
            coordinates {(FE,107.801) (SE,54.695542)};
        
        \end{axis}
    \end{tikzpicture}
            \subcaption{(b) Emissions.}
            \label{fig:Figure6_b}
        \end{subfigure}
        % \hspace{0.5cm} % Space between subfigures
        \begin{subfigure}[t]{0.33\textwidth}
            \centering
                \vspace{-3.85cm}\begin{tikzpicture}
            \begin{axis}[
                boxplot/draw direction=y,
                enlarge x limits=0.1,
                ylabel style={font=\footnotesize, yshift=-0.2cm},
                ylabel={Utilization (\%)},
                yticklabel style={font=\footnotesize},
                width=5cm,
                height=5cm,
                xtick=\empty, % No x-axis labels or ticks
                 % Option to show mean
                boxplot/median=0.6pt, % Thickness of the median line
            % boxplot/median style={solid, gray}, % Color and style of the median line
                boxplot/every median/.style={thick,solid,gray}
                ]        
        % Boxplot for 2E-LRP
        \addplot+[
            boxplot prepared={
                median=57.6878612716763, 
                upper quartile=90.57803468208093, 
                lower quartile=57.6878612716763, 
                upper whisker=99.88439306358381, 
                lower whisker=57.45664739884393
            },
            fill=none,
            draw=RoyalBlue,
            thick
        ] coordinates {};

        % Boxplot for 2E-LRP with Direct Shipment
        \addplot+[
            boxplot prepared={
                median=99.76878612716763, 
                upper quartile=99.88439306358381, 
                lower quartile=97.86127167630057, 
                upper whisker=99.88439306358381, 
                lower whisker=96.76300578034682
            },
            fill=none,
            draw=BlueGreen,
            thick
        ] coordinates {};

        \end{axis}
    \end{tikzpicture}      
            \vspace{0.4cm}\subcaption{(c) Truck utilization.}
            \label{fig:Figure6_c}
        \end{subfigure}
    }

    % Third row
    \vskip\baselineskip
    \makebox[\textwidth][c]{
        \begin{subfigure}[t]{0.33\textwidth}
            \centering
                        \begin{tikzpicture}
            \begin{axis}[
                boxplot/draw direction=y,
                enlarge x limits=0.1,
                ylabel style={font=\footnotesize, yshift=-0.2cm},
                ylabel={Utilization (\%)},
                yticklabel style={font=\footnotesize},
                width=5cm,
                height=5cm,
                xtick=\empty, % No x-axis labels or ticks
                 % Option to show mean
                boxplot/median=0.6pt, % Thickness of the median line
            % boxplot/median style={solid, gray}, % Color and style of the median line
                boxplot/every median/.style={thick,solid,gray}
                ]        
        % Boxplot for 2E-LRP
        \addplot+[
            boxplot prepared={
                median=97.4, 
                upper quartile=100.0, 
                lower quartile=93.5, 
                upper whisker=100.0, 
                lower whisker=85.7
            },
            fill=none,
            draw=RoyalBlue,
            thick
        ] coordinates {};

        % Boxplot for 2E-LRP with Direct Shipment
        \addplot+[
            boxplot prepared={
                median=98.7, 
                upper quartile=100.0, 
                lower quartile=96.1, 
                upper whisker=100.0, 
                lower whisker=90.9
            },
            fill=none,
            draw=BlueGreen,
            thick
        ] coordinates {};

        \end{axis}
    \end{tikzpicture}  
            \subcaption{(d) CB utilization.}
            \label{fig:Figure6_d}
        \end{subfigure}
        % \hspace{0.5cm} % Space between subfigures

        \begin{subfigure}[t]{0.33\textwidth}
            \centering
                        \vspace{-3.5cm}\begin{tikzpicture}
            \begin{axis}[
                boxplot/draw direction=y,
                enlarge x limits=0.1,
                ylabel style={font=\footnotesize, yshift=-0.1cm},
                ylabel={Distance (km)},
                yticklabel style={font=\footnotesize},
                width=5cm,
                height=5cm,
                xtick=\empty, % No x-axis labels or ticks
                 % Option to show mean
                boxplot/median=0.6pt, % Thickness of the median line
            % boxplot/median style={solid, gray}, % Color and style of the median line
                boxplot/every median/.style={thick,solid,gray}
                ]        
        % Boxplot for 2E-LRP
        \addplot+[
            boxplot prepared={
                median=27.19,
                %mean=25.8812,
                upper quartile=31.2065,
                lower quartile=19.4795,
                upper whisker=38.096,
                lower whisker=14.076,
            },
            fill=none,
            draw=RoyalBlue,
            thick
        ] coordinates {};

        % Boxplot for 2E-LRP with Direct Shipment
        \addplot+[
            boxplot prepared={
                median=26.34,
                %mean=27.2276,
                upper quartile=30.3,
                lower quartile=20.69275,
                upper whisker=40.065,
                lower whisker=10.17,
            },
            fill=none,
            draw=BlueGreen,
            thick
        ] coordinates {};
        \end{axis}
    \end{tikzpicture}
            \vspace{0.1cm}\subcaption{(e) Distance in \gls{fe}.}
            \label{fig:Figure6_e}
        \end{subfigure}
        
        % \hspace{0.5cm} % Space between subfigures
        
        \begin{subfigure}[t]{0.33\textwidth}
            \centering
                    \begin{tikzpicture}
            \begin{axis}[
                boxplot/draw direction=y,
                enlarge x limits=0.1,
                ylabel style={font=\footnotesize, yshift=-0.1cm},
                ylabel={Distance (km)},
                yticklabel style={font=\footnotesize},
                width=5cm,
                height=5cm,
                xtick=\empty, % No x-axis labels or ticks
                 % Option to show mean
                boxplot/median=0.6pt, % Thickness of the median line
            % boxplot/median style={solid, gray}, % Color and style of the median line
                boxplot/every median/.style={thick,solid,gray}
                ]        
        % Boxplot for 2E-LRP
        \addplot+[
            boxplot prepared={
                median=9.744,
                upper quartile=16.157,
                lower quartile=6.45625,
                upper whisker=30.277,
                lower whisker=1.683,
            },
            fill=none,
            draw=RoyalBlue,
            thick
        ] coordinates {};

        % Boxplot for 2E-LRP with Direct Shipment
        \addplot+[
            boxplot prepared={
                median=13.342,
                upper quartile=23.817,
                lower quartile=8.08,
                upper whisker=38.945,
                lower whisker=2.377,
            },
            fill=none,
            draw=BlueGreen,
            thick
        ] coordinates {};

        \end{axis}
    \end{tikzpicture}
            \subcaption{(f) Distance in \gls{se}.}
            \label{fig:Figure6_f}
        \end{subfigure}
    }
    
    % Legend
    \vskip 2em % Add vertical spacing
    \begin{tikzpicture} 
        \begin{axis}[
            hide axis,
            xmin=10,
            xmax=50,
            ymin=0,
            ymax=0.4,
            legend style={
                draw=white!15!black,
                legend cell align=left,
                legend columns=2, % Set 2 columns
                column sep=0.5cm, % Adjust column spacing
                row sep=0.5em % Adjust row spacing
            }
        ]
            \addlegendimage{fill=RoyalBlue, draw=RoyalBlue, line width=1.5pt} 
            \addlegendentry{\footnotesize \acrshort{2e-lrpmd}};
            \addlegendimage{fill=BlueGreen, draw=BlueGreen, line width=1.5pt}
            \addlegendentry{\footnotesize \acrshort{2e-lrpmd-ds}};
        \end{axis}
    \end{tikzpicture}
    % }  
    \caption{Comparison of \gls{2e-lrp} with mobile depots (MD) and \acrshort{2e-lrpmd} with direct shipment~(DS) settings.}
    \label{fig:Figure6}
\end{figure}
\textbf{Insight 2.} The \acrshort{2e-lrp} with mobile depots and direct shipment results in a 4.7\% reduction in total costs, 11\% in emissions, and an increase in \gls{fe} vehicle utilizations by 42\%.

\subsubsection{Sensitivity Analysis} \label{sensitivity-analysis}
\begin{figure}[htbp]
    \centering
    \includegraphics[scale = 0.8]{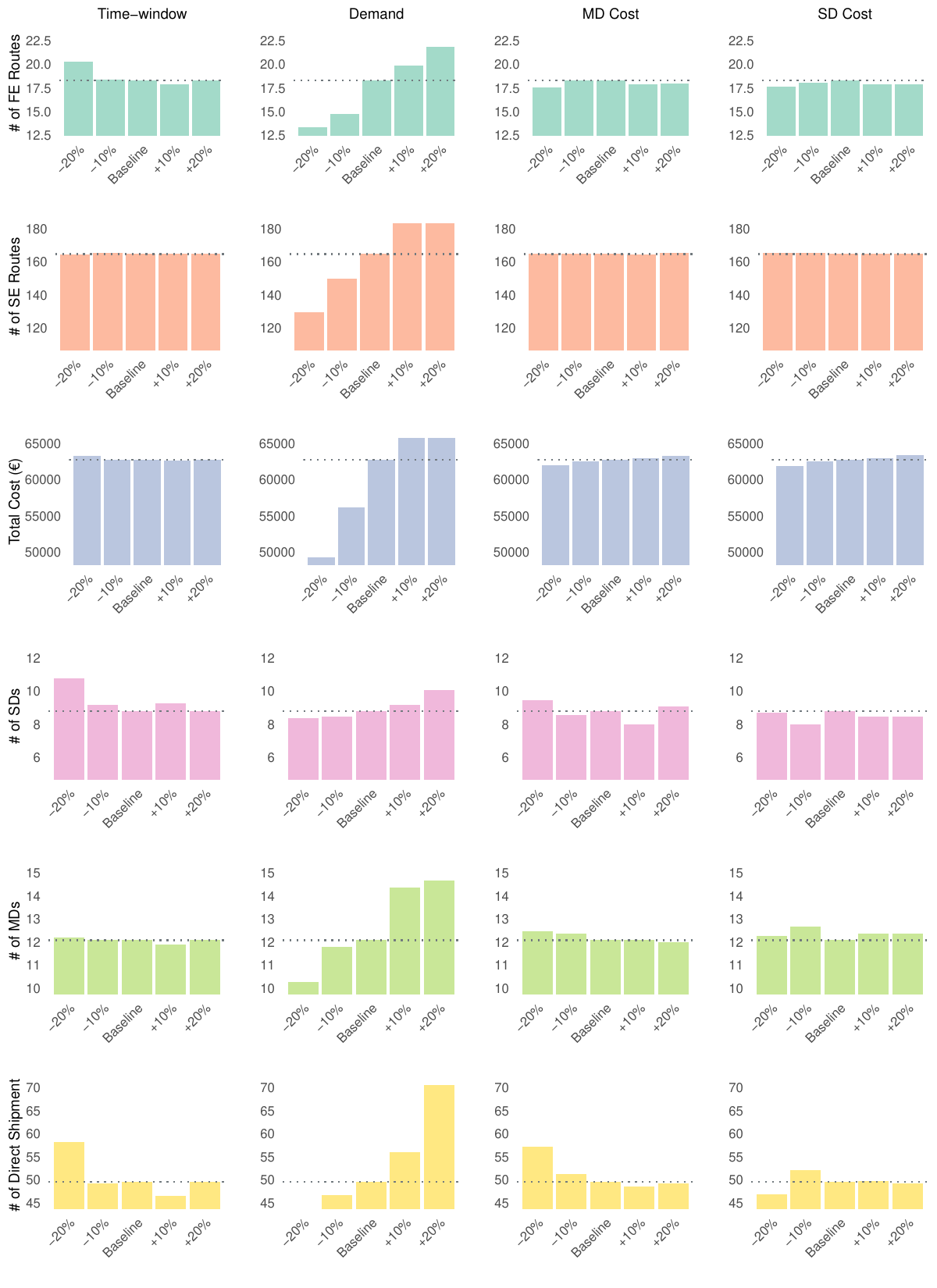}
    \caption{Sensitivity analysis of the case study instances for \gls{2e-lrp} with mobile depots and direct shipment. Columns indicate the adjusted input parameters: time windows, demand, mobile micro-depot (MD) cost, and stationary micro-depot (SD) cost. Rows represent the \glspl{kpi}: number of \gls{fe} routes, number of \gls{se} routes, total cost(\euro), number of SDs, number of MDs, and the number of direct shipments. The values shown are the average values out of 10 runs for each scenario.}
    \label{fig:sens-table}
\end{figure}

We conduct a sensitivity analysis to evaluate how variations in input parameters impact the performance of \gls{2e-lrp} with mobile depots and direct shipment. Using a baseline scenario defined by the parameters in Table~\ref{tab:Table2}, we vary customer time windows, customer demand, mobile micro-depot costs, and stationary micro-depot costs by $\pm{10}\%$ and $\pm{20}\%$. The analysis examines \glspl{kpi}, including \gls{fe} and \gls{se} routes, total cost, the number of stationary and mobile micro-depots, and direct shipments, as shown in Figure~\ref{fig:sens-table}.

As shown in the first column of Figure~\ref{fig:sens-table}, the customers' time windows have an impact on the number of \gls{fe} routes, the number of stationary micro-depots, and the number of direct shipments. As the customers' time windows decrease by 10\% and 20\%, the \gls{lsp} prefers locating more stationary micro-depots compared to the baseline solution. Moreover, the impact of direct shipment becomes more evident when time windows are reduced by 20\%. As a result, the total number of \gls{fe} routes increases compared to the baseline since the main depot directly serves more customers. A 10\% increase in time windows results in more stationary micro-depots but fewer mobile micro-depots, leaving the total number of micro-depots unchanged. However, with a 10\% decrease in time windows, we see a decrease in the number of directly shipped customers, which also reduces the number of \gls{fe} routes. A 20\% increase in the time windows does not have an impact on the \glspl{kpi}.

The second column of Figure~\ref{fig:sens-table} shows that all \glspl{kpi} correlate with variations in customer demand. As demand increases, the numbers of \gls{fe} and \gls{se} routes, total cost, and the number of micro-depots increase as expected. We also observe that the number of direct shipments increases as demand increases. This outcome is also expected and reasonable from an \gls{lsp}'s perspective. As demand increases, it will be more appropriate to serve customers directly from the main depot considering the capacity restrictions on cargo bikes. 

Analyzing the impact of mobile micro-depot costs on the \glspl{kpi} reveals minor changes in the number of \gls{fe} routes and total costs. A 20\% reduction in mobile micro-depot costs increases the number of stationary and mobile micro-depots used. Similarly, a decrease in stationary micro-depot costs leads to an increase in the number of mobile micro-depots. The \gls{lsp} prioritizes the micro-depot configuration that minimizes total costs, which may result in favoring stationary micro-depots despite lower mobile micro-depot costs or vice versa. By changing mobile or stationary micro-depot costs, we do not see a correlation with the number of \gls{se} routes.

\vspace{-0.4cm}\subsubsection{Benefit of Mobile Micro-Depots} \label{depot_configuration}
In the following, we analyze the trade-off between different micro-depot configurations in a two-echelon distribution network setting. Specifically, we create three scenarios where we allow locating only stationary micro-depots, only mobile micro-depots, and a hybrid of both stationary and mobile micro-depots. In the hybrid scenario, we run experiments for a fixed number of stationary micro-depots $\mSatellite^{f}$ between 0 and 40, optimizing additional mobile micro-depots accordingly. 

Figure ~\ref{fig:Figure8} shows the percentage change in total cost across different micro-depot configurations for both \gls{2e-lrp} and \gls{2e-lrp} with direct shipment settings. The configurations range from exclusively mobile micro-depots to exclusively stationary micro-depots, with mixed configurations in between. The configuration with 32 mobile micro-depots represents the scenario with only mobile micro-depots, while the configuration with 14 stationary micro-depots corresponds to the scenario with only stationary micro-depots.

\FloatBarrier

\textbf{Insight 3.} The number of \gls{fe} routes is sensitive to all input parameters, while the number of \gls{se} routes is only highly sensitive to changes in demand. 

\textbf{Insight 4.} The number of stationary micro-depots is sensitive to all input parameters, while the number of mobile micro-depots is only highly sensitive to changes in demand. 

\textbf{Insight 5.} Changing time windows and customer demand highly affects the number of direct shipments.

The results reveal that cost reductions occur primarily in mixed configurations, with the lowest costs achieved by opening 3 stationary and 24 mobile micro-depots for \gls{2e-lrp} (Figure~\ref{fig:Figure8_a}) and 4 stationary and 21 mobile micro-depots for \gls{2e-lrp} with direct shipment (Figure~\ref{fig:Figure8_b}). The analysis also shows that mobile micro-depots increase at a greater rate than the decrease of stationary micro-depots, indicating an over-proportional shift toward mobile micro-depots in mixed configurations. However, in scenarios dominated by mobile micro-depots, the cost-saving effect decreases, suggesting diminishing returns as the configuration departs from a balanced approach. This analysis suggests that a hybrid configuration of mobile and stationary micro-depots is the most strategic approach, offering cost efficiency.

\definecolor{RoyalBlue}{RGB}{65, 105, 225} % RGB values for Royal Blue
\definecolor{BlueGreen}{RGB}{13, 152, 186} % 

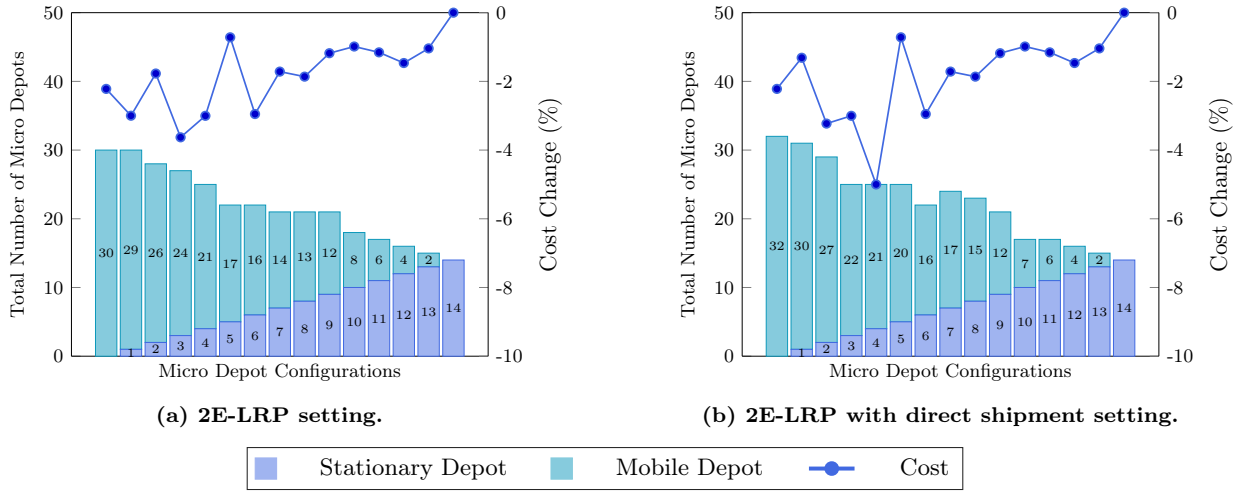
\begin{figure}[h]
    \centering
    \makebox[\textwidth][c]{
        \begin{subfigure}[t]{0.45\textwidth}
            \centering
            \scalebox{0.8}{\definecolor{RoyalBlue}{RGB}{65, 105, 225} % RGB values for Royal Blue
\definecolor{BlueGreen}{RGB}{13, 152, 186} % 
% \begin{figure}[htbp]
% \centering
\begin{tikzpicture}
  % Parameters for manual scaling
  \pgfmathsetmacro{\mincost}{50000}  % Minimum cost
  \pgfmathsetmacro{\maxcost}{62000}  % Maximum cost
  \pgfmathsetmacro{\maxtotaldepots}{45}  % Maximum total depots

  \begin{axis}[
    ybar stacked,                         % Stacked bar chart
    ymin=0,                               % Start y-axis at 0
    ymax=50,                              % Max for left y-axis
    ylabel={Total Number of Micro Depots},% Left y-axis label
    ylabel style={font=\footnotesize},
    xlabel={Micro Depot Configurations},  % X-axis label
    xlabel style={font=\footnotesize},
    xtick=\empty,                         % No x-ticks for the bar plot
    enlarge x limits=0.1,                 % Add padding on the x-axis
    ylabel near ticks,                    % Place y-axis label near ticks
    yticklabel style={font=\footnotesize},
    symbolic x coords={0, 1, 2, 3, 4, 5, 6, 7, 8, 9, 10, 11, 12, 13, 14}, % Configurations on x-axis
    axis y line*=left,                    % Left y-axis for bar chart
    bar width=10pt,                       % Width of the bars
    nodes near coords,                    % Add node labels on the bars
    nodes near coords align={center},     % Center align labels
    every node near coord/.style={        % Style for labels on the bars
      black,                              % Set label color to black
      font=\tiny                          % Set font size to tiny
    },
    every axis plot/.append style={font=\tiny}, % Font size for the plot
    xlabel style={font=\footnotesize},      % Font size for x-axis label
    ylabel style={font=\footnotesize},      % Font size for y-axis label
    ]

    % Stationary Depots (bar plot)
    \addplot+[ybar, RoyalBlue, fill=RoyalBlue!50] plot coordinates {
      (0, 0) (1, 1) (2, 2) (3, 3) (4, 4) (5, 5) (6, 6) 
      (7, 7) (8, 8) (9, 9) (10, 10) (11, 11) (12, 12) (13, 13) (14, 14)
    };

    % Mobile Depots (bar plot)
    \addplot+[ybar, BlueGreen, fill=BlueGreen!50] plot coordinates {
      (0, 30) (1, 29) (2, 26) (3, 24) (4, 21) (5, 17) (6, 16) 
      (7, 14) (8, 13) (9, 12) (10, 8) (11, 6) (12, 4) (13, 2) (14, 0)
    };

  \end{axis}

  % Add secondary y-axis for scaled costs
  \begin{axis}[
      axis y line*=right,                 % Right y-axis for line plot
      axis x line=none,                   % No x-axis for the second plot
      ylabel={Cost Change (\%)},            % Right y-axis label
      ylabel style={font=\footnotesize},
      ymin=0, ymax=50,                    % Match left y-axis range for scaling
      ylabel near ticks,                  % Place y-axis label near ticks
      yticklabel style={font=\footnotesize},
      ytick={0, 10, 20, 30, 40, 50},              % Custom ticks for right y-axis
      yticklabels={-10, -8, -6, -4, -2, 0},
      ]

      % % Line plot for costs, scaled to fit the left y-axis
      %normalized values
      %data$NormalizedCosts <- (data$Costs - min_cost) / (max_cost - min_cost) * max(data$TotalDepots) + data$TotalDepots

    % Line plot for costs, scaled to fit the left y-axis
    %normalized values converted to % change
      \addplot+[thick, color=RoyalBlue, mark=*] coordinates {
        (0,38.90476)
        (1,35)
        (2,41.14227)
        (3,31.85647)
        (4,34.98685)
        (5,46.414)
        (6,35.24666)
        (7,41.42257)
        (8,40.6887)
        (9,44.1004)
        (10,45.07344)
        (11,44.22372)
        (12,42.6707)
        (13,44.79871)
        (14,50)};

  \end{axis}
\end{tikzpicture}} % Rescale here
            \subcaption{(a) \gls{2e-lrp} setting.}
            \label{fig:Figure8_a}
        \end{subfigure}        
        \hfill
        \begin{subfigure}[t]{0.45\textwidth}
            \centering
            \scalebox{0.8}{\definecolor{RoyalBlue}{RGB}{65, 105, 225} % RGB values for Royal Blue
\definecolor{BlueGreen}{RGB}{13, 152, 186} % 
% \begin{figure}[htbp]
% \centering
\begin{tikzpicture}
  % Parameters for manual scaling
  \pgfmathsetmacro{\mincost}{50000}  % Minimum cost
  \pgfmathsetmacro{\maxcost}{64000}  % Maximum cost
  \pgfmathsetmacro{\maxtotaldepots}{45}  % Maximum total depots

  \begin{axis}[
    ybar stacked,                         % Stacked bar chart
    ymin=0,                               % Start y-axis at 0
    ymax=50,                              % Max for left y-axis
    ylabel={Total Number of Micro Depots},% Left y-axis label
    ylabel style={font=\footnotesize},
    xlabel={Micro Depot Configurations},  % X-axis label
    xlabel style={font=\footnotesize},
    xtick=\empty,                         % No x-ticks for the bar plot
    enlarge x limits=0.1,                 % Add padding on the x-axis
    ylabel near ticks,                    % Place y-axis label near ticks
    yticklabel style={font=\footnotesize},
    symbolic x coords={0, 1, 2, 3, 4, 5, 6, 7, 8, 9, 10, 11, 12, 13, 14}, % Configurations on x-axis
    axis y line*=left,                    % Left y-axis for bar chart
    bar width=10pt,                       % Width of the bars
    nodes near coords,                    % Add node labels on the bars
    nodes near coords align={center},     % Center align labels
    every node near coord/.style={        % Style for labels on the bars
      black,                              % Set label color to black
      font=\tiny                          % Set font size to tiny
    },
    every axis plot/.append style={font=\tiny}, % Font size for the plot
    xlabel style={font=\footnotesize},      % Font size for x-axis label
    ylabel style={font=\footnotesize},      % Font size for y-axis label
    ]

    % Stationary Depots (bar plot)
    \addplot+[ybar, RoyalBlue, fill=RoyalBlue!50] plot coordinates {
      (0, 0) (1, 1) (2, 2) (3, 3) (4, 4) (5, 5) (6, 6) 
      (7, 7) (8, 8) (9, 9) (10, 10) (11, 11) (12, 12) (13, 13) (14, 14)
    };

    % Mobile Depots (bar plot)
    \addplot+[ybar, BlueGreen, fill=BlueGreen!50] plot coordinates {
      (0, 32) (1, 30) (2, 27) (3, 22) (4, 21) (5, 20) (6, 16) 
      (7, 17) (8, 15) (9, 12) (10, 7) (11, 6) (12, 4) (13, 2) (14, 0)
    };

  \end{axis}

  % Add secondary y-axis for scaled costs
  \begin{axis}[
      axis y line*=right,                 % Right y-axis for line plot
      axis x line=none,                   % No x-axis for the second plot
      ylabel={Cost Change (\%)},            % Right y-axis label
      ylabel style={font=\footnotesize},
      ymin=0, ymax=50,                    % Match left y-axis range for scaling
      ylabel near ticks,                  % Place y-axis label near ticks
      yticklabel style={font=\footnotesize},
      ytick={0, 10, 20, 30, 40, 50},              % Custom ticks for right y-axis
      yticklabels={-10, -8, -6, -4, -2, 0},
      ]

      % % Line plot for costs, scaled to fit the left y-axis
      %normalized values
      %data$NormalizedCosts <- (data$Costs - min_cost) / (max_cost - min_cost) * max(data$TotalDepots) + data$TotalDepots
      %normalized values converted to % change
      \addplot+[thick, color=RoyalBlue, mark=*] coordinates {
        (0,38.90476)
        (1,43.43719)
        (2,33.85647)
        (3,34.98685)
        (4,25)
        (5,46.414)
        (6,35.24666)
        (7,41.42257)
        (8,40.6887)
        (9,44.1004)
        (10,45.07344)
        (11,44.22372)
        (12,42.6707)
        (13,44.79871)
        (14,50)};

  \end{axis}
\end{tikzpicture}} % Rescale here
            \subcaption{(b) \gls{2e-lrp} with direct shipment setting.}
            \label{fig:Figure8_b}
        \end{subfigure}
    }
    
    % Horizontal legend spanning across the figure
    \vskip 0.5em % Add vertical spacing
    \begin{tikzpicture} 
        \begin{axis}[%
            hide axis,
            xmin=10,
            xmax=50,
            ymin=0,
            ymax=0.4,
                legend style={
                draw=white!15!black,
                legend cell align=left,
                legend columns=3, % Set 2 columns
                column sep=0.5cm, % Adjust column spacing
                row sep=0.5em % Adjust row spacing
            }
            ]
            \addlegendimage{only marks, mark=square*, fill=RoyalBlue!50, draw=RoyalBlue!50, mark size=4pt}
            \addlegendentry{\footnotesize Stationary Depot};
            \addlegendimage{only marks, mark=square*, fill=BlueGreen!50, draw=BlueGreen!50, mark size=4pt}
            \addlegendentry{\footnotesize Mobile Depot};
            \addlegendimage{fill=RoyalBlue, draw=RoyalBlue, line width=1.5pt, mark=*}
            \addlegendentry{\footnotesize Cost};
        \end{axis}
    \end{tikzpicture}
    \caption{Percentage change in total cost of different micro-depot configurations in \gls{2e-lrp} and \acrshort{2e-lrp} with direct shipment settings.}
    \label{fig:Figure8}
\end{figure}

\textbf{Insight 6.} We can achieve costs savings of approximately 5\% with hybrid micro-depot configurations. 
% \FloatBarrier

We further show the changes in \glspl{kpi} across the above-mentioned scenarios for \gls{2e-lrp} and \gls{2e-lrp} with direct shipment settings in Figure~\ref{fig:Figure9}. The best results for the hybrid scenario are based on the results from Figure~\ref{fig:Figure8} for \gls{2e-lrp} and \gls{2e-lrp} with direct shipment settings. We define the reference scenario as the \gls{2e-lrp} configuration using only stationary micro-depots, and all other values are reported as percentage changes relative to this baseline. Figure~\ref{fig:Figure9_a} represents the total cost change. In the mobile-only scenario, the cost reductions are 3.07\% and 3.21\% for \gls{2e-lrp} and \gls{2e-lrp} with direct shipment settings, respectively.~Cost savings in the stationary-only and mobile-only scenarios remain marginal when direct shipment is allowed. In the \gls{2e-lrp} setting, the hybrid configuration offers no additional cost advantage over the mobile-only configuration, as their cost reductions are similar. However, in a direct shipment setting, the hybrid configuration achieves a cost reduction that is twice as much as the cost reduction from the mobile-only or hybrid configurations without direct shipment. Specifically, the hybrid configuration with direct shipment reduces total costs by 5.89\% compared to the reference scenario. Therefore, locating both stationary and mobile micro-depots while allowing direct shipment in a two-echelon distribution setting provides the cost-optimal solution.

\definecolor{RoyalBlue}{RGB}{65, 105, 225} % RGB values for Royal Blue
\definecolor{BlueGreen}{RGB}{13, 152, 186} %
\definecolor{SoftCyan}{RGB}{135, 206, 250}
\definecolor{Turquoise}{RGB}{64, 224, 208}
\definecolor{MintGreen}{RGB}{152, 251, 152}
\definecolor{SeaGreen}{RGB}{46, 139, 87}
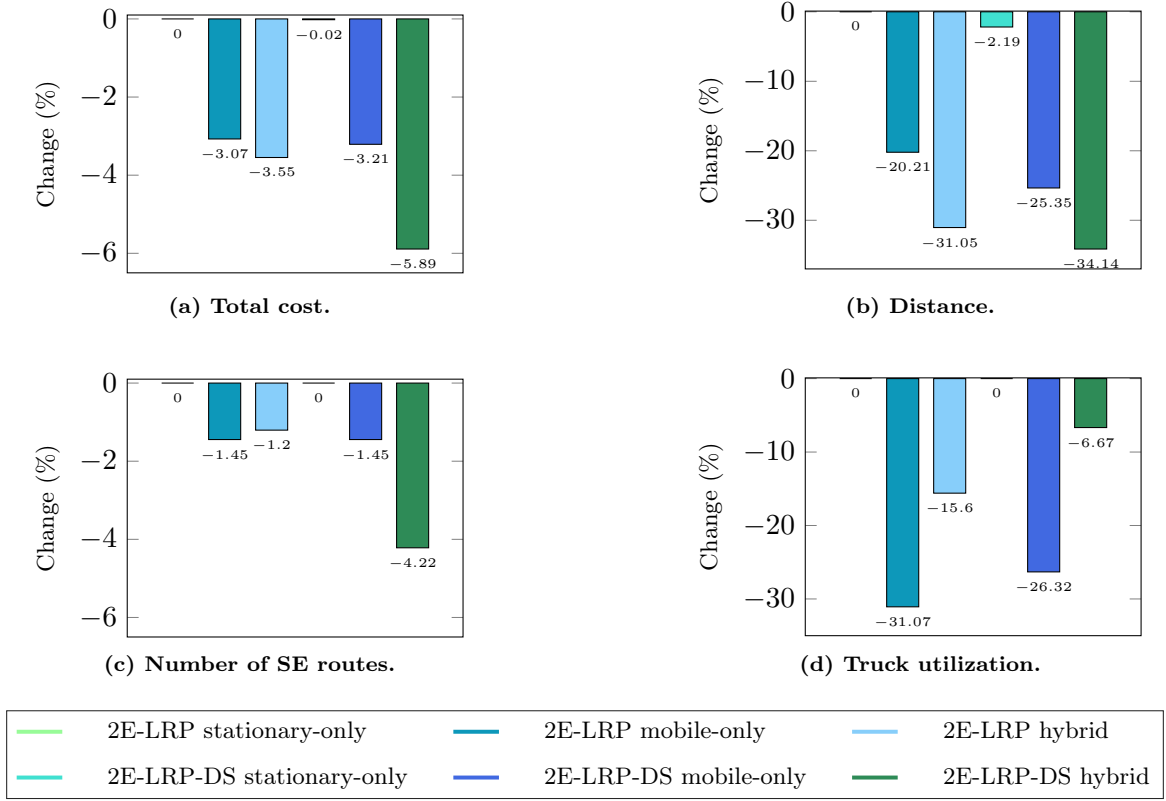
\begin{figure}[h]
    \centering
    % \scalebox{0.9}{%
    % First row
    \makebox[\textwidth][c]{
        \begin{subfigure}[t]{0.45\textwidth}
            \centering
            \definecolor{RoyalBlue}{RGB}{65, 105, 225} % RGB values for Royal Blue
\definecolor{BlueGreen}{RGB}{13, 152, 186} %
\definecolor{SoftCyan}{RGB}{135, 206, 250}
\definecolor{Turquoise}{RGB}{64, 224, 208}
\definecolor{MintGreen}{RGB}{152, 251, 152}
\definecolor{SeaGreen}{RGB}{46, 139, 87}

\begin{tikzpicture}
\begin{axis}[
    ylabel={Change (\%)},                                   % Y-axis label
    ylabel style={font=\footnotesize},
    enlarge x limits=3,                               % Add space around the x-axis
    legend style={at={(1.05,0.5)},anchor=west},         % Separate legend placement
    ybar,                                               % Bar style
    symbolic x coords={wodsstationary, wodsmobile, wodsmixed, dsstationary, dsmobile, dsmixed}, % Categories
    xtick=data,                                         % Use x-coordinates for ticks
    xtick=\empty,                                       % Hide x-axis labels
    ymin=-6.5, ymax=0.1,                                  % Y-axis range
    bar width=12pt,                                     % Adjust bar width
    width=6cm, height=5cm,                               % Adjust plot size
    nodes near coords,
    every node near coord/.append style={
        font=\fontsize{2}{3}\selectfont,                                   % Use small font
        yshift=0pt,                                        % Shift text down slightly
        anchor=north                                       % Anchor text above or below bar
    },
    nodes near coords style={/pgf/number format/.cd, fixed} % Disable scientific notation
]

% Add bar data for the dataframes with different colors
\addplot[fill=MintGreen] 
    coordinates {(wodsstationary,0)};
\addplot[fill=BlueGreen] 
    coordinates {(wodsmobile,-3.072880448)};
\addplot[fill=SoftCyan] 
    coordinates {(wodsmixed,-3.547521019)};
\addplot[fill=Turquoise] 
    coordinates {(dsstationary,-0.02)};
\addplot[fill=RoyalBlue] 
    coordinates {(dsmobile,-3.209657772)};
\addplot[fill=SeaGreen] 
    coordinates {(dsmixed,-5.890408793)};

\end{axis}
\end{tikzpicture}
            \subcaption{(a) Total cost.}
            \label{fig:Figure9_a}
        \end{subfigure}        
        % \hspace{1cm} % Space between subfigures
        \hfill
        \begin{subfigure}[t]{0.45\textwidth}
            \centering
            \definecolor{RoyalBlue}{RGB}{65, 105, 225} % RGB values for Royal Blue
\definecolor{BlueGreen}{RGB}{13, 152, 186} %
\definecolor{SoftCyan}{RGB}{135, 206, 250}
\definecolor{Turquoise}{RGB}{64, 224, 208}
\definecolor{MintGreen}{RGB}{152, 251, 152}
\definecolor{SeaGreen}{RGB}{46, 139, 87}

\begin{tikzpicture}
\begin{axis}[
    ylabel={Change (\%)},                                   % Y-axis label
    ylabel style={font=\footnotesize},
    enlarge x limits=3,                               % Add space around the x-axis
    legend style={at={(1.05,0.5)},anchor=west},         % Separate legend placement
    ybar,                                               % Bar style
    symbolic x coords={wodsstationary, wodsmobile, wodsmixed, dsstationary, dsmobile, dsmixed}, % Categories
    xtick=data,                                         % Use x-coordinates for ticks
    xtick=\empty,                                       % Hide x-axis labels
    ymin=-37, ymax=0.1,                                  % Y-axis range
    bar width=12pt,                                     % Adjust bar width
    width=6cm, height=5cm,                               % Adjust plot size
    nodes near coords,
    every node near coord/.append style={
        font=\fontsize{2}{3}\selectfont,                                   % Use small font
        yshift=0pt,                                        % Shift text down slightly
        anchor=north                                       % Anchor text above or below bar
    },
    nodes near coords style={/pgf/number format/.cd, fixed} % Disable scientific notation
]

% Add bar data for the dataframes with different colors
\addplot[fill=MintGreen] 
    coordinates {(wodsstationary,0)};
\addplot[fill=BlueGreen] 
    coordinates {(wodsmobile,-20.21196212)};
\addplot[fill=SoftCyan] 
    coordinates {(wodsmixed,-31.04641294)};
\addplot[fill=Turquoise] 
    coordinates {(dsstationary,-2.188185797)};
\addplot[fill=RoyalBlue] 
    coordinates {(dsmobile,-25.35139493)};
\addplot[fill=SeaGreen] 
    coordinates {(dsmixed,-34.13778993)};
\end{axis}
\end{tikzpicture}
            \subcaption{(b) Distance.}
            \label{fig:Figure9_b}
        \end{subfigure}
    }
    
    % Second row
    \vskip\baselineskip
    \makebox[\textwidth][c]{

        \begin{subfigure}[t]{0.45\textwidth}
            \centering
            \definecolor{RoyalBlue}{RGB}{65, 105, 225} % RGB values for Royal Blue
\definecolor{BlueGreen}{RGB}{13, 152, 186} %
\definecolor{SoftCyan}{RGB}{135, 206, 250}
\definecolor{Turquoise}{RGB}{64, 224, 208}
\definecolor{MintGreen}{RGB}{152, 251, 152}
\definecolor{SeaGreen}{RGB}{46, 139, 87}

\begin{tikzpicture}
\begin{axis}[
    ylabel={Change (\%)},                                   % Y-axis label
    ylabel style={font=\footnotesize},
    enlarge x limits=3,                               % Add space around the x-axis
    legend style={at={(1.05,0.5)},anchor=west},         % Separate legend placement
    ybar,                                               % Bar style
    symbolic x coords={wodsstationary, wodsmobile, wodsmixed, dsstationary, dsmobile, dsmixed}, % Categories
    xtick=data,                                         % Use x-coordinates for ticks
    xtick=\empty,                                       % Hide x-axis labels
    ymin=-6.5, ymax=0.1,                                  % Y-axis range
    bar width=12pt,                                     % Adjust bar width
    width=6cm, height=5cm,                               % Adjust plot size
    nodes near coords,
    every node near coord/.append style={
        font=\fontsize{2}{3}\selectfont,                                   % Use small font
        yshift=0pt,                                        % Shift text down slightly
        anchor=north                                       % Anchor text above or below bar
    },
    nodes near coords style={/pgf/number format/.cd, fixed} % Disable scientific notation
]

% Add bar data for the dataframes with different colors
\addplot[fill=MintGreen] 
    coordinates {(wodsstationary,0)};
\addplot[fill=BlueGreen] 
    coordinates {(wodsmobile,-1.445783133)};
\addplot[fill=SoftCyan] 
    coordinates {(wodsmixed,-1.204819277)};
\addplot[fill=Turquoise] 
    coordinates {(dsstationary,0)};
\addplot[fill=RoyalBlue] 
    coordinates {(dsmobile,-1.445783133)};
\addplot[fill=SeaGreen] 
    coordinates {(dsmixed,-4.21686747)};

\end{axis}
\end{tikzpicture}
            \subcaption{(c) Number of \gls{se} routes.}
            \label{fig:Figure9_c}
        \end{subfigure}

        % \hspace{1cm} % Space between subfigures
        \hfill
        
        \begin{subfigure}[t]{0.45\textwidth}
            \centering
            \definecolor{RoyalBlue}{RGB}{65, 105, 225} % RGB values for Royal Blue
\definecolor{BlueGreen}{RGB}{13, 152, 186} %
\definecolor{SoftCyan}{RGB}{135, 206, 250}
\definecolor{Turquoise}{RGB}{64, 224, 208}
\definecolor{MintGreen}{RGB}{152, 251, 152}
\definecolor{SeaGreen}{RGB}{46, 139, 87}

\begin{tikzpicture}
\begin{axis}[
    ylabel={Change (\%)},                                   % Y-axis label
    ylabel style={font=\footnotesize},
    enlarge x limits=3,                               % Add space around the x-axis
    legend style={at={(1.05,0.5)},anchor=west},         % Separate legend placement
    ybar,                                               % Bar style
    symbolic x coords={wodsstationary, wodsmobile, wodsmixed, dsstationary, dsmobile, dsmixed}, % Categories
    xtick=data,                                         % Use x-coordinates for ticks
    xtick=\empty,                                       % Hide x-axis labels
    ymin=-35, ymax=0.1,                                  % Y-axis range
    bar width=12pt,                                     % Adjust bar width
    width=6cm, height=5cm,                               % Adjust plot size
    nodes near coords,
    every node near coord/.append style={
        font=\fontsize{2}{3}\selectfont,                                   % Use small font
        yshift=0pt,                                        % Shift text down slightly
        anchor=north                                       % Anchor text above or below bar
    },
    nodes near coords style={/pgf/number format/.cd, fixed} % Disable scientific notation
]

% Add bar data for the dataframes with different colors
\addplot[fill=MintGreen] 
    coordinates {(wodsstationary,0)};
\addplot[fill=BlueGreen] 
    coordinates {(wodsmobile,-31.07177033)};
\addplot[fill=SoftCyan] 
    coordinates {(wodsmixed,-15.59735509)};
\addplot[fill=Turquoise] 
    coordinates {(dsstationary,0)};
\addplot[fill=RoyalBlue] 
    coordinates {(dsmobile,-26.31578947)};
\addplot[fill=SeaGreen] 
    coordinates {(dsmixed,-6.666666667)};

\end{axis}
\end{tikzpicture}
            \subcaption{(d) Truck utilization.}
            \label{fig:Figure9_d}
        \end{subfigure}
    }
    
    % Legend
    \vskip 1em % Add vertical spacing
    \begin{tikzpicture} 
        \begin{axis}[
            hide axis,
            xmin=10,
            xmax=50,
            ymin=0,
            ymax=0.4,
            legend style={
                draw=white!15!black,
                legend cell align=left,
                legend columns=3, % Set 2 columns
                column sep=0.5cm, % Adjust column spacing
                row sep=0.5em % Adjust row spacing
            }
        ]

            \addlegendimage{fill=MintGreen, draw=MintGreen, line width=1.5pt} 
            \addlegendentry{\footnotesize 2E-LRP stationary-only};
            \addlegendimage{fill=BlueGreen, draw=BlueGreen, line width=1.5pt}
            \addlegendentry{\footnotesize 2E-LRP mobile-only};
            \addlegendimage{fill=SoftCyan, draw=SoftCyan, line width=1.5pt} 
            \addlegendentry{\footnotesize 2E-LRP hybrid};
            \addlegendimage{fill=Turquoise, draw=Turquoise, line width=1.5pt}
            \addlegendentry{\footnotesize 2E-LRP-DS stationary-only};
            \addlegendimage{fill=RoyalBlue, draw=RoyalBlue, line width=1.5pt} 
            \addlegendentry{\footnotesize 2E-LRP-DS mobile-only};
            \addlegendimage{fill=SeaGreen, draw=SeaGreen, line width=1.5pt}
            \addlegendentry{\footnotesize 2E-LRP-DS hybrid};
        \end{axis}
    \end{tikzpicture}
    % }  
    \caption{Comparison of micro-depot configurations in \gls{2e-lrp} and \acrshort{2e-lrp} with direct shipment (DS) settings.}
    \label{fig:Figure9}
\end{figure}

\textbf{Insight 7.} A hybrid micro-depot configuration results in a 5.89\% reduction in total costs, 34.14\% in distances and 4.22\% in the number of \gls{se} routes compared to stationary-only configurations in a \gls{2e-lrp} setting. 

As illustrated in Figure~\ref{fig:Figure9_b}, the reductions in the total distance for both mobile-only and hybrid configurations are more prominent compared to the stationary-only configuration. By allowing direct shipment, distances are further reduced for all three scenarios. The hybrid configuration with direct shipment achieves the highest reduction, lowering total distance by 34.14\%. Figure~\ref{fig:Figure9_c} shows the percentage change in the number of \gls{se} routes, where mobile-only and hybrid configurations demonstrate better cargo bike utilization. Notably, the hybrid configuration with direct shipment reduces \gls{se} routes by 4.22\% compared to stationary-only configurations, offering superior route consolidation and making it the preferred micro-depot configuration. Finally, Figure~\ref{fig:Figure9_d} shows the average truck utilization in the \gls{fe}. In contrast to other \glspl{kpi}, stationary-only configurations generate the maximum utilized trucks. Mobile-only configurations, on the other hand, yield the lowest truck utilization rates. Hybrid configurations exhibit a 15.6\% decrease in truck utilization for the \gls{2e-lrp} setting. By allowing direct shipment, truck utilizations decrease by 6.67\% compared to the stationary-only configurations.

The hybrid configuration in \gls{2e-lrp} with direct shipment demonstrates a balance between cost efficiency, distance traveled, and adaptability to dynamic changes in demand locations. Although it may not achieve the highest truck utilization, its significant improvements in other \glspl{kpi} establish it as the most comprehensive and effective solution for \glspl{lsp}.
\FloatBarrier
\section{Conclusion} \label{conclusion}
In this paper, we proposed an extension of \gls{2e-lrp}, which considers the concurrent location of stationary and mobile micro-depots and allows direct shipment from the main depot to customers in city logistics. To address this problem, we developed an \gls{alns}-based solution approach, integrating a set cover problem to identify optimal network configurations and micro-depot locations. We evaluated our \gls{alns} on Prodhon's and Nguyen's \gls{2e-lrp} benchmark instances, achieving high-quality solutions with short computational times.
We introduced a decomposition-based cluster-first-route-second approach to tackle large-scale instances, which outperformed baseline algorithms by running 15 times faster while maintaining good solution quality. Using case study instances based on the city of Munich, our comparative analysis demonstrated that allowing direct shipment in the two-echelon distribution setting reduces total costs and emissions by 4.7\% and 11\%, respectively, while increasing truck utilization by 42\%. Our sensitivity analysis indicated that direct shipment is highly sensitive to input parameters. Additionally, we analyzed trade-offs between different micro-depot configurations, showing that hybrid configurations combining stationary and mobile depots achieve lower costs compared to using a single-depot type. Specifically, a hybrid configuration with direct shipment reduced total costs by 5.9\% compared to the traditional \gls{2e-lrp} setting.
Our study focuses solely on cost minimization. Future research could systematically investigate the environmental impact of the distribution system, particularly with respect to emissions. Additionally, learning-based approaches could be explored to identify efficient solution configurations.

\FloatBarrier
%
%% finally the bib file
\singlespacing{
%\footnotesize
\bibliographystyle{model5-names}%\biboptions{authoryear}
\bibliography{./alns}} % if more than one, comma separated
\newpage
%% and the appendices
\onehalfspacing
\begin{appendices}
	\normalsize
	\appendix       
        \section{Mathematical Formulation of 2E-LRP with Mobile Depots and Direct Shipment}\label{math-model}
            In this section, we define the mathematical formulation of \gls{2e-lrp} with mobile depots and direct shipment.
            
            \renewcommand\baselinestretch{1}
            \textbf{Sets} 
            
                \begin{tabular}{rl}
                $\mSetVertices$ & set of vertices \\
                 $\mSetArcs{}$ &  set of arcs \\
                    $\mSetHubLocations$ & set of micro-depots \\
                    $\mSetCustomers$ & set of customers \\
                    $\mSetTypeVehicle{f}$ & set of first echelon vehicles \\
                    $\mSetTypeVehicle{s}$ & set of second echelon vehicle \\
                    \end{tabular}\\
            
            \textbf{Parameters} 
            
                \begin{tabular}{p{1cm}p{15cm}}
                    \mCostLocation{}{\mNodeFirst} & micro-depot cost for using location $\mNodeFirst$ \\
                    \mCostArc{\mVehicleSecond}{\mNodeFirst \mNodeSecond}& travel cost on arc $(\mNodeFirst, \mNodeSecond)$ for vehicle $\mVehicleSecond$ \\
                    \mFixedCost{\mVehicleSecond} & fixed cost of vehicle $\mVehicleSecond$ \\
                    \mFreightCapacity{}{\mVehicleSecond} & capacity of vehicle $\mVehicleSecond$ \\
                    \mHubCapacity{\mSatellite}{} & capacity of micro-depot $\mSatellite$ \\
                    \mDemand{}{\mNodeFirst} & demand of customer $\mNodeFirst$ \\
                    $[\mEarliestBOS{\mNodeFirst},\mLatestBOS{\mNodeFirst}]$ & begin of service time window at location $\mNodeFirst$ \\
                    \mTravelTimeOfArc{f}{\mNodeFirst}{\mNodeSecond\mVehicleFirst} & travel time of vehicle $\mVehicleFirst$ on arc $(\mNodeFirst, \mNodeSecond)$ when departing from $\mNodeFirst$ on first echelon\\
                    \mTravelTimeOfArc{s}{\mNodeFirst}{\mNodeSecond\mVehicleFirst} & travel time of vehicle $\mVehicleSecond$ on arc $(\mNodeFirst, \mNodeSecond)$ when departing from $\mNodeFirst$ on second echelon\\
                \end{tabular}\\
         
            \textbf{Decision Variables}
            
                \begin{tabular}{p{1.5cm}p{12.5cm}}
                    \mDecisionHubLocation{}{\mSatellite} & decision on whether to use a micro-depot $\mSatellite$\\
                    \mDecisionArc{f}{\mNodeFirst\mNodeSecond\mVehicleFirst} & decision on whether vehicle $\mVehicleFirst$ is traveling on arc $(\mNodeFirst,\mNodeSecond)$ on first echelon\\
                    \mDecisionArc{s}{\mNodeFirst\mNodeSecond\mTypeVehicle \mSatellite} & decision on whether vehicle $\mTypeVehicle$ is traveling on arc $(\mNodeFirst,\mNodeSecond)$ served from micro-depot $\mSatellite$ on second echelon\\
                    \mDecisionVehicle{f}{\mVehicleFirst\mSatellite} & decision on whether vehicle $\mVehicleFirst$ is used and visit micro-depot $\mSatellite$ on first echelon\\
                    \mDecisionVehicle{s}{\mTypeVehicle \mSatellite} & decision on whether vehicle $\mTypeVehicle$ is used and served from micro-depot $\mSatellite$ on second echelon\\
                    \mBeginOfService{f}{\mNodeFirst\mTypeVehicle} & begin of service by vehicle $\mTypeVehicle$ at micro-depot $\mNodeFirst$ on first echelon\\
                    \mBeginOfService{s}{\mNodeFirst\mTypeVehicle} & begin of service by vehicle $\mTypeVehicle$ at customer $\mNodeFirst$ on second echelon\\
                    \mArrival{end}{\mTypeVehicle} & the time vehicle $\mTypeVehicle$ arrives at the main depot \\
                    \mDeparture{s}{\mNodeFirst\mTypeVehicle} & departure time of vehicle $\mTypeVehicle$ from location $\mNodeFirst$ on second echelon\\
                    \mArrival{s}{\mNodeFirst\mTypeVehicle} & the time vehicle $\mTypeVehicle$ arrives at location $\mNodeFirst$ on second echelon\\
                \end{tabular}\\

            \textbf{Objective Function}
            \begin{multline}
            \min\SumSet{\mSatellite\in\mSetHubLocations}\mCostLocation{}{\mSatellite} \mDecisionHubLocation{}{\mSatellite} + \SumSet{\mVehicleFirst\in
            	\mSetTypeVehicle{f}}~~\SumSet{(\mNodeFirst,\mNodeSecond)\in\mSetArcs{f}} \mCostArc{}{\mNodeFirst\mNodeSecond} \mDecisionArc{}{\mNodeFirst\mNodeSecond\mVehicleFirst} + \SumSet{\mTypeVehicle\in
            	\mSetTypeVehicle{s}}~~\SumSet{(\mNodeFirst,\mNodeSecond)\in\mSetArcs{s}} \mCostArc{}{\mNodeFirst\mNodeSecond} \mDecisionArc{}{\mNodeFirst\mNodeSecond\mTypeVehicle }{}
            	 + \SumSet{\mVehicleFirst\in\mSetTypeVehicle{f}} ~\SumSet{\mSatellite\in\mSetHubLocations} \mFixedCost{\mVehicleFirst} \mDecisionVehicle{}{\mVehicleFirst\mSatellite}
            	 	 + \SumSet{\mTypeVehicle\in\mSetTypeVehicle{s}} ~ \SumSet{\mSatellite\in\mSetHubLocations} \mFixedCost{\mTypeVehicle} \mDecisionVehicle{}{\mTypeVehicle\mSatellite} \\
            \label{obj: cost}
            \end{multline}
            
            \textbf{Routing Constraints for the 1\textsuperscript{st} Echelon}
            %assignment constraint e1
            \begin{multline}
            \quad\quad\mDecisionArc{f}{\mNodeFirst\mSatellite\mVehicleFirst} \leq \mDecisionHubLocation{}{\mSatellite}  \hfill \forall \mNodeFirst\in \mSetVertices^{f}, \mSatellite \in \mSetHubLocations, \mVehicleFirst \in \mSetTypeVehicle{f}
            \label{constr:assignment_constr_E1}
            \end{multline}
            
            %vehicle departs from depot constraint e1
            \begin{multline}
            \quad\quad\mDecisionArc{f}{\mNodeFirst\mSatellite\mVehicleFirst} \leq \mDecisionArc{f}{0\mSatellite\mVehicleFirst}  \hfill \forall \mNodeFirst\in \mSetVertices^{f}, \mSatellite \in \mSetHubLocations, \mVehicleFirst \in \mSetTypeVehicle{f}
            \label{constr:vehicle_departs_from_depot_constr_E1}
            \end{multline}
            
            % everything that goes in, goes out (flow balance constraint)
            \begin{multline}
            \quad\quad 	\SumSet{\mNodeSecond \in \mSetVertices^{f}} \mDecisionArc{f}{\mNodeSecond\mNodeFirst\mVehicleFirst} - \SumSet{\mNodeSecond \in \mSetVertices^{f}} \mDecisionArc{f}{\mNodeFirst\mNodeSecond\mVehicleFirst} = 0 \hfill \forall \mNodeFirst \in \mSetVertices^{f}, \mVehicleFirst \in \mSetTypeVehicle{f}
            \label{constr:flow_balance_E1}
            \end{multline}
            
            %depot_usage_constraint
            \begin{multline}
            \quad\quad\SumSet{\mSatellite \in \mSetHubLocations}~~\SumSet{\mVehicleFirst \in \mSetTypeVehicle{f}}\mDecisionArc{f}{0\mSatellite\mVehicleFirst} \geq 1 \hfill
            \label{constr:vehicle_depot_constr_E1}
            \end{multline}
            
            %vehicle_usage_constraint_E1
            \begin{multline}
            \quad\quad\SumSet{\mNodeSecond \in \mSetVertices^{f}} \mDecisionArc{f}{\mNodeFirst\mNodeSecond\mVehicleFirst} \leq  \SumSet{\mSatellite \in \mSetHubLocations}\mDecisionVehicle{f}{\mVehicleFirst\mSatellite}  \hfill \forall \mNodeFirst \in \mSetVertices^{f}, \mVehicleFirst \in \mSetTypeVehicle{f} 
            \label{constr:vehicle_usage_constraint_E1}
            \end{multline}

            \textbf{Routing Constraints for the 2\textsuperscript{nd} Echelon}
            
            %assignment constraint e2
            \begin{multline}
            \quad\quad\mDecisionArc{s}{\mSatellite\mNodeSecond\mVehicleSecond\mSatellite} \leq \mDecisionHubLocation{}{\mSatellite} \hfill \forall \mNodeSecond\in \mSetVertices^{s}, \mSatellite \in \mSetHubLocations, \mVehicleSecond \in \mSetTypeVehicle{s}
            \label{constr:assignment_constr_E2}
            \end{multline}
            
            %vehicle_departs_from_depot_E2
            \begin{multline}
            \quad\quad\SumSet{\mSatellite \in \mSetHubLocations}~~\SumSet{\mNodeSecond \in \mSetCustomers}\mDecisionArc{s}{\mSatellite\mNodeSecond\mVehicleSecond\mSatellite} \leq 1  \hfill \forall \mVehicleSecond \in \mSetTypeVehicle{s}
            \label{constr:vehicle_departs_from_depot_constr_E2}
            \end{multline}
            
            % everything that goes in, goes out (flow balance constraint)
            \begin{multline}
            \quad\quad 	\SumSet{\mNodeSecond \in \mSetCustomers} \mDecisionArc{s}{\mSatellite\mNodeSecond\mVehicleSecond\mSatellite} - \SumSet{\mNodeSecond \in \mSetCustomers} \mDecisionArc{s}{\mNodeSecond\mSatellite\mVehicleSecond\mSatellite} = 0 \hfill \forall \mSatellite \in\mSetHubLocations, \mTypeVehicle \in \mSetTypeVehicle{s}
            \label{constr:flow_balance_E2_1}
            \end{multline}
            
            % everything that goes in, goes out (flow balance constraint)
            \begin{multline}
            \quad\quad 	\SumSet{\mNodeSecond \in \mSetVertices^{s}} \mDecisionArc{s}{\mNodeSecond\mNodeFirst\mVehicleSecond\mSatellite} - \SumSet{\mNodeSecond \in \mSetVertices^{s}} \mDecisionArc{s}{\mNodeFirst\mNodeSecond\mVehicleSecond\mSatellite} = 0 \hfill \forall \mNodeFirst \in\mSetCustomers, \mSatellite \in\mSetHubLocations, \mTypeVehicle \in \mSetTypeVehicle{s}
            \label{constr:flow_balance_E2_2}
            \end{multline}
            
            %vehicle_usage_constraint_E2
            \begin{multline}
            \quad\quad\SumSet{\mSatellite \in \mSetHubLocations}~~\SumSet{\mNodeSecond \in \mSetVertices^{s}} \mDecisionArc{s}{\mNodeFirst\mNodeSecond\mVehicleSecond\mSatellite} \leq \SumSet{\mSatellite \in \mSetHubLocations}\mDecisionVehicle{s}{\mVehicleSecond\mSatellite}  \hfill \forall i\in \mSetVertices^{s}, \mVehicleSecond \in \mSetTypeVehicle{s}
            \label{constr:vehicle_usage_constraint_E2}
            \end{multline}
            
            %customer_visited_echelon
            \begin{multline}
            \quad\quad\SumSet{\mNodeFirst \in \mSetVertices^{f}}~~\SumSet{\mVehicleFirst\in\mSetTypeVehicle{f}} \mDecisionArc{f}{\mNodeFirst\mNodeSecond\mVehicleFirst} + \SumSet{\mNodeFirst \in \mSetVertices^{s}}~~\SumSet{\mVehicleSecond\in\mSetTypeVehicle{s}}~~\SumSet{\mSatellite \in \mSetHubLocations} \mDecisionArc{s}{\mNodeFirst\mNodeSecond\mVehicleSecond\mSatellite} = 1 \hfill \forall \mNodeSecond \in \mSetCustomers
            \label{constr:customer_visited_echelon}
            \end{multline}
            
            \textbf{Capacity Constraints}
            
            %cap_constr_E1
            \begin{multline}
            \quad\quad \SumSet{\mNodeFirst \in \mSetVertices^{f}}~~\SumSet{\mNodeSecond \in \mSetVertices^{f}} \mDemand{}{\mNodeFirst} \mDecisionArc{f} {\mNodeFirst\mNodeSecond\mVehicleFirst} \leq \mFreightCapacity{}{\mVehicleFirst} \SumSet{\mSatellite \in \mSetHubLocations}\mDecisionVehicle{f}{\mVehicleFirst\mSatellite} \hfill \forall \mVehicleFirst\in\mSetTypeVehicle{f} 
            \label{constr:cap_constr_E1}
            \end{multline}
            
            %cap_constr_E2
            \begin{multline}
            \quad\quad \SumSet{\mNodeFirst \in \mSetVertices^{s}}~~\SumSet{\mNodeSecond \in \mSetVertices^{s}} \mDemand{}{\mNodeFirst} \mDecisionArc{s}{\mNodeFirst\mNodeSecond\mVehicleSecond\mSatellite} \leq \mFreightCapacity{}{\mTypeVehicle} \SumSet{\mSatellite \in \mSetHubLocations}\mDecisionVehicle{s}{\mTypeVehicle \mSatellite} \hfill \forall \mVehicleSecond\in\mSetTypeVehicle{s}
            \label{constr:cap_constr_E2}
            \end{multline}
            
            %hub_capacity_constr
            \begin{multline}
            \quad\quad \SumSet{\mNodeFirst \in \mSetVertices^{s}}~~\SumSet{\mNodeSecond \in \mSetVertices^{s}}~~ \SumSet{\mVehicleSecond\in\mSetTypeVehicle{s}}\mDemand{}{\mNodeSecond} \mDecisionArc{s}{\mNodeFirst\mNodeSecond\mVehicleSecond\mSatellite} \leq \mHubCapacity{}{\mSatellite} \SumSet{\mVehicleSecond\in\mSetTypeVehicle{s}}\mDecisionVehicle{s}{\mTypeVehicle \mSatellite} \hfill \forall \mSatellite \in \mSetHubLocations
            \label{constr:hub_capacity_constr}
            \end{multline}
            
            %sum_vehicle_capacity
            \begin{multline}
            \quad\quad \SumSet{\mVehicleSecond\in\mSetTypeVehicle{s}} \mFreightCapacity{}{\mTypeVehicle} \mDecisionVehicle{s}{\mTypeVehicle \mSatellite} \leq 
            \mHubCapacity{}{\mSatellite} \mDecisionHubLocation{}{\mSatellite} \hfill
            \forall \mSatellite \in \mSetHubLocations 
            \label{constr:sum_vehicle_capacity}
            \end{multline}
            
            \textbf{Time Calculation Constraints for the 1\textsuperscript{st} Echelon}
            %departure_time_from_depot
            \begin{multline}
            \quad\quad \mBeginOfService{f}{0\mVehicleFirst} \geq \mEarliestBOS{\mVehicleFirst} - M (1 - \SumSet{\mSatellite\in\mSetHubLocations} \mDecisionArc{f}{0\mSatellite\mVehicleFirst}) \hfill
            \forall \mVehicleFirst\in\mSetTypeVehicle{f}
            \label{constr:departure_time_from_depot}
            \end{multline}
            
            %arrival_time_to_micro-depot
            \begin{multline}
            \quad\quad \mBeginOfService{f}{\mSatellite\mVehicleFirst} \geq \mBeginOfService{f}{\mNodeFirst\mVehicleFirst} +  \mTravelTimeOfArc{\mVehicleFirst}{\mNodeFirst}{\mSatellite} - M (1 - \mDecisionArc{f}{\mNodeFirst\mSatellite\mVehicleFirst}) 
            \hfill \forall \mVehicleFirst\in\mSetTypeVehicle{f}, \mNodeFirst \in \mSetVertices^{f}, \mSatellite \in \mSetHubLocations
            \label{constr:arrival_time_to_micro-depot}
            \end{multline}
            
            %arrival_time_to_depot_1
            \begin{multline}
            \quad\quad \mBeginOfService{f}{\mSatellite\mVehicleFirst} + \mTravelTimeOfArc{\mVehicleFirst}{\mSatellite}{0} - M (1 - \mDecisionArc{f}{\mSatellite0\mVehicleFirst}) \leq \mArrival{end}{\mVehicleFirst} \hfill \forall \mVehicleFirst\in\mSetTypeVehicle{f}, \mNodeFirst \in \mSetVertices^{f}, \mSatellite \in \mSetHubLocations
            \label{constr:arrival_time_to_depot_1}
            \end{multline}
            
            %arrival_time_to_depot_2
            \begin{multline}
            \quad\quad \mLatestBOS{\mVehicleFirst} + M (1 - \SumSet{\mSatellite\in\mSetHubLocations} \mDecisionArc{f}{\mSatellite0\mVehicleFirst})\geq \mArrival{end}{\mVehicleFirst} \hfill \forall \mVehicleFirst\in\mSetTypeVehicle{f}
            \label{constr:arrival_time_to_depot_2}
            \end{multline}
            
            %arrival_time_to_micro-depot_2
            \begin{multline}
            \quad\quad \mBeginOfService{f}{\mSatellite\mVehicleFirst} \leq M \SumSet{\mNodeFirst \in \mSetVertices^{f}} \mDecisionArc{f}{\mNodeFirst\mSatellite\mVehicleFirst} \hfill \forall \mSatellite\in\mSetHubLocations, \mVehicleFirst\in\mSetTypeVehicle{f}
            \label{constr:arrival_time_to_micro-depot_2}
            \end{multline}
            
            %arrival_time_to_depot_bound
            \begin{multline}
            \quad\quad \mArrival{end}{\mVehicleFirst} \leq M \SumSet{\mSatellite\in\mSetHubLocations} \mDecisionArc{f}{\mSatellite0\mVehicleFirst} \hfill \forall \mVehicleFirst\in\mSetTypeVehicle{f}
            \label{constr:arrival_time_to_depot_bound}
            \end{multline}

            \textbf{Time Calculation Constraints for the 2\textsuperscript{nd} Echelon}
            
            %%%%%Echelon 2%%%%%%%%
            %arrival_time_after_micro-depot
            \begin{multline}
            \quad\quad \mBeginOfService{s}{\mNodeSecond\mVehicleSecond} \geq \mDeparture{s}{\mSatellite\mTypeVehicle} + \mTravelTimeOfArc{\mTypeVehicle}{\mSatellite}{\mNodeSecond} - M (1 - \mDecisionArc{s}{\mSatellite j\mTypeVehicle \mSatellite}) \hfill \forall \mNodeSecond \in \mSetCustomers, \mSatellite\in\mSetHubLocations, \mTypeVehicle\in\mSetTypeVehicle{s}
            \label{constr:arrival_time_after_micro-depot_e2}
            \end{multline}
            
            %arrival_time_to_customer
            \begin{multline}
            \quad\quad \mBeginOfService{s}{\mNodeSecond\mVehicleSecond} \geq  \mBeginOfService{s}{\mNodeFirst\mVehicleSecond} + \mTravelTimeOfArc{\mTypeVehicle}{\mNodeFirst}{\mNodeSecond} - M (1 - \mDecisionArc{s}{\mNodeFirst\mNodeSecond\mTypeVehicle \mSatellite}) \hfill \forall i , \mNodeSecond \in \mSetCustomers, \mTypeVehicle\in\mSetTypeVehicle{s}, \mSatellite\in\mSetHubLocations
            \label{constr:arrival_time_to_customer_e2}
            \end{multline}
            
            %arrival_time_to_micro-depot
            \begin{multline}
            \quad\quad \mArrival{s}{\mSatellite\mTypeVehicle} \geq  \mBeginOfService{s}{\mNodeFirst\mVehicleSecond} + \mTravelTimeOfArc{\mTypeVehicle}{\mNodeFirst}{\mSatellite} - M (1 - \mDecisionArc{s}{\mNodeFirst\mSatellite\mTypeVehicle \mSatellite}) \hfill \forall \mNodeFirst \in \mSetCustomers, \mSatellite\in\mSetHubLocations, \mTypeVehicle\in\mSetTypeVehicle{s}
            \label{constr:arrival_time_to_micro-depot_e2}
            \end{multline}
            
            %arrival_time_to_micro-depot_2
            \begin{multline}
            \quad\quad \mArrival{s}{\mSatellite\mTypeVehicle} \leq  \mLatestBOS{\mSatellite} + M (1 - \SumSet{\mNodeSecond\in \mSetVertices^{s}}\mDecisionArc{s}{\mNodeSecond\mSatellite\mTypeVehicle \mSatellite}) \hfill \forall  \mSatellite\in\mSetHubLocations, \mTypeVehicle\in\mSetTypeVehicle{s}
            \label{constr:arrival_time_to_micro-depot_2_e2}
            \end{multline}
            
            %departure_time_from_micro-depot
            \begin{multline}
            \quad\quad \mDeparture{s}{\mSatellite\mTypeVehicle} \leq M \SumSet{\mNodeFirst\in\mSetCustomers} \mDecisionArc{s}{\mSatellite\mNodeFirst\mVehicleSecond\mSatellite} \hfill \forall \mSatellite\in\mSetHubLocations, \mTypeVehicle\in\mSetTypeVehicle{s}
            \label{constr:departure_time_from_micro-depot_e2}
            \end{multline}
            
            %return_time_to_micro-depot
            \begin{multline}
            \quad\quad \mArrival{s}{\mSatellite\mTypeVehicle} \leq   M \SumSet{\mNodeFirst\in\mSetCustomers} \mDecisionArc{s}{\mSatellite\mNodeFirst\mVehicleSecond\mSatellite} \hfill \forall \mSatellite\in\mSetHubLocations, \mTypeVehicle\in\mSetTypeVehicle{s}
            \label{constr:return_time_to_micro-depot_e2}
            \end{multline}

            % starting service inside time window
            \begin{multline}
            \quad\quad \mEarliestBOS{\mNodeFirst} \leq \mBeginOfService{\mTypeVehicle}{\mNodeFirst} \leq \mLatestBOS{\mNodeFirst}\hfill \forall \mNodeFirst\in \mSetCustomers, \forall \mTypeVehicle\in\mSetTypeVehicle{s}
            \label{constr:time_windows}
            \end{multline}
            
            \textbf{Integrality and Nonnegativity Constraints}
            % domains
            \begin{multline}
            \quad\quad \mDecisionArc{f}{\mNodeFirst\mNodeSecond\mVehicleFirst}\in\{0,1\} \hfill \forall \mNodeFirst,\mNodeSecond \in \mSetVertices^{f}, \mVehicleFirst \in\mSetTypeVehicle{f}
            \label{constr: arc_decision_e1}
            \end{multline}
            
            \begin{multline}
            \quad\quad \mDecisionArc{s}{\mNodeFirst\mNodeSecond\mVehicleSecond\mSatellite}\in\{0,1\} \hfill \forall \mNodeFirst,\mNodeSecond \in \mSetVertices^{s}, \mVehicleSecond \in\mSetTypeVehicle{s}, \mSatellite\in\mSetHubLocations
            \end{multline}
            
            \begin{multline}
            \quad\quad \mDecisionVehicle{f}{\mVehicleFirst\mSatellite} \in\{0,1\} \hfill  \forall \mVehicleFirst \in\mSetTypeVehicle{f},\mSatellite\in\mSetHubLocations
            \end{multline}
            
            \begin{multline}
            \quad\quad \mDecisionVehicle{s}{\mVehicleSecond\mSatellite} \in\{0,1\} \hfill  \forall \mVehicleSecond \in\mSetTypeVehicle{s}, \mSatellite\in\mSetHubLocations
            \end{multline}
            
            \begin{multline}
            \quad\quad \mResidualCargo{}{\mVehicleFirst}{\mNodeFirst} \geq 0 \hfill \forall \mNodeFirst \in \mSetVertices^{f}, \mVehicleFirst \in\mSetTypeVehicle{f}
            \end{multline}
            
            \begin{multline}
            \quad\quad\mBeginOfService{f}{\mNodeFirst\mVehicleFirst}\geq 0,\;\;  \mArrival{f}{\mNodeFirst\mVehicleFirst}\geq 0 \hfill \forall \mNodeFirst\in \mSetVertices^{f}, \forall \mVehicleFirst\in\mSetTypeVehicle{f}
            \label{constr:domains}
            \end{multline}
            
            \begin{multline}
            \quad\quad\mBeginOfService{s}{\mNodeFirst\mVehicleSecond}\geq 0, \;\; \mArrival{s}{\mNodeFirst\mVehicleSecond}\geq 0,\;\; \mDeparture{s}{\mNodeFirst\mTypeVehicle} \geq 0 \hfill \forall \mNodeFirst\in \mSetCustomers, \forall \mVehicleSecond\in\mSetTypeVehicle{s}
            \label{constr:domains}
            \end{multline}
            
            \begin{multline}
            \quad\quad\mDecisionHubLocation{}{\mSatellite}\in\{0,1\} \hfill \forall \mSatellite\in\mSetHubLocations
            \label{constr: domains}
            \end{multline}

            The objective function \ref{obj: cost} consists of the total fixed costs of opening the micro-depots, the traveling costs in both echelons and the fixed costs of the vehicles used in both echelons. Constraint~\ref{constr:assignment_constr_E1} ensures that a vehicle can visit a micro-depot only if the micro-depot is open. Constraint~\ref{constr:vehicle_departs_from_depot_constr_E1} ensures that a vehicle can visit a micro-depot only if it departs from the main depot. Constraint~\ref{constr:flow_balance_E1} is the flow balance constraint on the first echelon. Constraint~\ref{constr:vehicle_depot_constr_E1} ensures that there must be at least an outgoing arc from the main depot. Constraint~\ref{constr:vehicle_usage_constraint_E1} ensures that a vehicle can travel on arc $(i,j)$ only if the vehicle is used.
            
            Constraints \ref{constr:assignment_constr_E2} - \ref{constr:customer_visited_echelon} define the routing constraints for the second echelon. \ref{constr:assignment_constr_E2} ensures that a vehicle can visit a micro-depot only if the micro-depot is open. \ref{constr:vehicle_departs_from_depot_constr_E2} ensures that a vehicle cannot depart from more than one micro-depot. \ref{constr:flow_balance_E2_1} and \ref{constr:flow_balance_E2_2} are the flow balance constraints on the second echelon. \ref{constr:vehicle_usage_constraint_E2}  ensures that a vehicle can travel on arc $(i,j)$ only if the vehicle is used. \ref{constr:customer_visited_echelon} ensures that a customer is served either in the first or second echelon.
            
            Constraints \ref{constr:cap_constr_E1} - \ref{constr:sum_vehicle_capacity} define the capacity constraints for vehicles at both first and second echelon and micro-depots.
            
            Constraints \ref{constr:departure_time_from_depot} - \ref{constr:arrival_time_to_depot_bound} define the time calculations on the first echelon. Constraint \ref{constr:departure_time_from_depot} ensures that the departure time from the depot should be at least the earliest working time of the vehicle. Constraint \ref{constr:arrival_time_to_micro-depot} is the sub-tour elimination constraint. Constraints \ref{constr:arrival_time_to_depot_1} and \ref{constr:arrival_time_to_depot_2} ensure that the arrival time at the depot should be after serving the last node in the route and the latest working time of the vehicle. Constraint \ref{constr:arrival_time_to_micro-depot_2} ensures that the arrival time at the micro-depot can be positive only if the micro-depot is visited. Finally, constraint \ref{constr:arrival_time_to_depot_bound} ensures that the return time to main depot gets value only if there is an outgoing arc. 
            
            Constraints \ref{constr:arrival_time_after_micro-depot_e2}–\ref{constr:time_windows} define the time calculations for the second echelon. Constraint \ref{constr:arrival_time_after_micro-depot_e2} ensures that a vehicle arrives at a customer after departing from the micro-depot. Constraint \ref{constr:arrival_time_to_customer_e2} eliminates sub-tours. Constraints \ref{constr:arrival_time_to_micro-depot_e2} and \ref{constr:arrival_time_to_micro-depot_2_e2} calculate the arrival time at the micro-depot after completing the route. Constraint \ref{constr:departure_time_from_micro-depot_e2} ensures that the departure time from the micro-depot is set only when there is an outgoing arc, while Constraint \ref{constr:return_time_to_micro-depot_e2} ensures that the return time to the micro-depot is set only when there is an incoming arc. Lastly, constraint~\ref{constr:time_windows} is the time-window constraint for each customer.
            
            The domains of the decision variables are defined by \ref{constr: arc_decision_e1} - \ref{constr: domains}.

        \section{Algorithmic Component Analysis}\label{ablation}
            We conducted a statistical analysis to identify the effect of each algorithmic component on the solution quality. The results are shown in Table \ref{tab:Table6}. We computed results on large-scale instances with more than 100 customers. We started by using all the algorithmic components and computed a reference value, the average of the best solutions out of ten runs for the selected instances. Then, we removed every single component of our algorithm and computed the deviation, $\Delta \bar{\lambda}$, between the reference value and the value derived for each configuration. If removing a component decreased solution quality ($\Delta \bar{\lambda} > 0$), we kept the respective component in our algorithm.~Conversely, if $\Delta \bar{\lambda} \leq 0$, we excluded the component. We tested our algorithm in three blocks: basic algorithmic components, local search operators, and destroy operators. Starting with the basic algorithmic components, we identified the configuration that provided the best results, updated the reference value, and proceeded to the next block. This iterative approach ensured that the algorithm retained only the most effective components.
            
            \textbf{Basic algorithmic components:} Within this block, we tested all the basic algorithmic components containing the \textit{local search}, the \textit{simulated annealing acceptance} criterion, the \textit{adaptive learning} component, and the \textit{set cover} component. We found that removing any of the components decreased the solution quality. Therefore, we kept all the basic algorithmic components in our algorithm.
            
            \textbf{Local search operators:} Within this block, we tested all the local search operators: \textit{relocate intra-route, relocate inter-route, 2-opt, 2-opt$^{*}$} and \textit{exchange}. We found that removing the \textit{relocate intra-route} operator increased the solution quality. Thus, we removed the \textit{relocate intra-route} operator from our local search.
            
            \textbf{Destroy operators:} In this block, we tested all the destroy operators: \textit{Random removal, random string removal, furthest removal, micro-depot removal, partial micro-depot removal, partial micro-depot swap, FE route removal} and \textit{SE route removal}. We found that removing the \textit{partial micro-depot removal} operator increased the solution quality. Thus, we removed this operator from our configuration.
            
            \textbf{Repair operators:} In this block, we tested all the repair operators: \textit{Greedy insertion, regret insertion} and \textit{merge routes in first echelon}. We found that removing any of the operators decreased the solution quality. Therefore, we kept all the repair operators in our algorithm.

            \begin{landscape}
            \begin{table}[htbp]
              \centering
              \caption{Statistical analysis of algorithmic components and operators}
              \resizebox{\columnwidth}{!}{%
                \begin{tabular}{lcccccccccccccccccccc}
                \toprule
                \textbf{Basic algorithmic components} &       &       &       &       &       &       &       &       &       &       &       &       &       &       &       &       &       &       &       &  \\
                \midrule
                Local search &       & $\bullet$     & $\bullet$     & $\bullet$     & $\bullet$     & $\bullet$     & $\bullet$     & $\bullet$     & $\bullet$     & $\bullet$     & $\bullet$     & $\bullet$     & $\bullet$     & $\bullet$     & $\bullet$     & $\bullet$     & $\bullet$     & $\bullet$     & $\bullet$     & $\bullet$ \\
                Simulated annealing acceptance & $\bullet$     &       & $\bullet$     & $\bullet$     & $\bullet$     & $\bullet$     & $\bullet$     & $\bullet$     & $\bullet$     & $\bullet$     & $\bullet$     & $\bullet$     & $\bullet$     & $\bullet$     & $\bullet$     & $\bullet$     & $\bullet$     & $\bullet$     & $\bullet$     & $\bullet$ \\
                Adaptive learning & $\bullet$     & $\bullet$     &       & $\bullet$     & $\bullet$     & $\bullet$     & $\bullet$     & $\bullet$     & $\bullet$     & $\bullet$     & $\bullet$     & $\bullet$     & $\bullet$     & $\bullet$     & $\bullet$     & $\bullet$     & $\bullet$     & $\bullet$     & $\bullet$     & $\bullet$ \\
                Set cover & $\bullet$     & $\bullet$     & $\bullet$     &       & $\bullet$     & $\bullet$     & $\bullet$     & $\bullet$     & $\bullet$     & $\bullet$     & $\bullet$     & $\bullet$     & $\bullet$     & $\bullet$     & $\bullet$     & $\bullet$     & $\bullet$     & $\bullet$     & $\bullet$     & $\bullet$ \\
                \midrule
                \textbf{Local search operators} &       &       &       &       &       &       &       &       &       &       &       &       &       &       &       &       &       &       &       &  \\
                \midrule
                Relocate intra-route & $\bullet$     & $\bullet$     & $\bullet$     & $\bullet$     &       & $\bullet$     & $\bullet$     & $\bullet$     & $\bullet$     &       &       &       &       &       &       &       &       &       &       &  \\
                Relocate inter-route & $\bullet$     & $\bullet$     & $\bullet$     & $\bullet$     & $\bullet$     &       & $\bullet$     & $\bullet$     & $\bullet$     & $\bullet$     & $\bullet$     & $\bullet$     & $\bullet$     & $\bullet$     & $\bullet$     & $\bullet$     & $\bullet$     & $\bullet$     & $\bullet$     & $\bullet$ \\
                2-opt & $\bullet$     & $\bullet$     & $\bullet$     & $\bullet$     & $\bullet$     & $\bullet$     &       & $\bullet$     & $\bullet$     & $\bullet$     & $\bullet$     & $\bullet$     & $\bullet$     & $\bullet$     & $\bullet$     & $\bullet$     & $\bullet$     & $\bullet$     & $\bullet$     & $\bullet$ \\
                2-opt* & $\bullet$     & $\bullet$     & $\bullet$     & $\bullet$     & $\bullet$     & $\bullet$     & $\bullet$     &       & $\bullet$     & $\bullet$     & $\bullet$     & $\bullet$     & $\bullet$     & $\bullet$     & $\bullet$     & $\bullet$     & $\bullet$     & $\bullet$     & $\bullet$     & $\bullet$ \\
                Exchange & $\bullet$     & $\bullet$     & $\bullet$     & $\bullet$     & $\bullet$     & $\bullet$     & $\bullet$     & $\bullet$     &       & $\bullet$     & $\bullet$     & $\bullet$     & $\bullet$     & $\bullet$     & $\bullet$     & $\bullet$     & $\bullet$     & $\bullet$     & $\bullet$     & $\bullet$ \\
                \midrule
                \textbf{Destroy Operators} &       &       &       &       &       &       &       &       &       &       &       &       &       &       &       &       &       &       &       &  \\
                \midrule
                Random removal & $\bullet$     & $\bullet$     & $\bullet$     & $\bullet$     & $\bullet$     & $\bullet$     & $\bullet$     & $\bullet$     & $\bullet$     &       & $\bullet$     & $\bullet$     & $\bullet$     & $\bullet$     & $\bullet$     & $\bullet$     & $\bullet$     & $\bullet$     & $\bullet$     & $\bullet$ \\
                Random string removal & $\bullet$     & $\bullet$     & $\bullet$     & $\bullet$     & $\bullet$     & $\bullet$     & $\bullet$     & $\bullet$     & $\bullet$     & $\bullet$     &       & $\bullet$     & $\bullet$     & $\bullet$     & $\bullet$     & $\bullet$     & $\bullet$     & $\bullet$     & $\bullet$     & $\bullet$ \\
                Furthest removal & $\bullet$     & $\bullet$     & $\bullet$     & $\bullet$     & $\bullet$     & $\bullet$     & $\bullet$     & $\bullet$     & $\bullet$     & $\bullet$     & $\bullet$     &       & $\bullet$     & $\bullet$     & $\bullet$     & $\bullet$     & $\bullet$     & $\bullet$     & $\bullet$     & $\bullet$ \\
                Micro-depot removal & $\bullet$     & $\bullet$     & $\bullet$     & $\bullet$     & $\bullet$     & $\bullet$     & $\bullet$     & $\bullet$     & $\bullet$     & $\bullet$     & $\bullet$     & $\bullet$     &       & $\bullet$     & $\bullet$     & $\bullet$     & $\bullet$     & $\bullet$     & $\bullet$     & $\bullet$ \\
                Partial micro-depot removal & $\bullet$     & $\bullet$     & $\bullet$     & $\bullet$     & $\bullet$     & $\bullet$     & $\bullet$     & $\bullet$     & $\bullet$     & $\bullet$     & $\bullet$     & $\bullet$     & $\bullet$     &       & $\bullet$     & $\bullet$     & $\bullet$     &       &       &  \\
                Partial micro-depot swap & $\bullet$     & $\bullet$     & $\bullet$     & $\bullet$     & $\bullet$     & $\bullet$     & $\bullet$     & $\bullet$     & $\bullet$     & $\bullet$     & $\bullet$     & $\bullet$     & $\bullet$     & $\bullet$     &       & $\bullet$     & $\bullet$     & $\bullet$     & $\bullet$     & $\bullet$ \\
                FE route removal & $\bullet$     & $\bullet$     & $\bullet$     & $\bullet$     & $\bullet$     & $\bullet$     & $\bullet$     & $\bullet$     & $\bullet$     & $\bullet$     & $\bullet$     & $\bullet$     & $\bullet$     & $\bullet$     & $\bullet$     &       & $\bullet$     & $\bullet$     & $\bullet$     & $\bullet$ \\
                SE route removal & $\bullet$     & $\bullet$     & $\bullet$     & $\bullet$     & $\bullet$     & $\bullet$     & $\bullet$     & $\bullet$     & $\bullet$     & $\bullet$     & $\bullet$     & $\bullet$     & $\bullet$     & $\bullet$     & $\bullet$     & $\bullet$     &       & $\bullet$     & $\bullet$     & $\bullet$ \\
                \midrule
                \textbf{Repair Operators} &       &       &       &       &       &       &       &       &       &       &       &       &       &       &       &       &       &       &       &  \\
                \midrule
                Greedy insertion & $\bullet$     & $\bullet$     & $\bullet$     & $\bullet$     & $\bullet$     & $\bullet$     & $\bullet$     & $\bullet$     & $\bullet$     & $\bullet$     & $\bullet$     & $\bullet$     & $\bullet$     & $\bullet$     & $\bullet$     & $\bullet$     & $\bullet$     &       & $\bullet$     & $\bullet$ \\
                Regret insertion & $\bullet$     & $\bullet$     & $\bullet$     & $\bullet$     & $\bullet$     & $\bullet$     & $\bullet$     & $\bullet$     & $\bullet$     & $\bullet$     & $\bullet$     & $\bullet$     & $\bullet$     & $\bullet$     & $\bullet$     & $\bullet$     & $\bullet$     & $\bullet$     &       & $\bullet$ \\
                Merge routes in first echelon & $\bullet$     & $\bullet$     & $\bullet$     & $\bullet$     & $\bullet$     & $\bullet$     & $\bullet$     & $\bullet$     & $\bullet$     & $\bullet$     & $\bullet$     & $\bullet$     & $\bullet$     & $\bullet$     & $\bullet$     & $\bullet$     & $\bullet$     & $\bullet$     & $\bullet$     &  \\
                \midrule
                $\Delta\bar{\lambda}$    & 3.39  & 0.22  & 2.63  & 2.01  & \textbf{-0.15} & 0.35  & 0.19  & 0.12  & 0.68  & 0.34  & 0.20  & 0.24  & 1.00  & \textbf{-0.04} & 0.16  & 0.16  & 0.00 & 0.40 & 0.92 & 0.16\\
                \bottomrule
                \end{tabular}%
                }
                \begin{tablenotes} \footnotesize
                    \item The table shows the average deviation for the best results achieved out of 10 runs over selected 30 instances $\Delta\bar{\lambda}$ for the respective component analysis. 
                \end{tablenotes}
              \label{tab:Table6}%
            \end{table}%
            \end{landscape}
        
        \section{Capacitated Facility Location Problem (CFLP)}\label{cflp}
            The \gls{cflp} is defined as follows:
            
            \textbf{Sets}
            
                \begin{tabular}{lll}
                    $C$ & set of customers & $i \in C$ \\
                    $\mSetDepots{}$ & set of micro-depots & $j \in \mSetDepots{}$\\
                \end{tabular}\\
            
            \textbf{Parameter}
            
                \begin{tabular}{p{1cm}p{15cm}}
                    $c_{ij}$ & travel cost between micro-depot $j$ and customer $i$\\
                    $c_{0j}$ & travel cost between main depot and micro-depot $j$\\
                    $f_{j}$ & fixed opening cost of micro-depot $j$\\
                    $d_{i}$ & demand of customer $i$ \\
                    $q_{j}$ & capacity of micro-depot $j$ \\
                \end{tabular}\\
            
            \textbf{Decision Variable}
            
                \begin{tabular}{p{1.5cm}p{12.5cm}}
                    $x_{ij}$ & binary variable indicating whether customer $i$ is assigned to depot $j$\\
                    $y_{j}$ & binary variable indicating whether opening micro-depot $j$\\
                \end{tabular}\\
            
            \textbf{Model}
            
            \begin{equation}
                \min  \SumSet{i\in I} \SumSet{j \in \mSetDepots{}} c_{ij} x_{ij} + \SumSet{j \in \mSetDepots{}} (c_{0j} + f_{j}) y_{j}
            \end{equation}
            
            subject to
            
            Each customer must be assigned to exactly one micro-depot
            \begin{equation}
                \SumSet{j\in \mSetDepots{}} x_{ij} = 1  \quad\quad \forall i\in C
            \end{equation} 
            
            Customers can only be assigned to open depots
            \begin{equation}
                x_{ij} \leq y_{j}  \quad\quad \forall i\in C, j \in \mSetDepots{}
            \end{equation} 
            
            The demand of customers assigned to a micro-depot cannot exceed its capacity
            \begin{equation}
                \SumSet{i\in I} d_{i} x_{ij} \leq q_{j} y_{j}  \quad\quad \forall j \in \mSetDepots{}
            \end{equation} 
        
        \section{Speed-up Techniques}\label{speed-up}
            We apply preprocessing to speed up the computation time when applying repair operators or the local search. For fast evaluation, we use resource extension functions (REFs) as proposed by \cite{Irnich2007}. We also apply forward and backward calculations proposed by \cite{Vidal2014} as route-evaluation components. During repair and local search, these components help generate information about sub-sequences in a preprocessing step. We then evaluate new sequences created by concatenating multiple sub-sequences, using the precomputed information from these sub-sequences. We define the REFs to calculate the load, distance, duration, and time windows as in Equations~\eqref{eq:ref_load}-\eqref{eq:ref_tw}, which enables the route evaluations in $\mathcal{O}(1)$~\citep{Vidal2014}. 
            
            \begin{align}      
                Q(\sigma_1 \oplus \sigma_2) &= Q(\sigma_1) + Q(\sigma_2)  \label{eq:ref_load}\\
                D(\sigma_1 \oplus \sigma_2) &= D(\sigma_1) + d_{\sigma_1(|\sigma_1|)\sigma_2(1)} + D(\sigma_2) \label{eq:ref_distance} \\ 
                T(\sigma_1 \oplus \sigma_2) &= T(\sigma_1) + d_{\sigma_1(|\sigma_1|)\sigma_2(1)} + T(\sigma_2) \label{eq:ref_travel_time} \\
                E(\sigma_1 \oplus \sigma_2) &= \max\{E(\sigma_1) + d_{\sigma_1(|\sigma_1|)\sigma_2(1)} + T(\sigma_2), E(\sigma_2)\}  \label{eq:ref_et}\\
                L(\sigma_1 \oplus \sigma_2) &= \min\{L(\sigma_1), L(\sigma_2) - d_{\sigma_1(|\sigma_1|)\sigma_2(1)} - T(\sigma_1)\}  \label{eq:ref_lt}\\
                F(\sigma_1 \oplus \sigma_2) &\equiv F(\sigma_1) \land F(\sigma_2) \land (E(\sigma_1) + d_{\sigma_1(|\sigma_1|)\sigma_2(1)} \leq L(\sigma_2)) \label{eq:ref_tw}
            \end{align}
\end{appendices}
\end{document}